
\input amstex
\documentstyle{amsppt}
\magnification=1200
\NoRunningHeads
\baselineskip 18pt plus 2pt
\parskip 6pt
\def\and{\quad\hbox{\rm and}\quad}
\def\spec{\hbox{\rm spec}}

\def\rd{{\Bbb R}^d}
\def\swdb{{\Cal S}(\rd,{B})}

 \def\swddb{{\Cal S}'(\rd,{B})}

\def\ddwb{{\Cal D}'(\rd,{B})}

\def\hswb{H_s(\rd,{B})}
\def\hsw{H_s(\rd,{\Cal W})}
\def\hswib{H_s(\rd,{B})}

\def\hse{H_s({\Cal E})}
\def\crdww{C^\infty(\rd \times \rd, {\Cal L}_{{\Cal A}}({\Cal
W}, {\Cal W}))}
\def\srdww{{\Cal S}'(\rd \times \rd, {\Cal L}_{{\Cal A}}({\Cal
W}, {\Cal W}))}
\def\imw{I^m(\rd \times \rd, {\roman {Diag}}(\Bbb R^d), {\Cal
L}_{{\Cal A}}({\Cal W}))}

\def\crdwo{C^\infty_0(\rd ,{\Cal W})}
\def\crd{C^\infty(\rd )}

\def\ce{C^\infty({\Cal E} )}
\def\swm{ S_{{\Cal W}}^m}
\def\prdm{\Psi DO^m_{{\Cal W}}(\rd)}
\def\pbm{\Psi DO^m_{{\Cal B}}(M)}
\def\pbmm{\Psi DO^{-m}_{{\Cal B}}(M)}
\def\pbmi{\Psi DO^{-\infty}_{{\Cal B}}(M)}

\def\plm{\Psi DO^{m}_{{\Cal L}({\Cal E})}(M)}
\def\trn{{\roman {tr}}_{{\roman N}}}
\def\dimn{{\roman {dim}}_{{\roman N}}}
\def\ldet{{\roman {logdet}}}
\def\spec{{\roman {spec}}}
\def\ceg{C^\infty({\Cal E}_\Gamma)}
\def\cegm{C^\infty({\Cal E}|_{\Gamma^-})}
\def\cegp{C^\infty({\Cal E}|_{\Gamma^+})}
\def\ceoverg{C^\infty({\Cal E}|_{\Gamma})}
\def\crdwo{C^\infty_0(\rd ,{\Cal W})}
\def\crd{C^\infty(\rd )}

\def\ce{C^\infty({\Cal E} )}
\def\swm{ S_{{\Cal W}}^m}
\def\prdm{\Psi DO^m_{{\Cal W}}(\rd)}
\def\pbm{\Psi DO^m_{{\Cal B}}(M)}
\def\pbmm{\Psi DO^{-m}_{{\Cal B}}(M)}
\def\pbmi{\Psi DO^{-\infty}_{{\Cal B}}(M)}

\def\plm{\Psi DO^{m}_{{\Cal L}({\Cal E})}(M)}
\def\trn{{\roman {tr}}_{{\roman N}}}
\def\dimn{{\roman {dim}}_{{\roman N}}}
\def\ldet{{\roman {logdet}}}
\def\spec{{\roman {spec}}}
\def\ceg{C^\infty({\Cal E}_\Gamma)}
\def\cegm{C^\infty({\Cal E}|_{\Gamma^-})}
\def\cegp{C^\infty({\Cal E}|_{\Gamma^+})}
\def\ceoverg{C^\infty({\Cal E}|_{\Gamma})}
\NoBlackBoxes
\def\hang{\hangindent\parindent}
\def\item{\par\hang\textindent}

\def\textindent#1{\indent\llap{#1\enspace}\ignorespaces}
\hsize 6.5truein\vsize 8.9truein
\catcode`\@=11
\def\leftheadline{\rlap{\ }\hfill\iftrue\topmark\fi\hfill}
\def\rightheadline{\hfill\expandafter\iffalse\botmark\fi\hfill\llap{\ }}
\def\output@{\shipout\vbox{%
 \iffirstpage@  \global\firstpage@false
  \pagebody \logo@ \makefootline%
 \else \ifrunheads@ \makeheadline \pagebody \makefootline
       \else \pagebody \makefootline \fi
 \fi}%
 \advancepageno \ifnum\outputpenalty>-\@MM\else\dosupereject\fi}
\catcode`\@=\active
\def\qme{\Lambda^q(M;{\Cal E})}
\def\qh{{\Cal H}_q(M;{\Cal E})}
\def\qmme{\Lambda^{q-1}(M;{\Cal E})}
\def\qmpe{\Lambda^{q+1}(M;{\Cal E})}
\def\qmep{\Lambda^{+,q}_t(M;{\Cal E})}
\def\qmem{\Lambda^{-,q}_t(M;{\Cal E})}
\def\qht{{\Cal H}^q_t(M;{\Cal E})}
\def\tdq{\widetilde{ \Delta^+_q}(t)}
\def\ldetn{{\roman {logdet}}_{{\roman N}}}
\def\trn{{\roman {tr}}_{{\roman N}}}


\topmatter
\title Analytic and  Reidemeister torsion for representations
in finite type Hilbert modules
\endtitle
\thanks * Partially supported by NSF
\endthanks
\abstract
For a closed Riemannian manifold $(M,g)$ we extend the definition
of analytic and Reidemeister torsion associated to an orthogonal representation
of $\pi_1(M)$ on a
$\Cal{ A-}$Hilbert module
$\Cal W$ of finite type where $\Cal A$ is a finite von Neumann algebra.

If $(M,\Cal W)$ is of determinant class we prove, generalizing
the Cheeger-M\"uller theorem, that the analytic and Reidemeister
 torsion are equal. In particular, this proves the conjecture
that for closed Riemannian manifolds with positive Novikov-Shubin
invariants, the $L_2-$analytic and $L_2-$Reidemeister torsions are equal.
\endabstract

\author
D.Burghelea (Ohio State University)\\
L.Friedlander (University of Arizona)*\\
T. Kappeler (Ohio State University)*\\
P. McDonald (Ohio State University)\\
\endauthor

\endtopmatter

\document
\heading
0. Introduction
\endheading

The purpose of this paper is to prove the equality of $L_2-$analytic and
$L_2-$Reidemeister torsion. Both torsions are numerical invariants
defined for closed manifolds of determinant class, in particular
for closed manifolds with positive Novikov-Shubin invariants.
For these manifolds
their equality has been conjectured by
Carey, Mathai, Lott, L\"uck, Rothenberg and others
 (cf e.g. [LL, conjecture 9.7]).
It is implicit in [L\"u3](cf [BFK3]) that any closed
manifold whose fundamental group is residually finite is of determinant class.
\newline The interest
of the conjecture comes, among other issues, from the geometric
significance of the $L_2-$analytic torsion and the fact that sometimes the
$L_2-$Reidemeister torsion can be computed numerically in a very efficient
way. Indeed, if M is a closed hyperbolic manifold of dimension 3, the
$L_2-$analytic torsion coincides, up to a factor $-1/3\pi,$
with the hyperbolic volume.
\newline We establish the conjecture by proving a more general result: Given
a closed Riemannian manifold $(M,g),$ we extend the notion of analytic and
Reidemeister torsions to orthogonal representations of the fundamental
group $\pi_1(M)$ on a $\Cal A-$Hilbert module $\Cal W$ of finite type
 where $\Cal A$ is a
finite von Neumann algebra, and prove the equality of the two torsions
when $(M,\Cal W)$ is of determinant class.
We point out that in the case where $\Cal A$ is $\Bbb R,$ we obtain a new
proof of the well known result due to, independently, Cheeger [Ch] and
M\"uller [M\"u].
\newline From
the analytic point of view, the main difficulty comes from the fact that
the Laplacians, associated to such representations, may have continuous
spectrum. In particular, 0 might be in the essential spectrum.
\newline In order to formulate our results more precisely, we introduce
the following notations.
\newline Let $M$ be a closed smooth manifold.
A generalized triangulation of $M$ is a pair
$\tau =(h,g')$ with the following properties:
\newline (T1) $h:M\to \Bbb R$ is a smooth Morse function which is
self-indexing ($h(x)$ = index$(x)$ for any critical point $x \text{ of
} h);$
\newline (T2) $g'$ is a Riemannian metric so that $\text{ grad }_{g'}h$
satisfies the Morse-Smale condition (for any two critical points $x$
and $y$ of $h,$ the stable manifold $W^+_x$ and the unstable manifold
$W^-_y$, with respect to $\text{ grad }_{g'}h$, intersect transversely);
\newline (T3) in a neighborhood of any critical point of
$h$ one can introduce local coordinates such that,
with $q$ denoting the index of  this critical point,
$$
h(x)=q-(x_1^2+...+x_q^2)/2+(x_{q+1}^2+...+x_d^2)/2
$$
and the metric $g'$ is Euclidean in these coordinates.

The unstable manifolds $W^-_x$ provide a partition of $M$ into open cells
where $W^-_x$ is an open cell of dimension equal to the index of $x$. The name
``generalized triangulation'' for the pair $(h,g')$ is justified as a
generalized
triangulation can be viewed as a generalization of a simplicial
triangulation.\footnote {Given a smooth simplicial triangulation
$\tau_{sim}$, one can
construct a generalized triangulation $\tau =(h,g')$ so that the unstable
manifolds $W^-_x$ corresponding to $\text{ grad }_{g'}h,$ with $x$ a
critical point of $h,$ are the open simplexes of $\tau_{sim}$
(cf.\cite {Po}).}

Let $(M,g)$ be a closed Riemannian manifold with infinite fundamental group
$\Gamma=\pi_1(M)$
and let $\tau=(h,g')$ be a generalized triangulation.
 Note that $\Gamma$ is  countable.  Let $p:\tilde M\to M$ be
the universal covering of $M$
and denote by $\tilde g$ and
$\tilde{\tau}=(\tilde h, \tilde g')$ the lifts of $g$ and $\tau$ on $\tilde M.$
The Laplace operator $\Delta_q$ acting on compactly supported, smooth
q-forms on $\tilde M$
is essentially selfadjoint. Its closure,
also denoted by $\Delta_q$, is therefore selfadjoint; it is defined on
a dense subspace of the $L_2$-completion  of the space
of smooth forms
with  compact support with respect to the scalar product induced by the metric
$\tilde g.$ Observe that $\Delta_q$ is $\Gamma-$equivariant and nonnegative.
We can therefore define the spectral projectors $Q_q(\lambda)$ of $\Delta_q$
corresponding to the interval $(-\infty ,\lambda \rbrack$. They are
$\Gamma$-equivariant
and admit a $\Gamma$-trace which we denote by $N_q(\lambda)$.
\newline Let $\Cal C^q(\tilde\tau):=l^2({\roman {Cr}}_q(\tilde h))$
where  $l^2({\roman {Cr}}_q(\tilde
 h))$ denotes the Hilbert space
of $l_2$-summable, real-valued sequences indexed by the countable
set ${\roman {Cr}}_q(\tilde h)$ of critical points of $\tilde h$ of index q.
The left action of $\Gamma$ on $ {\roman {Cr}}_q(\tilde h)$ makes
$l^2({\roman {Cr}}_q(\tilde h))$ the underlying
Hilbert space of an orthogonal  $\Gamma-$representation.
The intersections of the stable and the unstable manifolds of
${\roman {grad}}_{\tilde g'}\tilde h$ induce a bounded, $\Gamma$-equivariant,
linear map
$$\delta_q: \Cal C^q(\tilde\tau) \to \Cal C^{q+1}(\tilde\tau).$$
Let $\delta_q^*$ be the adjoint of $\delta_q$ and introduce
 $$\Delta_q^{{\roman {comb}}}:=
\delta^*_q\cdot\delta_q+ \delta_{q-1}\cdot\delta^*_{q-1}.$$
Observe that $\Delta^{{\roman {comb}}}_q$ is a $\Gamma$-equivariant, bounded,
nonnegative, selfadjoint operator on $C^q(\tilde\tau)$.
We can therefore define the spectral projectors
$Q^{{\roman {comb}}}_q(\lambda)$ of $\Delta_q^{{\roman {comb}}}$ corresponding
to the interval $(-\infty ,\lambda \rbrack$. These projectors
are $\Gamma$-equivariant and thus admit a $\Gamma$-trace, which we denote by
$N^{{\roman {comb}}}_q(\lambda)$ (cf section 1).

We say that
\newline (1) $(M,g)$ is of $a-determinant$ class if $-\infty <\int _{0^+}^1
(\log\lambda) dN_q(\lambda)$ for all $q.$
\newline (2) $(M,\tau)$ is of $c-determinant$ class if
$-\infty < \int_{0^+}^1(\log\lambda) dN^{{\roman {comb}}}_q(\lambda)$ for
all $q.$

Here $\int_{0^+}^1$ denotes the Stieltjes integral on the half open
interval $(0,1].$
The following result can be derived from work of Gromov-Shubin $[GS]$
(cf also $[Ef]$).
\proclaim{\bf Proposition {1} }($[Ef]$, $[GS])$
Let $M$ be a closed manifold equipped with a Riemannian metric $g$ and a
generalized triangulation $\tau$. Then:
\newline (1) $(M,g)$ is of $a-determinant$ class iff $(M,\tau)$ is of
$c-determinant$ class.
\newline (2) Let $(M',\tau ')$ be another manifold with generalized
triangulation $\tau '$. If $M$ and $M'$ are homotopy equivalent, then
$(M,\tau)$ is of $c-determiant$ class iff $(M',\tau ')$ is.
\endproclaim

\proclaim{\bf Definition}
A compact manifold M is of determinant class if for some generalized
triangulation $\tau$ (and then for any), $(M,\tau)$ is of $c-determinant$
class.
\endproclaim
If $M$ is of determinant class the logarithm of the $\zeta$-regularized
determinant, ${\roman {logdet}}_N \Delta_q,$ is a finite
real number for all $q$ and
one can introduce the $L_2$- analytic
torsion $T_{an}$ :
$$\log T_{an}:=\frac{1}{2} \sum_q (-1)^{q+1}q {\roman {logdet}}_N  \Delta_q.$$

Similarly, if $M$ is of determinant class,
${\roman {logdet}}_N \Delta_q^{{\roman {comb}}}$
is a finite real number for all $q$ and one can  define the combinatorial
torsion:
$$\log T_{{\roman {comb}}}:=\frac{1}{2} \sum_q (-1)^{q+1}q
{\roman {logdet}}_N  \Delta^{{\roman {comb}}}_q.$$
To define the $L_2$-Reidemeister torsion, $T_{{\roman {Re}}},$ it remains to
introduce its metric part. Notice that ${\roman {Null}} (\Delta_q)$ consists
of smooth forms and that integration of smooth $q$-
forms over a smooth q-chain  induces (cf a theorem by Dodziuk $[Do]$)
an isomorphism $\theta^{-1}_q: {\roman {Null}} (\Delta_q)\to \overline
H^q(\Cal C^*(\tilde\tau),
\delta_*)=  {\roman {Null}} (\Delta^{{\roman {comb}}}_q)$ of
$\Cal N(\Gamma)$-Hilbert modules where $\Cal N(\Gamma)$ is the
von Neumann algebra associated to $\Gamma.$ Define
$$\log {\roman {Vol}}_N(\theta_q):=\frac{1}{2} {\roman {logdet}}_N
(\theta_q^*\theta_q)$$
where we used that ${\roman {det}}_N(\theta_q^*\theta_q) > 0$ as
$(\theta_q^*\theta_q)$ is a selfadjoint, positive,
bounded, $\Gamma$-equivariant operator on the $\Gamma$- Hilbert space
${\roman {Null}} (\Delta^{{\roman {comb}}}_q)$ whose spectrum is bounded away
 from $0$. As a consequence (cf section 1)
${\roman {logdet}}_N (\theta_q^*\theta_q)$ is a well defined real number and
one may introduce
$$\log T_{{\roman {met}}}:=\frac{1}{2} \sum _q(-1)^{q} {\roman {logdet}}_N
(\theta_q^*\theta_q).$$
Combining the above definitions we define the $L_2$-Reidemeister torsion
$T_{{\roman {Re}}}$
$$\log T_{{\roman {Re}}}=\log T_{{\roman {comb}}}+\log T_{{\roman {met}}}.$$
The concepts of $L_2$ -analytic and $L_2$-Reidemeister torsion were considered
by Novikov-Shubin in 1986 $[NS1]$ (cf.also later work by Lott $[Lo]$,
L\"uck-Rothenberg $[LR]$,Mathai $[Ma]$ and Carey-Mathai $[CM]$).
The main objective of this paper is to prove the following

\proclaim{\bf Theorem {1} }
Let $M$ be a closed manifold of determinant class of odd dimension d.
Then, for any Riemannian metric $g$ and for any
generalized triangulation $\tau,$ both $T_{{\roman {an}}}$ and
$T_{{\roman {Re}}}$ are positive real numbers and
$$T_{{\roman {an}}}=T_{{\roman {Re}}}.$$
\endproclaim

We point out that the corresponding result for manifolds of even dimension
is trivial as both torsions are equal to 1.
\newline Before explaining the ideas of  the proof
 of Theorem 1 let us make a few remarks:
\newline1)Lott-L\"uck  have conjectured (cf [LL], conjecture 9.2) that all
compact manifolds have positive Novikov-Shubin invariants and,therefore, are of
determinant class. The conjecture has been verified for many compact manifolds
and
in particular for all compact manifolds whose fundamental group is free or free
abelian.
\newline 2)By assigning to each compact Riemannan manifold $(M,g)$ with $M$ of
determinant class the $L_2$-(analytic=Reidemeister) torsion if the fundamental
group  $\pi_1(M)$ is infinite  and the usual (analytic=Reidemeister) torsion if
$\pi_1(M)$ is finite one obtains a numerical invariant $T(M,g)$, which
satisfies the product formula
$$
\log T(M_1\times M_2,g_1\times g_2)=\chi (M_2) \log T(M_1,g_1) +\chi (M_1)
\log T(M_2,g_2),
$$
and for any $n$-sheeted covering  $(\tilde M,\tilde g)$ of $(M,g)$  satisfies
$$ \log T(\tilde M,\tilde g)=n \log T(M,g).$$  Here $\chi (M)$ denotes the
Euler-Poincar\'e characteristic of $M.$
For compact manifolds with trivial $L_2$ Betti numbers, in particular for
manifolds of the homotopy type of a mapping torus, this invariant is
independent of the metric and is in fact  a homotopy type invariant.
\newline This invariant was calculated for a large class of $3$-dimensional
manifolds ; its logarithm is zero for Seifert manifolds (cf [LR] ) and
$(-1/{3\pi}) Vol( M,g)$ for
a hyperbolic manifold $(M,g),$ (cf[Lo]). The calculation in [LR] was done for
the Reidemeister torsion and
in [Lo] for the analytic torsion.
\newline 3) Recently W. L\"uck ([L\"u2]) found an algorithm
to calculate the $L_2$-Reidemeister torsion  of a $3$-dimensional hyperbolic
manifold in terms of a balanced presentation of its fundamental group. By
Theorem 1 and by Remark 2) above the algorithm
also calculates the hyperbolic volume.

Rather than viewing the above theorem as an $L_2-$version, we derive Theorem 1
as a particular case of a generalization of the Cheeger-M\"uller
theorem. This generalization concerns the extension of the analytic and
 Reidemeister torsion associated to a closed Riemannian manifold and
a finite dimensional orthogonal representation of $\Gamma$ to orthogonal
$\Cal A-$representation of $\Gamma$
which are of finite von Neumann dimension where $\Cal A$
is a finite von Neumann algebra (cf $[Si]$ for a similar approach
in connection with an $L_2-$ index theorem).
A representation of this type is called an $(\Cal A,\Gamma^{op})-$Hilbert
module of finite type.

In order to formulate this generalization
we must introduce (cf section 2)
a calculus of elliptic pseudodifferential $\Cal A-$operators
acting on sections of a bundle of $\Cal A-$Hilbert modules of finite
type over a compact manifold (typically the spectrum of such an operator
is no longer discrete) and develop
a theory of regularized determinants for (nonnegative)
elliptic pseudodifferential $\Cal A$-operators of positive order.
Both ingredients are presented in section 2.
The new difficulties in this theory come from the fact that $0$ might be
in the essential spectrum of an elliptic pseudodifferential $\Cal A$-operator.

In order to prove this generalization we need a Mayer-Vietoris type
formula and an asymptotic expansion for the logarithm of the determinant
of elliptic operators with parameter; i.e. the extension of the
results of $[BFK2]$ to this calculus (cf section 3 and [Lee]).

Let us now describe the generalization of the above Theorem in more detail.
Assume that $A$ is an elliptic operator in the new calculus.
For an angle $\theta$ and $\epsilon > 0$ introduce the solid angle
$$V_{\theta,\epsilon} := \{ z \in \Bbb C : |z| < \epsilon \} \cup
\{ z \in \Bbb C \setminus 0 : arg(z) \in (\theta - \epsilon , \theta
+ \epsilon) \}.$$

\proclaim{\bf Definition}
(1) $\theta$ is an Agmon angle for $A$, if there exists $\epsilon > 0$ so that
 $${\roman {spec}}(A) \cap V_{\theta,\epsilon} = \emptyset.$$
(2) $\theta$ is a principal angle for $A$ if there exists $\epsilon > 0$
so that $${\roman {spec}}(\sigma _A(x,\xi)) \cap V_{\theta,\epsilon}=
\emptyset$$
for all $(x,\xi)\in S_x^*M$ where $S^*M$ denotes the cosphere bundle and
$\sigma_A(x,\xi)$ is the principal symbol of $A$.
\endproclaim

It is well know that (1) implies (2) in the above definition but the converse
statement is not true. If, in addition, $A$ is of order $m > 0$ and admits an
Agmon angle, $\theta,$ one can define the regularized determinant,
${\roman {det}}_{\theta,N} A \in \Bbb C$. In the sequel, $\theta$ will always
be chosen to be $\pi$ and we will drop the subscript $\pi$ in
${\roman {det}}_{\pi,N}$.  If $A$ is of order $m > 0$, nonnegative and if
$0 \in {\roman {spec}}(A)$ then the ellipticity of $A$ implies that the
nullspace, ${\roman {Null}}(A),$ is an $\Cal A$-Hilbert module of finite von
Neumann dimension, ${\roman {dim}}_N {\roman {Null}}(A).$ Consider the
1-parameter family, $A + \lambda,$ $\lambda$ being the spectral
parameter. For $\lambda > 0$, introduce
the function ${\roman {logdet}}_N (A+\lambda) - {\roman {dim}}_N
{\roman {Null}}(A)\log \lambda.$ We can view this function as an element
in the vector space $\bold D$ consisting of
equivalence classes $\lbrack f \rbrack$ of real analytic
functions $f:(0,\infty)\to\Bbb R$ with $f \sim g$ iff
$\lim_{\lambda\to 0}(f(\lambda)-g(\lambda))=0$. The elements of $\bold D$
represented by the constant
functions form a subspace of $\bold D$ which can be identified with
$\Bbb R$, the space of real numbers.

Given a closed Riemannian manifold $(M,g)$, an arbitrary $(\Cal
A,\Gamma^{op})-$ Hilbert
module of finite type, $\Cal W,$
and a generalized triangulation $\tau$
we define (cf section 4) $\log T_{{\roman {an}}}(M,g,\Cal W)$ and
$\log T_{{\roman {Re}}}(M,g,\tau,\Cal W)$ as elements of $\bold D$.
As above we consider the analytic resp. combinatorial Laplacians associated
to $(M,g,\Cal W)$ resp. $(M,\tau,\Cal W)$ and introduce
the notion of a triple $(M, g,\Cal W)$ resp.
$(M,\tau,\Cal W),$ of $a-determinant, $ resp. of $c-determinant$ class.
Proposition 1 can be generalized to say that these two notions are
equivalent and homotopy invariant (Proposition 2, section 5).
This allows us to introduce the notion of a pair $(M,\Cal W)$ to be of
determinant class.
We point out that for
$\Cal {A} = \Bbb {R}$ any pair $(M, \Cal W)$ is of determinant class.
\proclaim{\bf Theorem 2}
Let $M$ be a closed manifold of odd dimension d and $\Cal W$ an $(\Cal A,
\Gamma^{op})$-Hilbert module of finite type. If the pair $(M,\Cal W)$ is of
determinant class then, for any Riemannian metric
$g$ and any generalized triangulation
$\tau$ of M,
$\log T_{{\roman {an}}}( M,g,\Cal W)$ and
\newline $\log T_{{\roman {Re}}}(M,\tau,g,\Cal W)$ are  both
finite, real numbers and
$$\log T_{{\roman {an}}} = \log T_{{\roman {Re}}}.$$
\endproclaim

Let us make a few comments concerning Theorem 2. First note that if $M$
is of even dimension, then Theorem 2 is also true. In this case, however,
it is an immediate consequence of the Poincar\'e duality,which implies that
both torsions are equal to $1$.

If $\Cal {A} = \Bbb R$, the
$(\Bbb {R},\Gamma^{op})$-Hilbert modules of finite type are precisely the
orthogonal $\Gamma$-representations
and Theorem 2  reduces to the Cheeger-M\"uller Theorem ([Ch],[M\"u]) and, when
specialized to this situation, we thus obtain a new proof of their theorem.

It suffices to prove Theorem 2 for $\Cal W$ a free $\Cal A$-module. This
follows easily from the following three facts:
\newline 1): If $\Cal {A} = \Bbb R$ then $\Cal W$ is $\Bbb R$ free and the
pair $(M,\Cal W)$ is of determinant class.
\newline 2): If $\Cal W$ is $\Gamma^{op}$-trivial then $(M,\Cal W)$ is of
determinant class,
$\log T_{{\roman {an}}}(M,g,\Cal W)= \text {dim}_N \Cal W \cdot
\log T_{{\roman {an}}}(M,g)$ and $\log T_{{\roman {Re}}}(M,\tau,g,
\Cal W)= \text {dim}_N \Cal W\cdot\log T_{{\roman {Re}}}(M,\tau,g),$
where $\log T_{{\roman {an}}}(M,g)$ and
$\log T_{{\roman {Re}}}(M,\tau,g)$ denote the analytic and Reidemeister
torsions of the Riemannian manifold $(M,g)$ and the trivial
1-dimensional representation of $\pi_1(M)$.
\newline 3): Suppose that $\Cal W_1$ and $\Cal W_2$ are two
$(\Cal A,\Gamma^{op})$-Hilbert modules  of finite type and
$\Cal W=\Cal W_1\oplus\Cal W_2$. If $(M,\Cal W)$ and $(M,\Cal W_2)$ are
of determinant class then $(M,\Cal W_1)$ is of determinant class and the
analytic resp. Reidemeister torsion of  $(M,g,\Cal W)$
is  the product of the analytic resp. Reidemeister torsion of $(M,g,\Cal W_1)$
with the analytic resp. Reidemeister torsion $(M,g,\Cal W_2).$

In the case when $\Cal A=\Cal{N}(\Gamma)$, the von Neumann algebra of bounded
operators associated to $\Gamma,$ acting on $\Cal {W} = l_2(\Gamma)$, and
$\Cal W$ is viewed as a $(\Cal {N}(\Gamma), \Gamma^{op})-$Hilbert module of
finite type, then Theorem 2 reduces to Theorem 1.

Theorem 2 is derived from Corollary C (section 6.2), a relative version of
Theorem 2, using product formulas for the analytic torsion and the
Reidemeister torsion (section 4) and the metric
anomaly (Lemma 6.10). To state Corollary C let $M$ and $M'$ be
two closed manifolds of the same dimension with fundamental groups isomorphic
to $\Gamma$, and assume that they are
equipped with generalized triangulations $\tau=(g,h)$ and $\tau'=(g',h')$ such
that the functions $h$ and $h'$ have the same number of critical points for
each index. Then, for an arbitrary  $(\Cal A,\Gamma^{op})$- Hilbert module
of finite type $\Cal W$, with $(M,\Cal W)$ and $(M',\Cal W')$ of determinant
class
$$\log T_{{\roman {an}}} -\log T'_{{\roman {an}}} = \log T_{{\roman {Re}}}
- \log T'_{{\roman {Re}}}.
$$
An important ingredient in the proof of Corollary C is the Witten deformation
of the de Rham complex associated with a generalized triangulation
$\tau=(g,h)$.
In section 5 we extend the analysis of Helffer-Sj\"ostrand [HS1] of the
Witten complex
to the analogous complex constructed for differential forms on $M$ with
coefficients in $\Cal W$. Although the results and the estimates in the case
where $\Cal A$ is an arbitray finite von Neumann
algebra are similar to those obtained by Helffer-Sj\"ostrand
in the case $\Cal A=\Bbb R$,  additional arguments are necessary
since the spectrum of the Laplacians
$\Delta_q(t)$ is typically not discrete.

The Witten deformation permits us to
define smooth functions $\log T_{{\roman {an}}}(h,t)$,
$\log T_{{\roman {sm}}}(h,t)$ and $\log T_{{\roman {la}}}(h,t)$ with
$\log T_{{\roman {an}}}(h,0)=\log T_{{\roman {an}}}$ where
$\log T_{{\roman {an}}}(h,t)=\log T_{{\roman {sm}}}(h,t)+
\log T_{{\roman {la}}}(h,t) $ is a decomposition of
$\log T_{{\roman {an}}}(h,t)$ into a part $\log T_{{\roman {sm}}}(h,t)$ which
corresponds to the small spectrum of the Laplacians $\Delta_q(t)$ and
a complimentary part $\log T_{{\roman {la}}}(h,t).$
The results presented in sections 2 and 3 lead to the conclusion
that these three
functions have asymptotic expansion when $t\to\infty.$The free term of
such an expansion refers to the $0$'th order coefficient of
the expansion as $t\to \infty.$

The results presented in sections 3 and 5
permit us to show that the free term of
$$\log T_{{\roman {an}}}(h,t)-\log T_{{\roman {sm}}}(h,t)-
\log T_{{\roman {an}}}(h',t)-\log T_{{\roman {sm}}}(h',t)$$
is equal to
$$\log T_{{\roman {an}}} -\log T_{{\roman {Re}}}-\log T'_{{\roman {an}}}
- \log T'_{{\roman {Re}}}.$$
Using the Mayer Vietoris type formula (cf section 3) we show that the free
term of
$\log T_{{\roman {la}}}(h,t)-\log T_{{\roman {la}}}(h',t)$ is equal to zero
and we can therefore conclude
Corollary C.

The paper is organized as follows:

In section 1 we recall, for the convenience of the reader,
the concepts of a finite von Neumann algebra $\Cal A$,
an $\Cal A-$Hilbert module of finite type, a finite (von Neumann) dimensional
representation of a group, determinants in the von Neumann sense
and the torsion of a finite complex of $\Cal A-$Hilbert
modules of finite type. This section is entirely expository.

In section 2 and 3 we describe the theory of pseudodifferential
operators acting on sections of a given bundle $\Cal E \to M$ of $\Cal
A-$Hilbert modules of finite type. In particular, we extend Seeley's result
on  zeta-functions
for elliptic pseudodifferential operators and the corresponding
regularized determinants, as well as the results of $[BFK2]$, to the
extent needed in this paper, for this new class of operators. The calculus
of such operators is not new, but we failed to find a reference suited for our
needs (cf e.g. [FM],[Le] and [Lu] for related work).

In section 4 we define the analytic torsion and Reidemeister
torsion and we prove a product formula for each of them. These product
formulas are slight generalizations of the product formulas presented
in $[L]$ and $[CM],$ but for the convenience of the reader we include
the proofs.

In section 5, we discuss the Witten deformation of the deRham
complex of $M$ with coefficients in a representation
of finite von Neumann dimension. Moreover using the Witten deformation we
prove Proposition 2.

In section 6  we present the proof of  Theorem 2.

One can generalize the analytic and Reidemeister torsions associated to
$(M,g,\Cal W)$ to include additional data, for example a finite dimensional
hermitian vector bundle on $M$ equipped with a flat connection . By the
same methods as presented in this paper one can prove  a result which compares
these  two generalized torsions.
In the case $\Cal A=\Bbb R$ such type of result was first established in
([BZ]), compare also ([BFK1].

Using the same arguments as in ([L\"u1]), one can  extend Theorem 2 to compact
manifolds with boundary. Both  extensions are useful for the calculations of
the
$L_2$ torsions; together with some applications  they will be presented in a
forthcoming paper.

The authors wish to thank W. L\"uck for bringing the conjecture about the
equality of analytic and Reidemeister torsion to their attention,  J. Lott,
W. L\"uck and M.Shubin for informations on related  work and  D.Gong and
M.Shubin for comments on  a preliminary version of this paper.

\vfill \eject
\heading
        1. Linear algebra in the von Neumann sense
\endheading

In this section we collect for the convenience of the reader
a number of results concerning linear algebra in the von Neumann sense
(cf e.g. [CM],[Co],[Di],[GS],[LR] for reference).

\subheading{1.1 $\Cal A-$Hilbert modules}

\proclaim{\bf Definition 1.1}
A finite von Neumann algebra $\Cal A$ is a unital $\Bbb R$-algebra with a
 $*$-operation and a trace
${\roman {tr}}_N:\Cal A\to \Bbb R$ which satisfies the following properties:
\newline (T1) $\langle .,.\rangle : \Cal A\times \Cal A\to \Bbb R$, defined
by $\langle a,b\rangle ={\roman {tr}}_N (a b^*)$, is a scalar product and the
completion $\Cal A_2$ of $\Cal A$ with respect to this scalar product is a
 separable Hilbert space.
\newline (T2) $\Cal A$ is separable and weakly closed, when viewed as a
subalgebra of $\Cal{L}(\Cal A_2):=\Cal{L}(\Cal A_2,\Cal A_2)$, the space of
linear, bounded operators on $\Cal A_2$, where elements of $\Cal A$ are
identified with the corresponding left translations in $\Cal {L}(\Cal A_2)$
(a sequence $\{a_n\}_{n \ge 1}$ in $\Cal A$ converges weakly to
$a \in \Cal A_2$ if $lim_{n \to \infty} \langle a_nx,y\rangle
=\langle ax,y\rangle $ for all $x,y \in \Cal A_2).$
\newline (T3) The trace is normal, i.e for any monotone increasing net,
$(a_i)_{i\in I},$ such that $a_i\geq 0$ and
$a=\sup_{i\in I}a_i$ exists in $\Cal A$, one has ${\roman {tr}}_Na=
\sup_{i\in I} {\roman {tr}}_N a_i.$
Here $a_i\geq 0$
means that $a_i=a^*_i$ and $\langle a_ix,x\rangle \geq 0$ for all $x\in \Cal
A.$
\endproclaim

In the sequel, $\Cal A$ will always denote a finite von Neumann algebra.
Introduce the opposite algebra $\Cal A^{op}$ of $\Cal A$, where $\Cal A^{op}$
has the same underlying
vector space, $|\Cal A^{op}|=|\Cal A|$, $*$-operation, trace  and unit element
as $\Cal A$, but the
multiplication $"\cdot_{op}"$ of the  elements $a,b\in |\Cal A^{op}|$ is
defined by $a\cdot_{op}b=b\cdot a.$
Note that $\Cal A^{op}$ is a finite von Neumann algebra as well.
The right translation by elements of $\Cal A$ induces an embedding
$r:\Cal A^{op}\to \Cal L(\Cal A_2)$
which identifies $\Cal A^{op}$ with the subalgebra $\Cal L_{\Cal
A}(\Cal A_2)\subset\Cal L(\Cal A_2)$ of bounded $\Cal A $ -linear maps
(with respect to the $\Cal A$-module structure of $\Cal A_2$ induced
by left multiplication). Therefore we can introduce a trace on $\Cal
L_{\Cal A}(\Cal A_2)$, also denoted by ${\roman {tr}}_N$, defined for
$f \in \Cal L_{\Cal A}(\Cal A_2)$ by $${\roman {tr}}_N(f):={\roman
{tr}}_N(r^{-1}(f)).$$

\proclaim{\bf Definition 1.2}
(1) $\Cal W$ is an $\Cal A-$Hilbert module if
\newline (HM1) $\Cal W$ is a Hilbert space with inner product denoted
by $\langle .,.\rangle .$
\newline (HM2) $\Cal W$ is a $left$ $\Cal A-$module so that $\langle
a^*v,w \rangle = \langle v,aw \rangle ~~ (a \in \Cal{A}; v,w \in \Cal
W).$
\newline (HM3) $\Cal W$ is isometric to a closed submodule of $\Cal A_2 \otimes
V$ where $V$ is a separable Hilbert space and the tensor product
$\Cal A_2 \otimes V$ is taken in the category of Hilbert spaces.
\newline(2) $\Cal W$ is an $\Cal A-$Hilbert module of finite type if $\Cal W$
is an $\Cal A-$Hilbert module and
\newline (HM4) $\Cal W$ is isometric to a closed submodule $\Cal A_2 \otimes V$
where $V$ is a finite dimensional vector space.
\newline (3) A morphism $f:\Cal W_1\to \Cal W_2$ between $\Cal A$-Hilbert
modules
of finite type, $\Cal W_1$ and $\Cal W_2$, is a bounded, $\Cal A-$linear
 operator; f is an isomorphism if it is bijective
 and both $f$ and $f^{-1}$ are morphisms.
\endproclaim

Let $\Cal W$ be an $\Cal A$-Hilbert module and $v$ an element in $\Cal W$. The
map $i_v:\Cal A_2\to \Cal W,$ defined by $i_v(a)=av,$ extends to an
$\Cal A$-linear bounded map $\Cal A_2\to \Cal W.$

\proclaim{\bf Definition 1.3}
A  collection $\{e_1,...,e_l\},l\leq\infty,$ of elements of $\Cal W$ is called
a base  of $\Cal W$
if $$i:\oplus_{1\leq\nu\leq l}(\Cal A_2)_\nu\to \Cal W \tag {1.1}$$
 is an isomorphism where each
$(\Cal A_2)_\nu$ is a copy of $\Cal A_2$ and $i=\sum_{1\leq\nu\leq l}
 i_{e_\nu}.$ The base is called orthonormal if, in addition, $i$
 is an isometry. A Hilbert module is free if it
has a base.
\endproclaim

If $l=\infty$ the direct sum in (1.1) is a Hilbert direct sum in the category
of Hilbert spaces. Note that  $\{e_1,...,e_l\}$ is an orthonormal basis iff
the closed invariant subspace spanned  by
$e_1,...,e_l$ is  $\Cal W,$ and for any $i,j$ and  $a,b \in \Cal A,$
$\langle ae_i , be_j \rangle = \langle a, b \rangle \delta_{ij}.$ Further if
$\{e_1,...,e_l\}$ is a base
then $\{f_1,...,f_l\}$ with $f_\nu=i(i^*i)^ {-\frac {1}{2}}
(0,...,0,1,0,...,0)$
 is an orthonormal base of $\Cal W.$

\proclaim{\bf Proposition 1.4}
\newline (1) Any $\Cal A$-Hilbert module $\Cal W$ can be decomposed as
$\Cal W=\oplus_{1\leq \nu\leq l}\Cal W_\nu$
with $l\leq \infty $ and $\Cal W_\nu$ a closed invariant subspace of
 $\Cal A_2.$ Further
$\Cal W$ is of
finite type iff $l <\infty.$
\newline (2) If $\Cal W'$ is a closed invariant subspace of $\Cal W$ then
$\Cal W' \simeq \oplus_{1\leq \nu\leq l}\Cal W'_\nu$ where
$\Cal W'_\nu$ is a closed invariant subspace of $\Cal W_\nu.$
\endproclaim

Proof: Note that (1) follows from (2). To prove (2), denote by $\pi_\nu$
the projection of $\Cal W$ on
$\Cal W_\nu$ and consider the filtration of $\Cal W'$,
$$\Cal W'(0):=\Cal W'\supset\Cal W'(1)\supset \Cal W'(2)\supset....$$ where
$\Cal W'(\nu)$ is defined
inductively $\Cal W'(\nu)= \text {Null}(\pi_\nu)\cap \Cal W'(\nu-1)$. It
follows that
$\Cal W'\simeq\oplus_{1\leq \nu < l} \Cal W'(\nu)/\Cal W'(\nu+1)\subset
 \Cal W_\nu. \diamondsuit$

Let $\Cal W$ be an $\Cal A-$Hilbert module of finite type. The algebra
\newline $\Cal L_{\Cal A}(\Cal W):=\Cal L_{\Cal A}(\Cal W,\Cal W)$
of bounded      $\Cal A$-linear operators on $\Cal W$ is a finite
von Neumann algebra, whose trace is defined in the following way.
First assume that the module $\Cal W$ is free. Choose a
basis $\{e_1,...,e_l\}.$ With respect to this basis an operator
$A\in\Cal L_{\Cal A}(\Cal W)$ has a matrix representation
$(a_{ij})_{1\leq i,j\leq l}$,
$i,j=1,...,l$, with entries $a_{ij}$ in
$\Cal L_{\Cal A}(\Cal A_2)=\Cal A^{op}$. Define
${\roman {tr}}_N(A)=\sum_{i=1}^l a_{ii}.$ One shows that
${\roman {tr}}_N(A)$ is independent of the
chosen basis and therefore well defined.
In the general case $\Cal W$ is a closed invariant subspace of a free
$\Cal A$-
Hilbert module $\Cal V$ of finite type . We write $\Cal V=\Cal W\oplus\Cal
W^{\perp}$ and consider
$\tilde A=A\oplus 0\in \Cal L_{\Cal A}(\Cal V,\Cal V).$ Define
${\roman {tr}}_N(A)={\roman {tr}}_N(\tilde A).$ One shows that
${\roman {tr}}_N(A)$ is independent of the choice of $\Cal
V.~~~~\diamondsuit$

For an $\Cal A$-Hilbert module $\Cal W$ of finite type one defines
the dimension
${\roman {dim}}_N(\Cal W)$  in the von Neumann sense by ${\roman
{dim}}_N\Cal W:= {\roman {tr}}_N{\roman {Id}}_{\Cal W}.$
If $\Cal W$ is not of finite type one sets
 ${\roman {dim}}_N\Cal W:=\sup\{{\roman {dim}}_N\Cal W';  \Cal W'
{\text{ closed submodule of finite type}} \}.$
 The von Neumann dimension is
always a nonnegative real number or $+\infty.$

{\bf Remark 1.5} If $\Cal W_1$ and $\Cal W_2$ are $\Cal A-$Hilbert modules,
such that $\Cal W_1$ is a closed invariant subspace of $\Cal W_2$
and ${\roman {dim}}_N(\Cal W_1)={\roman {dim}}_N(\Cal W_2) <\infty$,
then $\Cal W_1=\Cal W_2$ .
The von Neumann dimension of a Hilbert direct sum is the sum (possibly
infinite) of the von Neumann dimension of the summands.

The following proposition is well known.

\proclaim{\bf Proposition 1.6 }
{\sl    Assume that $\Cal W_1$ and $\Cal W_2$ are $\Cal A-$ Hilbert
 modules of finite type.
\newline (1) If $f\in \Cal L_{\Cal A}(\Cal W_1,\Cal W_2)$ and
$g\in \Cal L_{\Cal A}(\Cal W_2,\Cal W_1)$
then ${\roman {tr}}_N(fg)^n={\roman {tr}}_N(gf)^n$ for any $ n \geq 1.$
\newline (2) If $f:\Cal W_1\to \Cal W_2$ is an isomorphism and
$\alpha_i\in \Cal L_{\Cal A}(\Cal W_i), i=1,2$ so that
$f\cdot \alpha_1=\alpha_2\cdot f$ then
${\roman {tr}}_N\alpha_1={\roman {tr}}_N\alpha_2.$ }
\endproclaim

If $\Cal A'$ and $\Cal A''$ are two finite von Neumann algebras the
tensor product
$\Cal A'\otimes \Cal A''$ is defined as the weak closure of the image
of the algebraic
tensor product of $\Cal A'$ and $\Cal A''$ in
 $\Cal L(\Cal A'_2\otimes\Cal A''_2).$
The algebra $\Cal A'\otimes \Cal A''$ is again
 a finite von Neumann algebra whose trace has the property
that ${\roman {tr}}_N(a'\otimes a'')={\roman {tr}}_Na' {\roman {tr}}_Na''.$
If $\Cal W'$ and $\Cal W''$ are $\Cal A'-$ resp. $\Cal A''-$ Hilbert modules
of finite type then
$\Cal W'\otimes \Cal W''$ is an $\Cal A'\otimes \Cal A''$ Hilbert module of
finite type;
moreover, given $f'\in \Cal L_{\Cal A'}(\Cal W')$ and $f''\in
 \Cal L_{\Cal A''}(\Cal W''),$
${\roman {tr}}_N(f'\otimes f'')={\roman {tr}}_Nf' {\roman {tr}}_Nf''.$

\proclaim{\bf Definition 1.7}
A morphism $f: \Cal W_1 \to \Cal W_2$ is a {\it weak isomorphism} iff
${\roman {Null}}(f) =0$ and $\overline{{\roman {Range}}(f)}=\Cal W_2.$
\endproclaim

Using polar decomposition a weak isomorphism $f: \Cal W_1 \to \Cal W_2$
can be factored as $f=gf'$ where $f':\Cal W_1\to \Cal W_1$
is a weak isomorphism and $g:\Cal W_1\to \Cal W_2$ is an
isometric isomorphism given by $f'=(f^*f)^{1/2}, g= f\cdot(f^*f)^{-1/2}.$

\subheading{1.2 Determinant in the von Neumann sense}

Throughout this subsection we consider only $\Cal A$-Hilbert modules
 of finite type.

\proclaim{\bf Definition 1.8}
\newline (1) $\pi$  is an Agmon angle for $f\in \Cal L_{\Cal A}(\Cal W)$
iff there exists
$\epsilon > 0$ so that ${\roman {spec}}(f) \cap V_{\pi,\epsilon}
= \emptyset$ with $V_{\pi,\epsilon}$ defined as in the introduction.
\newline (2) $\pi$  is a weak  Agmon angle for $f\in \Cal L_{\Cal A}(\Cal W)$
iff
${\roman {spec}}(f) \cap V_{\pi,0}= \emptyset$
\endproclaim

We will first treat the case where $\pi$ is an Agmon angle. In particular
this implies that $f$ is an
isomorphism.
Define the complex powers of $f, f^s\in \Cal L_{\Cal A}(\Cal W)$,
$s\in \Bbb C ,$
by the formula
$$f^s=\frac {1}{2\pi i} \int_{\gamma} \lambda^s
(\lambda- f)^{-1}d\lambda, \tag{1.2}$$
where $\lambda^s$ is a branch of the complex power $s$ defined on
$\Bbb C_\pi=\Bbb C\smallsetminus \{z=\rho e^{i\pi}; \rho\in [0,\infty)\}$
and $\gamma$ is a closed
contour in $\Bbb C_\pi$ which surrounds the compact set  $spec f$ in
 $\Bbb C_\pi$ with
counterclockwise orientation. By Cauchy's theorem the contour $\gamma$ in
$(1.2)$ can be replaced by the contour
 $\gamma_{\pi, \epsilon}= \gamma_1\cup\gamma_2\cup\gamma_3$ where
$\gamma_1:=\{z=\rho e^{i\pi};\infty\le\rho\le\epsilon/2\}$,
${\gamma}_2:=\{ z=\frac {\epsilon}{2} e^{i\alpha}; \pi \ge\alpha\ge -\pi\}$
and
$\gamma_3:=\{z=\rho e^{i(-\pi)};\frac {\epsilon}{2} \le\rho\le\infty\}.$
Notice that $f^s$ is an entire function in $s\in\Bbb C$ with values in
$\Cal L_{\Cal A}(\Cal W)$
and ${\roman {tr}}_N(f^s)$
is an entire function on $\Bbb C. $ Therefore, if $\pi$ is an Agmon angle for
 $f,$ we can define
the determinant  $\det_Nf$ in the von Neumann sense by
$${\roman {logdet}}_N f = \frac {d}{ds} \mid_{s=0}
{\roman {tr}}_N (f^s).\tag {1.3}$$
We remark that this notion of determinant is equivalent to the definition
introduced by Fuglede and Kadison [FK].

If $f\in \Cal L_{\Cal A}(\Cal W_1,\Cal W_2)$ is an isomorphism then $f^*f$
is a selfadjoint positive isomorphism and one shows that
 ${\roman {det}}_N (f^*f)> 0.$
Define
$${\roman {Vol}}_Nf:=({\roman {det}}_N (f^*f))^{1/2}.$$

\proclaim {\bf Proposition 1.9} (1) Suppose $f_t\in \Cal L_\Cal A(\Cal W),$
with $t$ in an interval $I\subset \Bbb R,$
is a family of class $C^1$ of morphisms and $\pi$ is an Agmon angle for all
of them. Then
${\roman {logdet}}_N(f_t)$ is of class $C^1$ and
$$\frac{d}{dt}{\roman {logdet}}_N(f_t)={\roman {tr}}_N((\frac{d}{dt} f_t)
f_t^{-1}). \tag {1.4}$$
(2) Suppose $f_i\in \Cal L_\Cal A(\Cal W_i),i=1,2$ with $\Cal W_1$ and
$\Cal W_2$
$\Cal A$- Hilbert modules of finite type and
$\alpha:\Cal W_1\to \Cal W_2$ is an isomorphism so that
$\alpha f_1=f_2 \alpha.$ Then the following statements hold:
\newline (a) $spec f_1=spec f_2$ and therefore $\pi$ is an Agmon angle
for $f_1$ iff it is an
Agmon angle for $f_2.$ In this case ${\roman {logdet}}_Nf_1=
{\roman {logdet}}_Nf_2.$
\newline (b) $f_1$ is an isomorphism iff  $f_2$ is an isomorphism. In this
case ${\roman {Vol}}_Nf_1={\roman {Vol}}_Nf_2.$
\newline (3) Suppose $f\in \Cal L_\Cal A(\Cal W_1\oplus\Cal W_2)$ is
of the form
$$
f=\left( \matrix
f_1&0\\
g&f_2
\endmatrix\right).
$$ Then the following statements hold:
\newline (a) $spec f=spec f_1\cup spec f_2$ and therefore $\pi$ is an
Agmon angle for $f$ iff it is an
Agmon angle for both $f_1$ and $f_2.$ In this case
$${\roman {logdet}}_N f={\roman {logdet}}_Nf_1 +{\roman {logdet}}_Nf_2.
\tag {1.5A}$$
(b) $f$ is an isomorphism iff $f_1$ and $f_2$ are both isomorphisms.
 In this case
$$\log {\roman {Vol}}_N f=\log {\roman {Vol}}_Nf_1 +\log {\roman {Vol}}_Nf_2.
\tag {1.5B}$$
\newline (4) Suppose  $\Cal W_1,\ \Cal W_2 \text{ and }\Cal W_3$
 are $\Cal A$-Hilbert modules of
finite type. If $f_1\in \Cal L_\Cal A(\Cal W_1,\Cal W_2)$ and $f_2\in
\Cal L_\Cal A(\Cal W_2,\Cal W_3)$ are isomorphisms then
$f_2\cdot f_1\in \Cal L_\Cal A(\Cal W_1,\Cal W_3)$ is an isomorphism and
$$\log {\roman {Vol}}_N (f_2\cdot f_1)=\log {\roman {Vol}}_Nf_1 +
\log {\roman {Vol}}_Nf_2. \tag {1.6}$$
\newline (5) If $\alpha_i\in \Cal L_\Cal A(\Cal W_i), i=1,2$ are isometries and
$f:\Cal W_1\to \Cal W_2$ is an isomorphism then $\alpha_2 f \alpha_1\in
\Cal L_\Cal A(\Cal W_1,\Cal W_2)$ is an isomorphism and
$$\log {\roman {Vol}}_N (\alpha_2 f\alpha_1)=\log {\roman {Vol}}_Nf .$$
\endproclaim

{\sl Proof} All the statements can be proved in an elementary way.
For the convenience of the reader we indicate how to prove formula (1.6).
We have to prove that
$${\roman {logdet}}_N f_1^*f_2^*f_2f_1 \ - {\roman {logdet}}_N f_1^*f_1 \ =
{\roman {logdet}}_N f_2^*f_2.
$$
Consider the 1-parameter family $C(t)$ of positive, selfadjoint operators,
\newline $C(t):= f_1^*(f_2^*f_2)^{t}f_1,$ defined on $\Cal W_1.$
Using formula (1.4) one verifies that
$$\frac {d}{dt} {\roman {logdet}}_N f_1^*(f_2^*f_2)^{t}f_1
= \frac {d}{dt} {\roman {logdet}}_N(f_2^*f_2^*)^t.
$$
This leads to the claimed formula,
$$\align
 & {\roman {logdet}}_N f_1^*f_2^*f_2f_1 - {\roman {logdet}}_Nf_1^*f_1  = \\
 & \int_0^1 \frac {d}{dt} {\roman {logdet}}_N(f_1^*(f_2^*f_2)^{t}f_1)dt = \\
 & \int_0^1 \frac {d}{dt} {\roman {logdet}}_N(f_2^*f_2)^{t} dt =
 {\roman {logdet}}_N(f_2^*f_2).\ \  \diamondsuit
\endalign$$

First let us treat the case where $\pi$ is a weak Agmon angle for
$f\in \Cal L_\Cal A(\Cal W)$.
In this case $\pi$ is an Agmon angle for $f+\lambda$ with any $\lambda >0.$
One verifies that $
{\roman {logdet}}_N(f+\lambda)$ is a real analytic function in
$\lambda \in (0,\infty).$ We define
${\roman {logdet}}_N f$ as the element in $\bold D$ (cf.Introduction),
represented by the real analytic
function
$${\roman {logdet}}_N (f+\lambda)-\log \lambda {\roman {dim}}_N
{\roman {Null}}(f).\tag {1.8}$$
We note that parts (2) and (3) of Proposition 1.9 extend to this case as well.

Let $f\in \Cal L_\Cal A(\Cal W_1,\Cal W_2)$ be a weak isomorphism.
For $\lambda \geq 0$ denote by
$\Cal P_f(\lambda)$ the set of all $\Cal A$-invariant closed subspaces
$\Cal L\subset \Cal W_1$
such that, for $x \in \Cal L,$
 $\|f(x)\| \leq \sqrt\lambda \| x \|.$ Following Gromov-Shubin
([GS]) introduce the function
$F_f:[0,\infty)\to [0,\infty)$ defined by
 $F_f(\lambda):=\sup \{ {\roman {dim}}_N \Cal L; \Cal L \in \Cal
P_f(\lambda)\}.$
Observe that the function $F_f(\lambda)$ is nondecreasing, left continuous,
 $F_f(0)=0$ and $F_f(\lambda)={\roman {dim}}_N(\Cal W)$ for
$\lambda \geq \| f \|.$ Note that $f$ is an isomorphism iff there exists
$\lambda_0> 0$ s.t.
$F_f(\lambda)=0$
for $\lambda < \lambda_0.$ The Novikov-Shubin invariant $\alpha (f)$ associated
to a weak isomorphism $f\in \Cal L_\Cal A(\Cal W_1,\Cal W_2)$ is defined by
$$\alpha(f):=\liminf_{\lambda\to 0}\frac{\log F_f(\lambda)}{\log \lambda}
\in [0,\infty]. \tag {1.10}$$
Note that $\alpha(f)=\infty $ if $f$ is an isomorphism.

 If $f\in \Cal L_\Cal A(\Cal W_1,\Cal W_2)$
is an arbitrary morphism let
$$\overline f:\Cal W'_1=\Cal W_1 /{\roman {Null}}(f)\to
\overline {{\roman {Range}} f}=\Cal W'_2.
 \tag {1.11}$$
Note that $\overline f$ is a weak isomorphism and define $\alpha(f)$ by
$$\alpha(f):=\alpha(\overline f). \tag {1.12}$$

\proclaim{\bf Proposition 1.10 }
{\sl    (1) For any weak isomorphism $f\in \Cal L_\Cal A(\Cal W_1,\Cal W_2)$
 $$F_f(\lambda)=F_{(f^*f)^{1/2}}(\lambda)=F_{f^*}(\lambda).$$
\newline (2) If $f:\Cal W_1\oplus \Cal W_2\to \Cal W_1\oplus \Cal W_2$ is a
weak isomorphism of the form
$$
f=\left( \matrix
f_1&0\\
g&f_2
\endmatrix\right),
$$
then $f_1$ and $f_2$ are both weak isomorphisms and
$$
\sup \{F_{f_1}(\lambda), F_{f_2}(\lambda)\}\leq F_f(\lambda)
\leq F_{f_1}(\lambda) +F_{f_2}(\lambda).
$$
\newline (3) If $f\in \Cal L_\Cal A(\Cal W)$ is nonnegative and selfadjoint,
define the spectral projectors $Q_f(\lambda)\in \Cal L_\Cal A(\Cal W)$
corresponding to the interval $[0,\lambda]$ and $N_f(\lambda)
:={\roman {tr}}_N Q_f(\lambda).$
For $\lambda\geq 0$, and $\overline f$ given by (1.11)
$$N_f(\lambda)={\roman {dim}}_N({\roman {Null}}(f))
+F_{\overline f}(\lambda).\tag {1.13}$$
\endproclaim

The function $N_(\lambda)$ is called the spectral function of $f.$
Note that $F_{\overline f}(\lambda)$ is nondecreasing and $F_{\overline f}(0)
=0.$
$F_{\overline f}(\lambda)$ can be used to represent
$\log\det_Nf $ as a
Stieltjes integral
\newline $\int^{\infty}_{0+} log(\mu+\lambda)dF_{\overline f}(\mu).$

Denote  by $\bold F$ the set of functions
 $F:[0,\infty)\to [0,\infty)$ satisfying
\newline (1) $F(0)=0$;
\newline (2) $F(\lambda)$ is nondecreasing;
\newline (3) $F$ is continuous to the left,
\newline and recall the following definitions of Gromov-Shubin  (cf [GS])

\proclaim{\bf Definition 1.11}
(1) Functions $F,G\in \bold F$ are said to be dilational equivalent,
denoted  $ F\overset d\to \thicksim G,$
iff there exists $C> 0$ such that for $\lambda \geq 0$
$$G(C^{-1}\lambda)\leq F(\lambda)\leq G(C\lambda)\tag {1.14}$$
(2) Functions $F,G\in \bold F$ are said to be dilational equivalent near zero,
denoted $ F\overset d\to{\underset 0\to \thicksim} G$
iff there exist $C> 0$ and $\lambda_0> 0$ such that (1.14)  holds for
$\lambda < \lambda_0.$
\endproclaim
We end this subsection with the following observation. Suppose that
$\psi:\Cal A'\to \Cal A''$
is a homomorphism of finite von Neumann algebras which preserves the units
and the traces.
Then $\psi$  is injective. In particular, if  $\Cal A''$ is an
$\Cal A'$ -Hilbert module of finite
type, then any $\Cal A''$-Hilbert module of finite type
$\Cal W$ is
an $\Cal A'$-Hilbert module of finite type and
$\Cal L_{\Cal A''}(\Cal W)\subseteq \Cal L_{\Cal A'}(\Cal W).$

{\bf Remark 1.12} Assume that ${\roman {dim}}_{N,\Cal A'}
\Cal A'' = r \in \Bbb R.$
\newline (1) ${\roman {dim}}_{N,\Cal A'}(\Cal W)=r {\roman {dim}}_{N,\Cal A''}
(\Cal W).$
\newline (2) If $f\in \Cal L_{\Cal A''}(\Cal W)$ then
 ${\roman {tr}}_{N,\Cal A'}(f)= r{\roman {tr}}_{N,\Cal A''}(f).$
\newline (3) If $\pi$ is a weak Agmon angle for $f$ then
$${\roman {logdet}}_{N,\Cal A'}(f)=r {\roman {logdet}}_{N,\Cal A''}(f).$$

\subheading{1.3 Cochain complexes of finite type and torsion in the
von Neumann sense}

\proclaim{\bf Definition 1.13}
A cochain complex in the category of $\Cal A-$Hilbert modules of finite type,
 $\Cal C=(\Cal C_i,d_i)$, consists of a collection of Hilbert modules of finite
 type $\Cal C_i$,
all but finitely many zero, and a collection of morphisms
$d_i:\Cal C_i \to \Cal C_{i+1}$ which satisfies $d_i d_{i-1}=0$.
 In the sequel we always assume
that $\Cal C_i=0$ for $i < 0$ and refer to such a complex as a
 cochain        complex of finite type over $\Cal A$, or simply as a
 cochain complex of finite type.
\endproclaim

The     reduced cohomology of $\Cal C$, $\overline H^i(\Cal C),$
is defined by
 $$\overline H^i(\Cal C)={\roman {Null}} (d_i) / {}
 \overline {{\roman {Range}}(d_{i-1})}.$$
 Define the Betti numbers and Euler-Poincar\'e characteristic of $\Cal C$ by
 $$\beta_i(\Cal C):= {\roman {dim}}_N \overline H^i(\Cal C)
 ;~~~~~       \chi(\Cal C):=\sum_i(-1)^
i \beta_i(\Cal C), \tag {1.15}$$
and introduce a weighted version of the Euler-Poincar\'e characteristic,
 $$\psi(\Cal C):=\sum_i (-1)^i i \beta_i(\Cal C).\tag {1.16}$$
 Denote by $d^*_i:\Cal C_{i+1}\to \Cal C_i$ the adjoint of
 $d_i$, and consider $\Delta_i=d^*_i d_i+ d_{i-1} d^*_{i-1}$.
The operator $\Delta_i$ is a selfadjoint and nonnegative
 morphism.

\proclaim{\bf Definition 1.14}
(1) Given two cochain complexes of finite type over $\Cal A$, $\Cal C'$ and
 $\Cal C''$, a morphism
$\bold f:\Cal C'\to \Cal C''$ is given by a collection of morphisms
$f_i:\Cal C'_i \to \Cal C''_i$which
commute with the differentials $d_j$.

(2) A homotopy $\underline{\bold t}$ between $\bold f$ and
 $\bold g$      is
given by a collection of morphims       $t_i:C'_i\to C''_{i-1}$ satisfying
 $$f_i-g_i=d''_{i-1}t_i+t_{i+1}d'_i. \tag{1.17}$$

(3) Two cochain complexes $\Cal C'\text{ and } \Cal C''$ are homotopy
equivalent
if      there exist morphisms (in the category of complexes)
$\bold f:\Cal C'\to \Cal C'' \text{ and }\bold g: \Cal C'' \to \Cal C'$ so that
$\bold g\cdot \bold f $ resp. $\bold f\cdot \bold g$ is homotopic to
$\bold {id}_{\Cal C'}$ resp. $\bold {id}_{\Cal C''}.$
\endproclaim

Given a morphism $\bold f:\Cal C'\to \Cal C''$ denote by
$\overline H(\bold f)^i$ the
induced morphisms of $\Cal A$-Hilbert modules $\overline H(\bold f)^i:
\overline H ^i(\Cal C') \to \overline H ^i(\Cal C'').$  Note that if
$ f_i:\Cal C'\to \Cal C'', i=1, 2$  are two homotopic morphisms then
$\overline H(\bold f_1)^i=\overline H(\bold f_2)^i$ for all $i.$
Given a finite type cochain complex $\Cal C=(\Cal C_i,d_i)$ each $\Cal C_i$
can be
decomposed as a direct sum of mutually orthogonal subspaces
$\Cal C_i=\Cal H_i\oplus\Cal C^+_i\oplus  \Cal C^-_i$ with
$$\Cal H_i=\text{Null}\Delta_i;\text {  }
\Cal C^+_i=\overline{d_{i-1}(\Cal C_{i-1})}, \ \ \ \ \Cal C^-_i=
\overline{d^*_{i+1}(\Cal C_{i+1})}.\tag {1.18}$$
 This decomposition is called the Hodge decomposition .
The map $d_i$ can then be described by a $3\times 3$ matrix of the form
$$
d_i=\left( \matrix0&0&0\\
0&0&{\underline d_i}\\
0&0&0
\endmatrix \right).
\tag {1.19}$$
where $\underline d_i:\Cal C^-_i \to \Cal C^+_{i+1}$ is a weak isomorphism
and the
combinatorial Laplacian
 $\Delta_i=d_{i-1}d^*_{i-1}+d^*_id_i$ then takes the form of the
diagonal matrix $diag(0,d_{i-1}d^*_{i-1},d^*_id_i).$

Let $\bold f:\Cal C^1\to \Cal C^2$ be a morphism. With respect to the Hodge
decompositions
of $\Cal C^1_i$ and $\Cal C^2_i,$ the morphism $f_i: \Cal C^1_i\to
\Cal C^2_i$ can be
written as
a $3\times 3$-matrix of the form
$$
f_i=\left( \matrix
f_{i,11}&0&f_{i,13}\\
f_{i,21}&f_{i,22}&f_{i,23}\\
0&0&f_{i,33}
\endmatrix \right).
\tag {1.20}$$
where $f_{i,11}\in \Cal L_\Cal A(\Cal H^1_i,\Cal H^2_i)$,
$f_{i,22}\in \Cal L_\Cal A(\Cal C^{1,+}_i,\Cal C^{2,+}_i)$
$f_{i,33}\in \Cal L_\Cal A(\Cal C^{1,-}_i,\Cal C^{2,-}_i)$
and
\newline $\underline d^2_i \cdot f_{i,33}=f_{i+1,22}\cdot \underline d^1_i.$

\proclaim{\bf Definition 1.15} A cochain complex $\Cal C$ is called perfect
 if, for any i, $\underline{d}_i$ is an isomorphism.
\endproclaim

\proclaim{\bf Lemma 1.16} (1) Given a cochain complex $\Cal C=(\Cal C_i,d_i)$
one can find a modification $\tilde d_i$ of $d_i$ so that
$\tilde\Cal C=(\Cal C_i,\tilde d_i)$
is perfect and has the same Hodge decomposition as $\Cal C.$
\newline (2) Given an isomorphism
$\bold f:\Cal C^1\to \Cal C^2$ of cochain complexes
$\Cal C^k=(\Cal C^k_i, d_i),k=1,2$ one can find modifications
$\tilde d^k_i$ of $d^k_i$
so that
$$f_{i+1}\tilde d^1_i=\tilde d^2_i f_i \tag {1.21}$$
 and the cochain complexes
$\tilde\Cal C^k=(\Cal C^k_i,\tilde d^k_i),$ with $ k=1,2,$ are perfect and
have the same Hodge decompositions as $\Cal C^k.$
\endproclaim

{\sl Proof} Statement (1) follows by choosing $\tilde d_i$ of the form
$$
\tilde d_i=\left( \matrix0&0&0\\
0&0&\tilde{\underline d_i}\\
0&0&0
\endmatrix \right).
$$
where $\tilde{\underline d_i}$ is the isometry in the polar decomposition of
$\underline d_i$  given by $\tilde{\underline d_i} = {\underline d_i}
({\underline d_i^*}\underline d_i)^{-\frac{1}{2}}.$
\newline (2) With respect to the Hodge decomposition  of $\Cal C^1_i$
and $\Cal C^2_i$
define $\tilde d^1_i$ as in (1) and choose $\tilde d^2_i$ to be of the form
$$
\tilde d^2_i=\left( \matrix0&0&0\\
0&0&\tilde{\underline d^2_i}\\
0&0&0
\endmatrix \right)
$$
with $\tilde{\underline d^2_i}:=f_{i+1,22}\cdot\tilde{\underline d^1_i}
\cdot f_{i,33}^{-1}.$

In section 6.1 we will need the following

\proclaim{\bf Proposition 1.17} Suppose $\Cal C(t)=(\Cal C_i(t),d_i(t))$
is a family of
cochain complexes of finite type depending on a parameter $t\geq 0$,
and $\bold f(t):\Cal C(t)\to \Cal C$ is an isomorphism of cochain complexes
for any $t$. Introduce
 $\log V(t):=\sum_{q=0}^d (-1)^q \log {\roman {Vol}}_N
 \overline H(\bold f(t))^q.$
Assume that $\Cal C_i(t)$ and $\Cal C_i$ are free modules, that
 $e_{i,1}(t),...,e_{i,l_i}(t)$ is
an orthonormal base for $\Cal C_i(t)$, $e_{i,1},...,e_{i,l_i}$ is
an orthonormal base for $\Cal C_i$ and that $f_i(t)$, when expressed
with respect to these bases, is
an $l_i\times l_i$-matrix with entries in $\Cal A^{op}$
of the form $Id+O(1/t)$.

 Then $\log V(t)=O(1/t).$
\endproclaim
{\sl Proof} In view of Lemma 1.16 it suffices to prove the result for the
case where
$\Cal C(t)$ and $\Cal C$ are all perfect. In view of (1.20)
$$\log {\roman {Vol}}_N\overline H(\bold f(t))^i
=\log {\roman {Vol}}_N(f_{i,11}(t))
=\log {\roman {Vol}}_N(\pi_if_i(t)I_i(t)) \tag {1.22}$$
where $I_i(t)$ denotes the inclusion of $\Cal H_i(t)$ in $\Cal C_i$ and $\pi_i$
the orthogonal projection of $\Cal C_i$ onto $\Cal H_i$. Denote  by $P_i(t)$
the
orthogonal projections of $\Cal C_i(t)$ onto $\Cal H_i(t). $ We obtain
$$ \log {\roman {Vol}}_N(\pi_if_i(t)I_i(t))=
 \frac {1}{2} {\roman {logdet}}_N(P_i(t)f_i^*(t)\pi_if_i(t)I_i(t))=$$
$$\frac {1}{2} {\roman {logdet}}_N (P_i(t)f_i^*(t)f_i(t)I_i(t)-
 P_i(t)f_i^*(t)(Id-\pi_i)f_i(t)I_i(t)).\tag {1.23}$$
In view of the assumptions on $f_i(t)$, one has  $P_i(t)f_i^*(t)f_i(t)I_i(t)
 = Id+O(1/t).$
Next we will show that
 $$P_i(t)f_i^*(t)(Id-\pi_i)f_i(t)I_i(t) =O(1/t). \tag {1.24}$$
For this purpose note that the assumption
concerning the asymptotics of $f_i(t)$ combined with the hypothesis
that all complexes are perfect
imply the existence of  $\epsilon> 0$ and $t_0> 0$ such that for
$t\geq t_0$, $spec \Delta_i(t)\setminus \{0\}\text { and }
spec \Delta_i\setminus \{0\}$
are contained in $[2\epsilon,\infty)$. Therefore
$$P_i(t)= \frac {1}{2\pi i} \int_{S_\epsilon}(z-\Delta_i(t))^{-1}dz,
\text { and } \pi_i=\frac {1}{2\pi i} \int_{S_\epsilon}(z-\Delta_i)^{-1}dz,$$
where $S_\epsilon$ is the circle in the
complex plane of radius $\epsilon,$
centered at the origin. Using the assumptions on $f_i(t),$
one concludes that
$\Delta_i\cdot f_i(t)=f_i(t)\cdot\Delta_i(t)+O(1/t),$ and in view of
(1.23) one obtains
$$(Id-\pi_i)f_i(t)=f_i(t)(Id-P_i(t))+O(1/t).\tag {1.25}$$
As $(Id-P_i(t))I_i(t)=0$, estimate  (1.25) implies (1.24);
formulas (1.23) and (1.24) lead to
$$\log {\roman {Vol}}_N(\overline H(\bold f(t))^i
= \frac {1}{2} {\roman {logdet}}_N(Id+O(1/t))=O(1/t).\tag {1.26}$$
Given a cochain complex $\Cal C$ of finite type, introduce,
 following Gromov-Shubin ($[GS]$,
cf. also end of section 1.2), the functions
 $F_{\Cal C,i}(\lambda)\in \bold F$ defined
by $F_{\Cal C,i}(\lambda):=F_{\underline d_i}(\lambda)$ and the numbers
$\alpha_i$ defined by
$\alpha_i:=\alpha(\underline d_i).$ The following result is due to
Gromov-Shubin ($[GS]$)

\proclaim{\bf Proposition 1.18} Suppose $\bold f:\Cal C'\to \Cal C''$  and
$\bold g:\Cal C'' \to \Cal C'$ are two morphisms of cochain complexes so that
$\bold {id}_{\Cal C'}$
is homotopic to $\bold g\bold f$ by a homotopy $\underline{\bold t}=\{t_i\}.$
Then
$F_{\Cal C',i}(\lambda)\leq
F_{\Cal C'',i}(4\| f_{i+1}\|^2 \|g_i\|^2\lambda)$ for
$\lambda\leq \frac{1}{4||t_{i+1}||^2}.$
\endproclaim

In particular if $\bold f:\Cal C'\to \Cal C''$ is an isomorphism then
$F_{\Cal C',i}(\lambda) {\overset d\to\sim} F_{\Cal C'',i}(\lambda)$ and if
 $\Cal C'$ and $\Cal C''$ are homotopy equivalent,
then $F_{\Cal C',i}(\lambda)\underset 0\to{\overset d\to\sim} F_{\Cal C'',i}
(\lambda),$ and therefore
$\alpha'_i=\alpha_i''.$

The torsion $\log T(\Cal C)$ is the element in $\bold D$
defined by
$$\log T(\Cal C)=\frac {1}{2} \sum_i (-1)^{i+1}i{\roman {logdet}}_N\Delta_i.$$

{\bf Remark 1.19} The spectral functions
$N_i(\lambda):=N_{\Delta_i}(\lambda)$ satisfies
$$N_i(\lambda)=\beta_i+F_{i-1}(\lambda)+F_i(\lambda) \tag {1.28}$$
and therefore $\log T(\Cal C)$
can be represented  by the real analytic function  in $\lambda$
$$\frac {1}{2} \sum_i(-1)^i\int_{0+}^{\infty}\log(\mu+\lambda)dF_i
(\mu). \tag {1.29}$$

\proclaim {\bf Definition 1.20} The cochain complex of finite type $\Cal C$
is of determinant
class iff
\newline $\int_{0+}^1\log(\lambda)dN_i(\lambda) > -\infty$ for all $i$.
\endproclaim

We point out that if $\Cal C$ is of determinant class, then $\log T(\Cal C)$ is
in $\Bbb R \subset \bold D,$
and a sufficient condition for $\Cal C$ to be of determinant class is that
$\alpha_k> 0$
for $0\leq k\leq d.$

It will be convenient for the proof of the product formula below to
introduce for $\lambda> 0$ and $s\in \Bbb C$ with $\Re s > 0$,
$$\zeta_{\Cal C} (\lambda,s) = \frac {1}{2} \sum_i (-1)^i i \frac {1}{\Gamma
 (s)} \int_0^{\infty} t^{s-1} e ^{(-t\lambda)}{\roman {tr}}_N
e^ {-t\Delta_i}dt. \tag {1.30}$$
This function is real analytic in $\lambda$, complex analytic in $s$ for
 $\Re s> 0$
and admits an analytic continuation to the entire $s$-plane  so that $s=0$
is a regular point.
Using that
$$\zeta_{\Cal C} (\lambda,s) = \frac {1}{2}
 \sum_i (-1)^i i {\roman {tr}}_N ((\Delta_i + \lambda )^{-s}). \tag {1.31}$$
one sees that $\log T(\Cal C)$ is also represented by the real analytic
function in $\lambda$
given by $$\frac{d}{ds}|_{s=0}\zeta_{\Cal C} (\lambda,s)-\psi(\Cal C)
\log\lambda \tag {1.32}$$
where $\psi(\Cal C)$ denotes the weighted Euler-Poincar\'e characteristic.
This perturbation  of the
$\zeta$-function was also studied in $[Go].$

Suppose that $\Cal A_i$ $(i=1,2)$ are finite von Neuman algebras. Note
that $\Cal A=\Cal A_1\otimes\Cal A_2$ (tensor product in the category of
finite von Neumann algebras) is also a finite von Neumann algebra.
If $\Cal W_i$ are $\Cal A_i$- Hilbert modules of finite type then
the tensor product $\Cal W_1\otimes\Cal W_2$
(tensor product in the category of Hilbert spaces) is an
$\Cal A$-Hilbert module of finite type and
$${\roman {dim}}_{N}(\Cal W_1\otimes\Cal W_2) = {\roman {dim}}_{N}\Cal W_1
\cdot {\roman {dim}}_{N}\Cal W_2.$$

Let $\Cal C'$ resp.
$\Cal C''$, be two cochain complexes of finite type over $\Cal A_1,$ resp.
$\Cal A_2$. Denote by $\Cal C = \Cal C' \otimes \Cal C''$ the tensor product
of these cochain complexes,     $$\Cal C_i =\sum_{p+r=i}\Cal C'_p
\otimes\Cal C_r'',
\sum_{p+r=i}\ d_i=d'_p\otimes id +(-1)^p id\otimes d_r''.$$
Then $\Cal C$ is a cochain complex of finite type over $\Cal A.$

\proclaim{\bf Proposition 1.21 } (cf [CM], [LR])
{\sl    Let $\Cal C'$, resp. $\Cal C''$ be two finite type cochain
complexes        over $\Cal A_1$ resp. $\Cal A_2.$
Then, with $\Cal C = \Cal C' \otimes \Cal C''$,
\roster

\item $$ {\overline H^i (\Cal C)}=\sum_{p+r=i}{\overline H^p (\Cal C')}\otimes
{\overline H^r (\Cal C'')}$$

\item $$\zeta_{\Cal C}(\lambda ,s)=\zeta_{\Cal C'}(\lambda ,s) \chi(\Cal C'')+
\zeta_{\Cal C''}(\lambda ,s) \chi(\Cal C') $$

\item $$\psi(\Cal C)=\psi(\Cal C') \chi(\Cal C'')+
\psi(\Cal C'') \chi(\Cal C').$$
\endroster
\endproclaim

{\sl Proof:} The proof of $(1)$ can be found e.g. in [LR,Theorem 3.16] and
$(3)$ follows from $(1)$. To prove $(2)$ (cf [CM]) let $\Cal C$ be
a cochain complex of finite type over $\Cal A_1$,and
$\Cal N$ be  a $\Cal A_2-$Hilbert module of finite type. We will show below
that
$$\sum_q(-1)^q {\roman {tr}}_{N}(e^{-t\Delta_q})=\chi(\Cal C). \tag {1.33}$$
Therefore if $\beta:\Cal N\to \Cal N$ is a morphism, then
$$\sum_q(-1)^q {\roman {tr}}_{N}(e^{-t\Delta_q}\otimes\beta)
={\roman {tr}}_{N}(\beta)\chi(\Cal C). \tag {1.34} $$

To prove (1.33) we use the matrix representation of
 $\Delta_q$ with respect to the Hodge
decomposition, $diag (0,\underline d_{q-1} \underline d_{q-1}^*,
\underline d_q^* \underline d_q)$ and
Proposition 1.7 (1); they give
$${\roman {tr}}_Ne^{-t\Delta_{q-1}}|_{\Cal C_q^-}
= {\roman {tr}}_Ne^{-t\Delta_q}|_{\Cal C_q^+} $$ and consequently
$$
{\roman {tr}}_{N}e^{-t\Delta_q}=
{\roman {tr}}_{N}(e^{-t\Delta_q|_{\Cal C^+_q}}) +
{\roman {tr}}_{N}(e^{-t\Delta_q|_{\Cal C^-_{q+1}}}) +
{\roman {dim}}_N {\roman {Null}}(\Delta_q)$$
which leads to (1.33).
Next, decompose $\Delta_q=\bigoplus_{p+r=q}\Delta_{p,r}$, where
$\Delta_{p,r}=\Delta'_p\otimes id +id\otimes\Delta_r''$ to obtain
$$
2 \zeta_{\Cal C}(\lambda,s)=\sum_{p,r}(-1)^{(p+r)}(p+r)\frac{1}{\Gamma(s)}
\int_0^{\infty} t^{s-1}e^{-t\lambda}{\roman {tr}}_{N}(e^{-t\Delta'_p}
\otimes e^{-t\Delta_r''})dt
$$
$$= \frac{1}{\Gamma(s)}\int_0^{\infty}t^{s-1}
e^{-t\lambda}{\roman {tr}}_N(\{\bigoplus_p (-1)^p
pe^{-t\Delta'_p}\}\otimes\{\bigoplus_r(-1)^r e^{-t\Delta_r''}\})dt
$$
$$+ \frac{1}{\Gamma(s)}
\int_0^{\infty} t^{s-1}e^{-t\lambda}{\roman {tr}}_N
(\{\bigoplus_p (-1)^pe^{-t\Delta'_p}\}\otimes
\{\bigoplus_r(-1)^rr e^{-t\Delta_r''}\})dt
$$
which, in view of (1.34), is equal to
$$
\frac{1}{\Gamma(s)}\int_0^{\infty} t^{s-1}e^{-t\lambda}
\left(\chi(\Cal C'') \sum_p (-1)^pp {\roman {tr}}_N e^{-t\Delta'_p}
+\chi(\Cal C')\sum_r (-1)^rr
{\roman {tr}}_Ne^{-t\Delta_r''} \right) dt$$
 $$
=2\zeta_{\Cal C'}(\lambda,s)\cdot \chi(\Cal C'')+2\zeta_{\Cal
 C''}(\lambda ,s) \cdot \chi (\Cal C').
$$

\proclaim{\bf Corollary 1.22} (cf [CM]) With the same assumptions as
 in Proposition 1.21,
the following    identity, viewed in the vector space $\bold D$, holds:
$$\log T(\Cal C) =\chi (\Cal C'') \log T(\Cal C') +
\chi (\Cal C') \log T(\Cal C'').$$
\endproclaim

\subheading{1.4 $(\Cal A,\Gamma^{op})-$Hilbert module}

\proclaim{\bf Definition 1.23}
$\Cal W$ is an $(\Cal A, \Gamma^{op})-$Hilbert module of finite type if
\newline (BM1) $\Cal W$ is an $\Cal A-$Hilbert module of finite type;
\newline (BM2) $\Cal W$ is a $\Gamma^{op}-$Hilbert module, defined
by a unitary representation of $\Gamma$;
\newline (BM3) the action of $\Cal A$ and $\Gamma^{op}$ commute.
\endproclaim

{\bf Example} Let $X$ be
a countable set. Denote by $l^2(X)$ the Hilbert space obtained
by completion of
${\Bbb {R}} (X)=\lbrace f:X \to {\Bbb {R}} ; supp(f) \text{ is finite}
\rbrace$ with
respect to the
scalar
product $$\langle f_1,f_2\rangle  :=\sum_{x\in X}f(x)g(x).$$
Let $\Gamma$ be a countable group and $\Bbb R(\Gamma)$ denote the unital
$\Bbb R
$ -algebra
with
multiplication defined by convolution and $*$-operation induced by
the map $, g \to g^{-1}$.The algebra $\Bbb R(\Gamma)$ has a finite trace given
by $tr(f) := f(e)$ where $e$ denotes the unit element in $\Gamma$,
and acts from the left by convolutions on $l_2(\Gamma).$ This algebra can be
viewed as a $*$-subalgebra
of $\Cal L_{\Gamma}(l_2(\Gamma),l_2(\Gamma)).$ Denote by $\Cal N(\Gamma)$ its
weak closure
in $\Cal L_{\Gamma}(l_2(\Gamma),l_2(\Gamma)).$ Then $\Cal N(\Gamma)$ is a
finite von Neumann algebra.

Let $\rho :\Gamma\times X\to X$ be a {\it left} action
of $\Gamma$
on the set $X$ with
finite isotropy groups. $\rho $ induces a left action of $\Gamma$ by
isometries which makes  $l_2(X)$ an
$\Cal N(\Gamma)$-Hilbert module ; if the quotient set $\Gamma\backslash X$
is finite, then this module is a Hilbert module of finite type. Suppose, in
addition, that $\Gamma'$ is another countable group
and $\rho':X\times\Gamma'\to X$ is a {\it right} action of $\Gamma'$ on $X$ so
that
$\rho$ and $\rho '$ commute. $\Gamma'$ induces an action by isometries on
$l_2(X)$ which makes $l_2(X)$
an $(\Cal N(\Gamma), {\Gamma'}^{op})$-Hilbert module of finite type.
As an example, consider the case $X=|\Gamma|$, the underlying set of
$\Gamma$,
$\Gamma=\Gamma'$ and $\rho$ and $\rho '$ given by
$\rho(g_1,g_2)=g_1 g_2,$ and $\rho '(g_2,g_1)=g_2 g_1^{-1}.$
Then $l_2(\Gamma)$ is an
$(\Cal N(\Gamma),\Gamma^{op})$-Hilbert
module of finite type, refered to as the {\it regular birepresentation.}
\heading
1.5 Bundles of $\Cal A-$Hilbert modules
\endheading
\proclaim{\bf Definition 1.24}
A smooth bundle $p: \Cal E \to M$ over a smooth manifold $M$
 is a bundle of $\Cal A-$Hilbert modules of finite type with fiber $\Cal W$ if
\newline (B1) $p:\Cal E \to M$ is a smooth bundle of infinite dimensional
topological vector spaces,
equipped with a Hermitian structure $\mu$ which makes each fiber
$p^{-1}(x),x\in
 M,$
into a separable Hilbert space;
\newline (B2) $\Cal E$ is equipped with a smooth fiberwise action
$\rho :\Cal A \times \Cal E \to M$
which makes each fiber $p^{-1}(x)$ an $\Cal A-$Hilbert module of finite type.
\newline (B3) $\Cal W$ is an $\Cal A$-Hilbert module of finite type and
$p: \Cal E \to M$ is locally isomorphic to $p_o : \Cal W \times M \to M$
where the local isomorphism intertwines $p, p_o$, the Hermitian structures
and the $\Cal A-$actions.
\endproclaim

{\bf Example} Let $M$ be a closed smooth manifold with fundamental group
 $\Gamma :=\pi_1(M)$
and let $\Cal W$ be an
$(\Cal A,\Gamma^{op})$-Hilbert module of finite type. Let $\tilde p:\Cal
 W\times\tilde M\to \tilde M$ be the trivial
smooth bundle of $\Cal A-$Hilbert modules;      $\tilde p$ is
$\Gamma$-equivariant with
respect to the diagonal action of $\Gamma$
on $\Cal W\times_{\Gamma}\tilde M $ and the left action of $\Gamma$ on $\tilde
 M$. Therefore $\tilde p$
induces $p: \Cal E=\Cal W\times_{\Gamma}\tilde M\to M$ which is a smooth bundle
 of
$\Cal A-$Hilbert modules of finite type. This bundle is the canonical bundle
over $M,$ associated
to $\Cal W$.

\vfil\eject

\heading
{2.Calculus of pseudodifferential operators acting on
${\Cal A}$-Hilbert bundles of finite type}
\endheading

In this section we construct a calculus of pseudodifferential
operators, called pseudodifferential $\Cal A$-operators,
on a compact manifold,
where $\Cal A$ is a finite von Neumann algebra
(cf e.g. [FM],[Le] ,[Lu] for related work).

\subheading{  2.1 Sobolev spaces, symbols and kernels}

Let $B$ be a Banach space. For $u \in \swdb,$ the space of functions $u:\rd
\longrightarrow {B}$ of Schwartz class, $\|u\|_s$ denotes the
Sobolev $s$-norm given by
$$\align
\|u\|_s^2 & :=  \int_{\rd}(1+|\xi|^2)^s \|\hat{u}(\xi)\|^2d\xi
\endalign$$
where $\hat{u}(\xi)$ denotes the Fourier transform of $u.$

\proclaim{Definition 2.1}{ (1)The Sobolev space $\hswb$ is the completion of
 $\swdb$ with
respect to the Sobolev $s$-norm; equivalently, it can be defined as
the space of all distributions $u\in \swddb$
with
$$
(1 + |\xi |^2)^{s/2} \hat{u} \in L_2(\rd ;B).$$
\newline (2) The space $H^{loc}_s(\Bbb R,B)$ is the space of all distributions
$u \in
\ddwb$ such that $\phi u \in \  \hswib$ for any $\phi \in
 C_0^\infty(\rd).$}
\endproclaim

Note that the Sobolev spaces $\hswb$ have the same properties
as the usual Sobolev spaces except that Rellich's compacity theorem
does not hold.

Let $\Cal W$ be an $\Cal A-$Hilbert module.
The space $\hsw$ is an $\Cal A$-Hilbert module whose dual can be identified
with $H_{-s}(\rd, {\Cal W}).$  Note also that $H^{loc}_s(\Bbb R,\Cal W)$ is an
$\Cal A$-module.
Extending the classical case ${\Cal A} = {\Bbb R},$ symbols are
defined as follows:

\proclaim{Definition 2.2 }{ (1) A function $a \in \crdww$ is a
symbol of order $m \in {\Bbb R},$ denoted by $a \in \swm =
S_{{\Cal W}}^m(\rd \times \rd),$ if the following conditions hold:
\newline (Sy1) $a(x,\xi)$ has compact support in $x;$
\newline (Sy2) for any multiindices, $\alpha$ and $\beta,$ there exists a
constant, $C_{\alpha \beta},$ such that
$$\|\partial_x^\alpha \partial_\xi^\beta a(x,\xi) \| \leq C_{\alpha
\beta}(1 + |\xi|)^{m - |\beta |} .\tag {2.1}$$
\newline (2) (cf [Sh1]) {A symbol $a \in \swm$ is {\sl
classical} if it admits an expansion of the form
\newline $\sum_{j \geq 0}
\psi(\xi) a_{m-j}(x,\xi)$ where $ \psi \in \crd$ satisfies $$\psi(\xi)
=\cases          0 & \hbox{  for } |\xi| \leq \frac{1}{2} \\
        1 & \hbox{  for }|\xi| \geq 1
\endcases
$$
\newline (Sy3) $a_{m-j} \in C^\infty(\rd \times (\rd \setminus
\{0\}), {\Cal L}_{{\Cal A}}({\Cal W}, {\Cal W}))$ has compact
$x$-support and is positively homogeneous of degree $m-j;$
\newline (Sy4) $a(x,\xi) - \sum_{j=0}^{N-1} \psi(\xi) a_{m-j}(x,\xi) \in
S_{{\Cal W}}^{m-N}$ for all $N \geq 0.$ }
\endproclaim

Subsequently, we always assume that all symbols are classical.

Given $a \in \swm,$ define a linear operator $A:\crdwo \longrightarrow
\crdwo$ by
$$
Au(x)= \frac{1}{(2\pi)^d}  \int_{\rd} d\xi \int_{\rd}
dxe^{i(x-y,\xi)} a(x,\xi) u(y). \tag {2.2}
$$
The principal symbol
of $A, \ \sigma_A(x,\xi) = a_m(x,\xi),$ is invariantly defined as a
smooth function on $T^*\rd \setminus \{0\}$ with values in $ {\Cal
L}_{{\Cal A}}({\Cal W}, {\Cal W}).$

The operator $A$ is said to be
a pseudodifferential $\Cal A$-operator of order $m,$ denoted $A \in \pbm,$
and can be extended to a bounded, linear
operator (any $s \in \Bbb R$)
$$\align
A & :  \hsw \longrightarrow H_{s-m}(\rd, {\Cal W}).
\endalign$$

The Schwartz kernel of $A, \ K_A(x,y),$ is given formally as an
oscillatory integral
$$
K_A(x,y) =  \int_{\rd} e^{i(x-y,\xi)} a(x,\xi) d\xi. \tag {2.3}
$$
To characterize the properties of this distributional kernel we
introduce, following \cite{H\"o, p.100},

\proclaim{Definition 2.3} {A distribution $K \in \srdww$ with values in
\newline ${\Cal L}_{{\Cal A}}({\Cal W},{\Cal W}),$ is conormal of order $m$
along the
diagonal $Diag(\Bbb R^d)= \{(x,x); x \in \Bbb R^d \},$ denoted $K\in
\imw,$ if for all $C^\infty$-vectorfields $V_1, \dots ,V_N$ ($N\geq 1$
arbitrary)  which are tangential
to $Diag(\Bbb R^d)$ and $\phi \in C_0^{\infty}(\rd \times \rd)$
$$\align
V_1 \cdots V_N \phi K & \in  \  H_{-m-\frac{2d}{4}}(\rd
\times \rd, {\Cal L}_{{\Cal A}}({\Cal W},{\Cal W})).
\endalign$$}
\endproclaim

Following \cite{H\"o}, one verifies the following

\proclaim{Lemma 2.4 }{Let $a \in \swm$ and let $K_A$ be the corresponding
distributional kernel as defined above.  Then $K_A \in \imw.$}
\endproclaim

We note that if $m< -d,$ the kernel $K_A(x,y)$ is
continuous.

\proclaim{Definition 2.5}{A pseudodifferential operator $A$ in $\prdm$ is
said to be
elliptic on $U_1 \subset \rd$ if the principal symbol $a_m(x,\xi)$ is
invertible for $\xi \in \rd \setminus \{0\}$ for all $x \in U_1$ and
$$
\|a_m(x,\xi)^{-1}\| \leq  C_1 (1 + |\xi |)^{-m} \ \hbox{ for }x \in
U_1, \ |\xi| \geq 1. \tag {2.4}
$$
}
\endproclaim

Note that as $a_m(x,\xi) \in  {\Cal L}_{{\Cal A}}({\Cal W}, {\Cal W})$
the inverse satisfies $a_m(x,\xi)^{-1} \in  {\Cal L}_{{\Cal A}}({\Cal
W}, {\Cal W}). $

\subheading{  2.2 Pseudodifferential operators on bundles of Hilbert modules}

Throughout this section let $(M,g)$ denote a compact Riemannian manifold of
dimension $d,$ possibly with boundary, ${\Cal A}$ a finite von Neumann
algebra, ${\Cal W}$ an ${\Cal A}$-Hilbert module of finite type and
$p:{\Cal E } \longrightarrow M$ a bundle of ${\Cal A}$-Hilbert modules
 with fiber ${\Cal W}.$

Introduce the Banach bundles of bounded linear operators ${\Cal L}
\longrightarrow M \times M$ and ${\Cal B} = {\Cal L}_{\Cal A}
\longrightarrow M \times M$ whose fibres at $(x,y) \in M \times M$ are
given by
$$
{\Cal L}_{xy}  =  {\Cal L}({\Cal E}_y, {\Cal E}_x); \ \ \ \ \
{\Cal B}_{xy}  =  {\Cal B}({\Cal E}_y, {\Cal E}_x).
$$
where $ {\Cal L}({\Cal E}_y, {\Cal E}_x) $ denotes the Banach space of
all bounded linear operators from the fiber ${\Cal E}_y$ to the fiber
${\Cal E}_x$ and
$$\align
{\Cal B}_{xy} & =  \{f \in {\Cal L}_{xy} ; f \hbox{ is }{\Cal
A}-\hbox{linear} \}.
\endalign$$
In a straightforward manner one may verify that
the Banach bundle $\omega :{\Cal B} \longrightarrow M \times M
$ has the following properties
\newline (Bu1) ${\Cal B}_{xy}$ is a weakly closed linear subspace of ${\Cal
L}_{xy};$
\newline (Bu2) if $b \in {\Cal B}_{xy},$ then $b^* \in {\Cal B}_{yx};$
\newline (Bu3) if $b \in {\Cal B}_{xy}, \ b' \in {\Cal B}_{yz},$ then
$b b' \in  {\Cal B}_{xz};$
\newline (Bu4) ${\roman {Id}} \in {\Cal B}_{xx};$
\newline (Bu5) if $a \in {\Cal B}_{xx}$ is invertible then $a^{-1} \in {\Cal
B}_{xx}.$

Let $U$ be an open connected subset of $M$ and $X= \Bbb R^d$ or, in case $U$ is
a neighborhood of a boundary point of $M$, $X= \Bbb R^d_+ := \{
(x_1,...,x_d); x_d \ge 0 \}$.

\proclaim {Definition 2.6}{A pair $(\phi, \Phi)$ of smooth maps
$\phi :U \to X$ and
$\Phi : \Cal E | U \to X \times \Cal W $ is said to be a coordinate
chart of $(M, \Cal E \to M)$ if $\phi$ is a chart of $M$ and $\Phi$ is
an $\Cal A$-trivialization of $\Cal E \to M$ over $U$.

In particular, $\Phi_x := \Phi |_{ p^{-1}(x)} : p^{-1}(x) \to \Cal W$ is
an isometry.}
\endproclaim

By a standard localizing procedure, one defines the Sobolev spaces
$\hse = H_s(M,{\Cal E})$ using the  definition of Section 2.1 and a smooth
partition of unity subordinate to an open cover of $M$ which comes from
an atlas of ${\Cal E} \longrightarrow M.$
(Equivalently, the Sobolev norms can be defined using a Riemannian metric on
$M$ and a connection on $\Cal E.$)
The inner product in $\hse$
will depend on the particular choice of the partition; however a
different choice of partition of unity will lead to an equivalent
inner product.

The Sobolev $s$-norm of an element $u \in \hse$ will be denoted by
$\|u\|_s.$We point out that for $s \geq t, \ \hse$ imbeds into $H_t({\Cal E}).$
This embedding, however, is not compact if ${\Cal W}$ is of infinite
dimension (i.e. Rellich's lemma does not hold) even when $M$ is closed.

 To simplify the exposition we again assume that $M$ is closed.
\proclaim{Definition 2.7} {(1) A linear operator
$A: \ce \longrightarrow \ce$ is an $\Cal A-$smoothing operator, if $A$
is of the form
$$(Au)(x) = \int _M K_A(x,y)u(y)dy$$
where the Schwartz kernel $K_A$ of $A$ is a smooth section of the bundle
$\Cal B \longrightarrow M\times M.$
\newline (2) A linear operator $A: \ce \longrightarrow \ce$ is a
pseudodifferential $\Cal A-$operator of order m if for some atlas
$(\phi _j, \Phi _j)_{j \in J}$ of $\Cal E\longrightarrow M$,
$A = \sum_j A_j + T$ where $T$ is an $\Cal A-$smoothing operator
and the operators $A_j$ are operators with support in the domain
of $\phi_j$ and,when expressed in local coordinates, pseudodifferential
$\Cal A-$operators of order m.}
\endproclaim

One shows that $A \in \pbm$ is a linear operator, $A: \ce
\longrightarrow \ce,$ which can be extended, for any $s \in {\Bbb R},$
to a bounded linear operator
$$\align
A & :  \hse \longrightarrow H_{s-m}({\Cal E}).
\endalign$$
The principal symbol $\sigma_A$ of
$A$ can be defined invariantly as a smooth function
$\sigma_A(x,\cdot): T^*M \setminus \{0\} \longrightarrow {\Cal
L}_{{\Cal A}}({\Cal E}_x).$

Note that $\pbmi = \cap_m \pbm$ is the space of $\Cal A-$smoothing operators.
While it is clear that operators in $\pbmi$ have smooth kernels, it is
in general not true that such operators are compact.

As in the classical theory one develops a calculus for these
pseudodifferential operators.  In particular, one shows that the
composition $A \circ B$ of two pseudodifferential operators $A$ and
$B$ as well as the adjoint $A^*$ (with respect to the Hermitian
structure on ${\Cal E} \longrightarrow M$) are pseudodifferential
operators of the expected order.

\subheading{  2.3 Elliptic pseudodifferential operators}

To simplify the exposition we assume that $M$ is closed.
\proclaim{Definition 2.8}{An operator $A \in \pbm$ is said to be elliptic if
the
principal symbol of $A, \ \sigma_A(x,\xi),$ is invertible for all $x
\in M$ and all $\xi \in T^*_xM \setminus \{0\}.$}
\endproclaim

As in the classical case one can construct a parametrix, $R(\mu),$ for
the operator $(\mu - A)$ when $A$ is elliptic and $\mu \in {\Bbb C}
\setminus \bigcup_{(x,\xi) \in  T^*M \setminus \{0\}} \spec
(\sigma_A(x,\xi)) .$

The operator $R(\mu)$ is an element of $\pbmm$ and represents an
inverse of $(\mu - A)$ up to smoothing operators.  Let $U$ be a chart of $M$
which belongs to an atlas of $\Cal E \longrightarrow M$.  Denote by $\phi$
and $\Phi$ the
diffeomorphisms
$$\align
\phi & :  \rd \longrightarrow U \subset M \\
\Phi & :  \rd \times {\Cal W} \longrightarrow {\Cal E}_{|_{U}}
\endalign$$
where $U$ is an open subset of $M$ and $\Phi$ trivializes the bundle
$p: {\Cal E} \longrightarrow M$ over $U$ such that $p \phi = \Phi
p_1$ with $p_1 : \rd \times {\Cal W} \longrightarrow \rd.$

The symbol of $R(\mu)$ in the chart $U$ has an asymptotic expansion
determined inductively as follows:
$$\align
r_{-m}(x,\xi,\mu) & =  (\mu - a_m(x,\xi))^{-1}
\endalign$$
and, for $j \geq 1,$
$$\align
r_{-m-j}(x,\xi,\mu) & = r_{-m}(x,\xi,\mu)\left(\sum_{k=0}^{j-1}
\sum_{|\alpha|+l+k = j} \frac{1}{\alpha ! } \partial^\alpha_\xi
a_{m-l}(x,\xi) D^\alpha_x r_{-m-k}(x,\xi,\mu)\right)
\tag{2.5}\endalign$$
where $\alpha$ is a mutiindex, $\alpha = (\alpha_1,\dots, \alpha_d), \
\alpha ! = \alpha_1 ! \alpha_2 ! \cdots \alpha_d !,$ and $D_x^\alpha
= (\frac{1}{i}\partial_x)^\alpha.$  The term $r_{-m-j}(x,\xi,\mu)$ is
an element of ${\Cal L}_{{\Cal A}}({\Cal W},{\Cal W})$ and is
positively homogeneous of degree $-m-j$ in $(\xi, \mu^{\frac{1}{m}}):$
$$\align
r_{-m-j}(x,\lambda \xi,\lambda^m \mu) & =  \lambda^{-m-j}
r_{-m-j}(x,\xi,\mu)
\endalign$$
for any $\xi \in \rd \setminus \{0\}$ and $\lambda > 0.$

In the classical case the parametrix of an invertible elliptic
pseudodifferential operator is readily used to show that the inverse
of an operator is also pseudodifferential.  For our more general
calculus, additional arguments are necessary.

\proclaim{Proposition 2.9}{Assume that $M$ is closed and that $A\in \pbm $
is elliptic.  If $A$ considered as a bounded linear operator,
$A:H_m({\Cal E}) \longrightarrow L_2({\Cal E}),$ is one-to-one and
onto, then $A^{-1} \in \pbmm.$  }
\endproclaim

{\sl Proof:}
 Denote by $B \in \pbmm$  a parametrix for $A$.  The operators $T_1 :=
AB - {\roman {Id}}$ and $T_2 = BA - {\roman {Id}}$ are in $\pbmi.$
 From this we
conclude that $A^{-1} = B - A^{-1}T_1.$  The statement
follows once we prove that $A^{-1}T_1 \in \pbmi$ or, equivalently,
that $A^{-1}T_1$ has a smooth Schwartz kernel, $K_{A^{-1}T_{1}}(x,y)
\in {\Cal B}_{xy}\subset  {\Cal L}({\Cal E}_y,{\Cal E}_x).$
This is proved by using a technique due to Shubin [Sh2].
$~~~\diamondsuit$

To state Lemma 2.10 we introduce the space $\plm$ of
pseudodifferential operators of order $m$ with coefficients in ${\Cal
L}({\Cal E}).$  They are defined as were the operators in $\pbm$ by
simply replacing the bundle ${\Cal B} \longrightarrow M \times M$ with
the bundle ${\Cal L}({\Cal E}) \longrightarrow M \times M$ whose
fibres ${\Cal L}({\Cal E})_{xy}$ are given by ${\Cal L}({\Cal E}_y,
{\Cal E}_x).$
\proclaim{Lemma 2.10}{Let $M$ be a closed manifold of dimension $d.$
 $\pbm$ is a closed subspace in $\plm$ with respect to the
topology provided by the operator norm, $|\| A \||,$ where $A \in
\plm$ is viewed as a bounded linear operator $A : H_m ({\Cal E})
\longrightarrow L_2({\Cal E}).$
\endproclaim
{\sl Proof} In view of Lemma 2.9 it suffice to prove the statement for $m <
-d.$  In that case the result follows by noting that in view of Lemma 2.4 the
Schwartz kernel $K_A$ of $A \in \pbm$ is a
continuous function of $(x,y) \in M \times M$ and $A$ with $K_A(x,y)
 \in {\Cal B}_{xy}.$  As ${\Cal
B}_{xy}$ is a closed subspace of ${\Cal L}({\Cal E}_y,
{\Cal E}_x),$ statement (1) follows.

As in the classical theory one proves the following estimates for the
resolvent:

\proclaim{Lemma 2.11}{(cf \cite{Se1},\cite{Sh1}) Assume that $A \in \pbm$
is an elliptic
operator of order $m \geq 0$ such that $A:H_m({\Cal E})
\longrightarrow L_2({\Cal E})$ is one-to-one and onto. Further assume that
$\pi$ is an Agmon angle for $A.$
\newline Then for
$\lambda < 0$ with $|\lambda |$ sufficiently large and for $0 \leq m'
\leq m,$  $$|\| (\lambda -
A)^{-1} \| |_{0 \to m'}  \leq
C_{m'} |\lambda|^{-1 + \frac{m'}{m}} \tag {2.7}$$for some constants
$C_{m'} > 0.$
}
\endproclaim

\subheading{  2.4 Zeta-functions and regularized determinants of an
invertible elliptic operator}

Let $(M,g)$ be a closed Riemannian manifold.  Assume that $A \in \pbm$
is elliptic and of positive order, $ m > 0,$ with $\pi$ as
an Agmon angle; i.e. there exists $\epsilon > 0$ such that
\newline (1) $V_{\pi,\epsilon} \cap \spec A = \emptyset;$
\newline As a consequence $\pi$ is also a principal angle
\newline (2) $V_{\pi,\epsilon} \bigcap \left(\cup_{x \in M, \ (x,\xi)
\in S_x^*M} \spec (\sigma_A(x,\xi))\right) = \emptyset.$

The solid angle $V_{\pi,\epsilon}$ is defined as in the introduction. Note that
(1) implies the invertibility of $A$ viewed as a bounded linear
operator, $A: H_m({\Cal E}) \longrightarrow L_2({\Cal E}).$  Moreover,
for ${\Re}s < 0,$ one can define the complex powers of $A$ by
$$\align
A^s & :=  \frac{1}{2\pi i} \int_{\gamma_{\pi,\epsilon}}
\mu^{-s}(\mu -A)^{-1} d\mu
\tag {2.8}\endalign$$
where $\gamma_{\pi, \epsilon}$ is a path in ${\Bbb C}$ as defined in the
 introduction.
For $s $ satisfying $0 \leq k-1 \leq {\Re}s < k \in {\Bbb N}$ one
defines
$$\align
A^s & =  A^k A^{s-k}.
\endalign$$
It follows from Proposition 2.9 using arguments due to
Seeley \cite{Se1}, that $A^s \in \Psi DO^s_{{\Cal B}}(M)$ (after
suitably generalizing the concept of order to complex numbers $s \in
{\Bbb C}$), depending holomorphically on $s.$  Moreover, for
$\Re s < -\frac{d}{m},
\ A^s$ has a von Neumann trace
$$\align
\trn(A^s) & :=  \int_M \trn K_{A^s}(x,x) dx
\endalign$$
where $K_{A^s}(x,y) \in {\Cal L}({\Cal E}_y,{\Cal E}_x)$ denotes the
Schwartz kernel of $A^s.$

For $\alpha \in C^\infty(M,{\Bbb C})$ and ${\Re}s > \frac{d}{m}$
one defines the generalized zeta-function
$$
\zeta_{\alpha,A}(s) =  \trn(\alpha A^{-s}). \ \ \ \ \
$$
As in \cite{Se1} (cf also \cite{Gi} Lemma 1.7.7) one shows
\proclaim{Theorem 2.12}{(cf \cite{Se1})
\newline (1) Assume $A \in \pbm$ where $m >0$ and $A$ is elliptic with
$\pi$ as an Agmon angle.  If $\alpha  \in
C^\infty(M,{\Bbb C}),$ then $ \zeta_{\alpha,A}(s)$ admits a
meromorphic extension to  the entire $s$-plane.  The extension has
at most simple poles and $s=0$ is a regular point.  The value of $
\zeta_{\alpha,A}(s) $ at $s = 0$ is given by
$$\align
 \zeta_{\alpha,A}(0) & =  \int_{M}\alpha(x) I_d(x)
\tag {2.9}\endalign$$
where $I_d(x)$ is a density on $M.$  In an appropriate coordinate
chart, $I_d(x)$ is given by
$$\align
I_d(x) & =  \frac{1}{m} \frac{1}{(2\pi)^d} \int_{|\xi| = 1} d\xi
\int_0^\infty \trn (r_{-m-d}(x,\xi,- \mu))d\mu.
\tag {2.10}\endalign$$
If $A$ is a differential operator and $d = {\roman dim}M$ is odd, then
$I_d(x) \equiv 0.$
\newline (2) Assume that $A(t):H_m({\Cal E}) \longrightarrow L_2({\Cal E})$
is a family of classical pseudodifferential operators of order $m$
depending in a $C^r$-fashion on a parameter $t$ varying in an
open set of ${\Bbb R}.$  Assume that $A(t)$ is elliptic and that $\pi$
is an Agmon angle for any $t,$ uniformly in $t.$  Then
$\zeta_{A(t)}(s)$ is a family of functions holomorphic in $s$ in a
neighborhood of $s = 0$ which depends in a $C^r$-fashion on $t.$}
\endproclaim

Theorem 2.12 above allows us to introduce the
$\zeta$-regularized determinant of an elliptic operator $A \in \pbm$
of order $m > 0$ with $\pi$ as an Agmon angle:
$$\align
{\roman {det}} A & :=  \exp \left\{ -\left. \frac{d}{ds}\right|_{s=0}\zeta_A(s)
\right\}.\tag {2.11}
\endalign$$

To treat the case where $A$ is not invertible, first note the following

\proclaim {Lemma 2.13} {Assume $A \in \pbm.$ Then the nullspace of
$A$, ${\roman {Null}}(A)$, is an $\Cal A$-Hilbert module of finite von
Neumann dimension, $dim_N({\roman {Null}}(A)).$}
\endproclaim
Assume that $A$ is an elliptic operator, $A \in \pbm,$ of order
$m > 0$ with $\pi$ as a weak Agmon angle
(i.e ${\roman {spec}}( A)\cup (-\infty,0)=\emptyset$) . Then the
operator $A+\lambda$ with $\lambda >
0$ has $\pi$ as an Agmon angle and $\ldet(A + \lambda)$ is a real analytic
function in
$\lambda$. Define $\ldet_N (A)$ to be the element in $\bold D$
 represented by the analytic function

$$\align
\ldet_N (A) & :=  \ldet_N(A + \lambda) - \dimn ({\roman {Null}}(A)) \log
\lambda.\tag {2.12}
\endalign$$
\proclaim{Definition 2.14} {$A$ is of {\sl determinant class} if
$$\align
 &   \lim_{\lambda \downarrow 0} \left( \ldet_N(A + \lambda) - \dimn
({\roman {Null}}(A)) \log \lambda \right) \tag {2.13}
\endalign$$
exists.  In this case, $\ldet_N A$ is a real number.}
\endproclaim

If $A$ is selfadjoint and nonnegative, there is a functional calculus
for $A.$  In particular, one can introduce the spectral projections
$Q(\lambda)$ corresponding to the intervals $(-\infty , \lambda].$
Using Proposition 2.9 and the assumption that $A$ is
nonnegative, one verifies that $Q(\lambda)$ is in $\pbmi$ for any value
of $\lambda \in {\Bbb R}.$  Denote the distibution kernel of $Q(\lambda)$ by
$K_\lambda$ and define the spectral function
$$\align
N_A(\lambda) & :=  \int_M \trn K_\lambda(x,x) dx.\tag {2.14}
\endalign$$
Note that $N_A(\lambda)$ is nonnegative, right continuous and monotone
increasing as a function of $\lambda.$  Moreover, $N_A(\lambda) = 0$
for $\lambda < 0$ and for an appropriate constant $C>0,$
$$\align
N_A(\lambda) & \leq  C|\lambda|^{\frac{d}{m}}.
\endalign$$
\proclaim{Proposition 2.15}{Assume that $A \in \pbm$ is an elliptic,
selfadjoint, nonnegative operator of order $m>0$ with $\pi$ as a
principal angle.  Then the following statements are equivalent:
\newline (1) $A$ is of determinant class.
\newline (2) $\int_{0^{+}}^{1} \log \lambda dN_A(\lambda) > - \infty.$
\newline Here the integral $\int_{0^+}^1$ denotes the Stieltjes integral
on the half
open interval $(0,1].$}
\endproclaim

The proof of the proposition follows from the heat kernel
representation of the determinant which we briefly discuss (cf
\cite{Gi}).  Let $\gamma$ be a path in ${\Bbb C}$ defined by the
composition $\gamma_- \circ \gamma_+$ of two straight half lines:
$$\align
\gamma_+ & :=  \{x+ i(x + 1); \ -1 \leq x \leq \infty \} \\
\gamma_- & :=  \{x- i(x + 1); \ -1 \leq x \leq \infty \}
\endalign$$
where $\gamma_+$ starts at infinity and $\gamma_-$ starts at $x = -1.$
Using Proposition 2.9 and Lemma 2.10 we may define
the following bounded linear operator on $L_2({\Cal E}):$
$$\align
e^{-tA} & :=  \frac{1}{2\pi i} \int_\gamma e^{-t\lambda} (\lambda -
A)^{-1} d\lambda.
\endalign$$
One verifies that $e^{-tA} \in \pbmi$ for $t > 0.$  Hence, $e^{-tA}$ has
a smooth kernel, denoted $K_A(x,y,t),$ with values in ${\Cal B}$ and
admits a finite von Neumann trace, $\trn e^{-tA},$ given by
$$\align
\trn e^{-tA} & =  \int_{-\infty }^\infty e^{-t\lambda} dN_A(\lambda).\tag{2.15}
\endalign$$
As in the classical case one shows that for $t \longrightarrow 0,$ the
kernel $K_A(x,y,t)$ has an expansion on the diagonal $x=y$ of the form
$$\align
K(x,x,t) & = \sum_{j=0}^{N-1} t^{\frac{j-d}{m}}l_j(x) +
O(t^{\frac{N-d}{m}}) \tag {2.16}
\endalign$$
where $N \geq 1$ is arbitrary and the densities $l_j(x),$ when
expressed in local coordinates, are given by expressions of the form
(cf \cite{Gi})
$$\align
l_j(x) & =  \frac{1}{2\pi i} \frac{1}{(2\pi)^d} \int_{\rd} d\xi
\int_\gamma d\mu e^{-\mu} r_{-m-j}(x,\xi,\mu)
\tag {2.17}\endalign$$
where $r_{-m-j}(x,\xi,\mu),$ defined inductively by
(2.5), are the elements of the symbol expansion of the
resolvent $(\mu -A)^{-1}$ of $A$ in local coordinates.

{\sl Proof of Proposition 2.15.}  Using the heat kernel representation of
 the zeta-function we deduce that
$$\align
\ldet (A+\lambda) = & - \left. \frac{d}{ds}\right|_{s=0} \frac{1}{\Gamma(s)}
\int_0^1  t^{s-1} \left( \trn e^{-(A+\lambda )t} dt \right) \\
                    & - \int_1^\infty  t^{-1} \left( \trn
e^{-(A+\lambda )t}dt \right)  .
\endalign$$
The expansion (2.16) is used to show that
$$\align
        & - \left. \frac{d}{ds}\right|_{s=0} \frac{1}{\Gamma(s)}
\int_0^1  t^{s-1} \left(
\trn e^{-(A+\lambda )t} - \dimn ({\roman {Null}} (A))e^{-\lambda t}
\right) dt
\endalign$$
is a continuous function of $\lambda$ for $\lambda \geq 0.$  To
analyze
$$\align
G(\lambda) = & \int_1^\infty   t^{-1} \left(
\trn e^{-(A+\lambda )t} - \dimn ({\roman {Null}} (A))e^{-\lambda t}
\right) dt
\endalign$$
we write, applying Fubini's theorem together with $\trn e^{-At} = \int
_{-\infty}^\infty e^{-\mu t} dN_A(\mu)$ and
\newline$ \dimn ({\roman {Null}}(A))= N_A(0),$
$$\align
G(\lambda) = & \int_{0^+}^\infty dN_A(\mu) \int_1^\infty   t^{-1}
e^{-(\mu+\lambda )t} dt\\
           = &  \int_{\lambda^+}^\infty dN_{A+\lambda} (\mu) \int_1^\infty
 t^{-1} e^{-\mu t}dt .
\endalign$$
For $0 < \lambda \leq 1,$ write $G(\lambda) = G_1(\lambda) + G_2(\lambda)$
where $G_1(\lambda)$ and $G_2(\lambda)$ are given by
$$G_1(\lambda) = \int_{1^+}^\infty dN_{A+\lambda}(\mu) \int_1^\infty t^{-1}
e^{-\mu t}dt$$
$$ G_2(\lambda) = \int_{\lambda ^+}^1 dN_{A+\lambda}(\mu)\int_1^\infty t^{-1}
e^{-\mu t}dt.$$
The function $G_1(\lambda)$ is estimated in a straightforward way.
Concerning $G_2(\lambda)$, note that
$$\int_1^\infty t^{-1}e^{-\mu t} dt =
-\log \mu + (1 - e^{-\mu})\log \mu + \int_{\mu}^\infty e^{-s} \log s ds$$
and that the function $(1-e^{-\mu})\log \mu + \int_\mu ^\infty e^{-s} \log sds$
is bounded for $\mu \in [0,1].$ $~~~\diamondsuit$

\subheading{2.5 Elliptic boundary value problems}

Let $(M,g)$ be a compact Riemannian manifold with boundary $\partial M
\neq \emptyset.$  For the purpose of this paper we need only consider
the Dirichlet problem for an elliptic selfadjoint positive
differential operator $A$ of order 2, $A \in DO^2_{{\Cal B}}(M).$

Introduce the operator
$$\align
A_D   : & C_D^\infty({\Cal E}) \longrightarrow \ce
\endalign$$
where $ C_D^\infty({\Cal E}) :=\{u \in \ce : u|_{\partial M} = 0\}.$
Assume that $\pi$ is a principal angle for $A.$

Following \cite{Se2} one constructs a paramatrix, $R_D(\mu),$ for
$\mu -A_D$ in a similar fashion as in the case $\partial M =
\emptyset,$ describing inductively the asymptotic expansion of
the parametrix symbol.  The constructions differ in that, in the case of a
manifold with boundary, each term in the expected symbol expansion
includes a term arising from the boundary condition.  These added
terms arising from the boundary conditions only depend on the symbol
expansion of $A$ and its derivatives along the boundary $\partial M.$
Having constructed a parametrix one argues as in Proposition 2 of
section 2.4 to conclude that $A^{-1} \in \pbmm.$  This allows one to
introduce complex powers of $A_D$ and, for $\Re s >
\frac{d}{2},$ the zeta-function $\zeta_{A_{D}}(s)$ and its generalized
version $\zeta_{\alpha,A_D}(s) $ (cf section 2.4).  Following Seeley's
arguments one obtains the analog of Theorem 2.12

\proclaim{Theorem 2.12$'$} {Let $(M,g)$ be a compact Riemannian manifold with
boundary $\partial M \neq \emptyset.$  Assume that $A$ is a
selfadjoint, positive, differential operator of order 2 in $DO^2_{{\Cal
B}}(M)$ (note that $\pi$ is an Agmon angle for $A$ and therefore a
principal angle as well).  Then the function
$\zeta_{\alpha,A_D}(s) $ admits a meromorphic continuation to the
entire $s$-plane.
The continuation has at most simple poles and $s=0$
is a regular point.  The value of $\zeta_{\alpha,A_D}(s) $ at $s=0$ is
given by
$$\align
\zeta_{\alpha,A_D}(0)  = & \int_M \alpha(x) I_d(x) + \int_{\partial M}
\alpha(x) B_d(x)
\tag {2.18}\endalign$$
where in a coordinate chart of $(M,{\Cal E}\longrightarrow M), \
I_d(x)$ is defined as in (2.10).  In a coordinate chart
of $(\partial M, {\Cal E}|_{\partial M} \longrightarrow \partial M), \
B_d(x)$ is given by a formula similiar to that found in \cite{Se2}
involving at most the first $d$ terms of the symbol expansion of $A$
and its derivatives up to order $d.$}
\endproclaim
Theorem 2.12$'$ allows us to introduce the $\zeta$-regularized determinant
of $A_D$ by
$$\align
\ldet_N A_D  = & -\left. \frac{d}{ds}\right|_{s=0} \zeta_{A_D}(s).
\tag {2.19}\endalign$$
\vfil\eject

\heading
3. Asymptotic expansion and the Mayer-Vietoris type formula for
determinants
\endheading

\subheading{3.1 1-parameter families with parameter}

Let $\Lambda_{0,\epsilon}$ denote the solid angle in ${\Bbb C}$ given
by $\Lambda_{0,\epsilon} = \{re^{i\theta} ; r \geq 0, \ 2\pi |\theta|
\leq \epsilon\}.$  Consider a family of pseudodifferential operators
 $A(t), \ t \in \Lambda_{0,\epsilon}$ with $A(t) \in \pbm.$
\proclaim{Definition 3.1} A(t) is a 1-parameter family of weight
$\chi$ in $\pbm$ if for any chart $\phi: X \longrightarrow U$ of an
atlas of ${\Cal E} \longrightarrow M$ (where $X = \rd,$ or in case $U$
is a neighborhood of a boundary point, $X = \rd_+$) and for all $h, \
h' \in C^\infty_0(U),$ the operator $h'Ah,$ when expressed in local
coordinates, has an $\Cal L_{\Cal A}(\Cal W,\Cal W)$-valued symbol
$a = a_{h,h';U}$
satisfying the following properties:
\newline (1) for any multiindices $\alpha, \ \beta$ there is a constant
$C_{\alpha, \beta} > 0$ such that
$$\align
\| \partial_x^\alpha \partial_\xi^\beta a(x,\xi,t) \| \leq &
C_{\alpha \beta}(1+|\xi| + |t|^{\frac{1}{\chi}})^{m-|\beta|}
\endalign$$
where $x \in X, \ \xi \in \rd, $ and $ t \in \Lambda_{0,\epsilon};$
\newline (2) $a$ has an asymptotic expansion
$$\align
a \sim & \sum_{j\geq 0} \psi(\xi) a_{m-j}(x,\xi,t)
\tag {3.1}\endalign$$
with $\psi \in C^\infty (\rd)$ satisfying
$$
\psi(\xi) = \cases 0 & \hbox{ if } |\xi| \leq \frac{1}{2} \\
                   1 & \hbox{ if } |\xi| \geq 1 \endcases
$$ and $a_{m-j} \in C^\infty(X, \rd \setminus \{0\},
\Lambda_{0,\epsilon}; {\Cal L}_{{\Cal A}}({\Cal W},{\Cal W}))$
depending in a $C^1$-fashion on the parameter $t,$ has compact
$x$-support and is positive homogeneous of degree $m-j$ in $\xi, \
t^{\frac{1}{\chi}} ,$ i.e.
$$\align
a_{m-j}(x,\tau \xi, \tau^{\frac{1}{\chi}} t) = &
\tau^{m-j}a_{m-j}(x,\xi,t)
\endalign$$
for all $\tau > 0.$
\endproclaim
In the case where $M$ is closed one proves (cf e.g. \cite{Sh1}) that
for any $s \in {\Bbb R}$ and $l \geq m, \ A(t)$ is a bounded linear
operator, $A(t): \hse \longrightarrow H_{s-l}({\Cal E}).$  Denote by
$|\| A(t) \| |_{s \rightarrow s-l}$ the operator norm of $A(t)$, viewed
as an operator $A(t): \hse \longrightarrow H_{s-l}({\Cal E}).$

\proclaim{Theorem 3.2}{ Let $M$ be closed.  The following estimates
hold:
\newline (1) if $l \geq 0,$ then $|\| A(t) \| |_{s \rightarrow s-l} \leq
C_{s,l}\left(1 + |t|^{\frac{1}{\chi}}\right)^m ;$
\newline (2) if $m \leq l\leq 0,$ then  $|\| A(t) \| |_{s \rightarrow s-l} \leq
C_{s,l}\left(1 + |t|^{\frac{1}{\chi}}\right)^{-(l-m)}. $
}
\endproclaim

\proclaim{Definition 3.3} A 1-parameter family A(t) in $\pbm$ is
elliptic with parameter, if for any chart $\phi: X \longrightarrow U$ of an
atlas of ${\Cal E} \longrightarrow M$ (where $X = \rd,$ or in case $U$
is a neighborhood of a boundary point, $X = \rd_+$) and for all $h, \
h' \in C^\infty_0(U),$ the operator $h'Ah,$ when expressed in local
coordinates, has principal symbol $a_m(x,\xi,t)$ with values in $ {\Cal
L}_{{\Cal A}}({\Cal W},{\Cal W})$ such that for all $x \in X$ with
$h(\phi(x)) h'(\phi(x)) \neq 0, \ a_m(x,\xi,t)$ is invertible for
$(\xi , t) \in (\rd \times \Lambda_{0,\epsilon}) \setminus \{(0,0)\}.$
\endproclaim

Let $M$ be closed.  For a 1-parameter family $A(t),$ elliptic with parameter,
one constructs a parametrix, $R(\mu,t),$ for $\mu - A(t):$  Given $\mu
\notin \bigcup_{t \in \Lambda_{0.\epsilon}} \spec (A(t)), \ R(\mu,t)$ is
a 1-parameter family in $\pbmm$ satisfying
$$\align
R(\mu,t)(\mu - A(t)) - {\roman Id} \in & \pbmi ~~~ \text{and}\\
(\mu - A(t))R(\mu,t) - {\roman Id} \in & \pbmi .
\endalign$$

In local coordinates the symbol of $R(\mu,t)$ is constructed
inductively:
$$\align
r_{-m}(x,\xi,t,\mu) = & (\mu - a_m(x,\xi,t))^{-1}
\endalign$$
$$\multline
r_{-m-j}(x,\xi,t,\mu) \\
= r_{-m}(x,\xi,t,\mu)\left(\sum_{k=0}^{j-1} \sum_{|\alpha| +l + k = j}
\frac{1}{\alpha !}  \partial_\xi^\alpha a_{m-l}(x,\xi,t) D_x^\alpha
r_{-m-k}(x,\xi,t,\mu)\right)
\endmultline \tag {3.2}$$
where $D_x = \frac{1}{i}\partial_x.$  The term $r_{-m-j}(x,\xi,t,\mu)$
is positive homogeneous of degree $-m-j$ in $(\xi, \
t^{\frac{1}{\chi}}, \ \mu^{\frac{1}{m}}).$
\subheading{3.2 Asymptotic expansion for determinants}

As in \cite{BFK2, Appendix}, one proves a result concerning the asymptotic
expansion of a 1-parameter family $A(t)$ in $\pbm, \ A(t)$ elliptic
with parameter.
\proclaim{Theorem 3.4}{Let $M$ be a closed manifold.  Assume that
$A(t)$ is a 1-paramter family in $\pbm,$ elliptic with paramter of
weight $\chi$ and having $\pi$ as an Agmon angle
uniformly in $t$ (cf \cite{BFK1, Theorem 1.1}).  Then the function
$\ldet_N A(t)$ admits an asymptotic expansion for $t \longrightarrow
\infty$ of the form
$$\align
\ldet_N A(t) \sim & \sum_{-\infty}^d \overline{a}_j |t|^{\frac{j}{\chi}}
+  \sum_{0}^d \overline{b}_j |t|^{\frac{j}{\chi}}\log |t|
\tag {3.3}\endalign$$
where $\overline{a}_j = \int_M a_j(x,\frac{t}{|t|}) dx, \
\overline{b}_j = \int_M b_j(x,\frac{t}{|t|}) dx, $ are defined by
smooth densities $ a_j(x,\frac{t}{|t|})$ and $
b_j(x,\frac{t}{|t|})$, which can be computed in terms of the
symbol of $A(t).$}
\endproclaim

In particular, with respect to a coordinate chart, $
a_0(x,\frac{t}{|t|})$ is given by
$$\align
 a_0\left(x,\frac{t}{|t|}\right) = & \left. \frac{d}{ds}\right|_{s=0}
 \frac{1}{(2\pi)^d}
\frac{1}{2\pi i} \int_{\rd} d\xi \int_\Gamma d\mu \mu^{-s} \trn\left(
r_{-m-d}\left(x,\xi, \frac{t}{|t|},\mu\right)\right) \\
                       = & - \frac{1}{(2\pi)^d} \int_{\rd} d\xi
\int_0^\infty d\mu \trn \left(r_{-m-d}\left(x,\xi, \frac{t}{|t|},-\mu\right)
\right).
\tag {3.4}\endalign$$
A similar result holds in the case where $M$ has a nonempty boundary,
$\partial M \neq \emptyset$ (cf \cite{BFK2}).  With the notation
introduced in section 2.6, one obtains

\proclaim{Theorem 3.5}{Let $(M,g)$ be a compact Riemannian  manifold
with boundary, $\partial M \neq \emptyset.$  Assume that
$A(t), \ t \in \Lambda_{0,\epsilon}$, is a 1-parameter family
of selfadjoint, positive differential operators in $\pbm,$
of order $m=2$, elliptic with paramter of
weight $\chi.$ Assume
that there exists $\epsilon' >0$ such that for all $t \in
\Lambda_{0,\epsilon'}, \ \spec A_D(t) \cap V_{\pi,\epsilon'} =
\emptyset.$  Then the function
$\ldet_N A_D(t)$ admits an asymptotic expansion for $t \longrightarrow
\infty$ of the form
$$\align
\ldet_N A_D(t) \sim & \sum_{j=-\infty}^d
\left(\overline{a}_j+\overline{a}_j^b\right)
|t|^{\frac{j}{\chi}} +
\sum_{j=0}^d \left(\overline{b}_j+\overline{b}_j^b\right)
|t|^{\frac{j}{\chi}}\log |t| \tag{3.5}
\endalign$$
where $\overline{a}_j$ and $\overline{b}_j$ are given as in Theorem
3.4. The quantities $\overline{a}_j^b$ and $\overline{b}_j^b$ are
contributions from the boundary and are of the form
$$
\overline{a}_j^b = \int_{\partial M}  a_j^b\left(x,\frac{t}{|t|}\right);
\ \ \ \overline{b}_j^b = \int_{\partial M}  b_j^b\left(x,\frac{t}{|t|}\right).
\tag {3.6}$$
In a coordinate chart of $(\partial M, {\Cal E}|_{\partial M}
\longrightarrow  \partial M)$ the densities $
a_j^b(x,\frac{t}{|t|})$ and $ b_j^b(x,\frac{t}{|t|})$ are given by a
formula each involving only finitely many terms in the symbol
expansion of $A(t)$ and finitely many of its derivatives.}
\endproclaim

\subheading{3.3 Mayer-Vietoris type formula}

We restrict ourselves to the case needed for this paper.  We assume
throughout this subsection that $(M,g)$ is a closed Riemannian
manifold.  Let $\Gamma$ be a smooth embedded hypersurface in $M$
with trivial normal bundle.
Consider an elliptic, selfadjoint, positive, differential operator $A$ of
order $2$, $A: \ce \longrightarrow \ce$, of Laplace -Beltrami type (i.e
the principal symbol is of the form $\sigma_{A}(x,\xi)=\|\xi\|^2Id_{\Cal E_x}$)
with $\spec (A) \subset
[\epsilon , \infty)$ for some $\epsilon >0.$  Denote by
$M_\Gamma$ the manifold whose interior is
 $M \setminus \Gamma $, and
whose boundary is $\partial M_\Gamma = \Gamma^+ \sqcup \Gamma^-,$
where $\Gamma^+$ and $\Gamma^-$ are both copies of $\Gamma.$  Let
$g_\Gamma$ be the Riemannian metric on $M_\Gamma$ obtained by pulling
back the metric $g$ and let ${\Cal E}_\Gamma \longrightarrow M_\Gamma$
be the pullback of the bundle ${\Cal E} \longrightarrow M.$  Consider
$A_\Gamma : \ceg \longrightarrow  \ceg$ with Dirichlet boundary
conditions.  Then $A_\Gamma$ is selfadjoint, positive and elliptic
with $\spec (A_\Gamma) \subset [\epsilon,\infty)$  so that $\pi$
is an Agmon angle for $A_\Gamma.$
Introduce the Dirichlet to Neumann operator, $R_{DN},$ associated to
the unit vectorfield normal to $\Gamma.$  This operator is defined as the
composition
$$\align
\ceoverg \overset\Delta_{{\roman {iag}}} \to\longrightarrow  & \cegp
\oplus \cegm
 \overset P_D \to \longrightarrow \ceoverg \\
          &  \overset N \to \longrightarrow \cegp \oplus \cegm
\overset \Delta_{{\roman {iff}}} \to \longrightarrow \ceoverg
\endalign$$
where $\Delta_{{\roman {iag}}}(f) = (f,f)$ is the diagonal operator, $P_D$
is the
Poisson operator associated to $A_\Gamma,$ $N$ is the first order
scalar differential operator induced by the normal unit vectorfield along
$\Gamma$ and $\Delta_{{\roman {iff}}}$ is the difference operator
$\Delta_{{\roman {iff}}}(f^+,f^-) = f^+ - f^-.$  As in \cite{BFK2}
one proves the following
\proclaim{Theorem 3.6}{Assume that $(M,g)$ is a closed Riemannian
manifold and $A$ is an elliptic, selfadjoint, positive differential
operator, $A: \ce \longrightarrow \ce$ of order 2 of Laplace-Beltrami type
with
 $\spec (A) \subset [\epsilon, \infty)$ for some
$\epsilon >0.$  Then $R_{DN}$ is an invertible pseudodifferential
operator in $\Psi DO^1_{{\Cal B}}(\Gamma).$
The inverse $R_{DN}^{-1}$ is given by
$$\align
R_{DN}^{-1} = & JA^{-1}( \cdot \otimes \delta_\Gamma)
\tag {3.7}\endalign$$
where $J$ is the trace operator $J: \hse \longrightarrow H_{s-1}({\Cal
E}|_\Gamma) $ and $\delta_\Gamma$ denotes the Dirac distribution along
$\Gamma$ (cf \cite{BFK2} (4.5)).  As a consequence one concludes
\newline (1) $R_{DN}$ is selfadjoint and positive with $\spec (R_{DN})
\subset [\epsilon',\infty)$ for some $\epsilon' >0.$  In particular,
$\pi$ is an Agmon angle for $R_{DN}.$
\newline (2) The principal symbol, $ \sigma(R_{DN}^{-1})$, of $R_{DN}^{-1}$ can
be computed in terms of the principal symbol $\sigma(A^{-1})$ of
$A^{-1}$ (cf \cite{BFK2} (4.6)):
$$\align
 \sigma(R_{DN}^{-1})(x',\xi') = & \frac{1}{2\pi} \int_{{\Bbb
R}}\sigma(A^{-1})(x',0,\xi',\eta) d\eta
\tag {3.8}\endalign$$
where $x=(x',w)$ are coordinates in a collar neighborhood of $\Gamma$
such that $x' $ are coordinates of $\Gamma$ and the normal vectorfield
along $\Gamma$ is represented by $\frac{\partial }{\partial w}.$
\newline (3) In a coordinate chart for $\Gamma$ which arises from a chart
belonging to an atlas of ${\Cal E}|_\Gamma \longrightarrow \Gamma,$
the symbol of $R_{DN}$ has an expression whose terms depend only on
the terms of the expansion of the symbol of $A$ and its derivatives in
an arbitrarily small neighborhood of $\Gamma$.
\newline (4)
$$\align
{\roman det}_N(A) = & \overline{c}  {\roman det}_N(A_\Gamma)
 {\roman det}_N(R_{DN})
\endalign$$
where
$$\align
\overline{c}  = & \exp \left(\int_\Gamma c(x) \right)
\endalign$$
and the density $c(x),$ when expressed in a coordinate chart of
$\Gamma$ which is contained in an atlas of ${\Cal E}|_\Gamma
\longrightarrow \Gamma,$ depends only on the first $d$ terms of the
symbol expansion of $A$ and their derivatives in an arbitrarily small
neighborhood of $\Gamma.$
\newline (5) Assume that instead  of the single operator $A,$ there is a
family $A(t) : \ce \longrightarrow \ce $ of differential operators of
order 2 of Laplace-Beltrami type with parameter $t \in \Lambda_{0,\epsilon'},
\ \epsilon'>0,$
of weight $\chi$ such that $A(t)$ is elliptic, selfadjoint and
positive for each $t.$  Introduce as above $A(t)_\Gamma, \ R_{DN}(t)$
and assume that $\spec (A(t)) \cap V_{\pi, \epsilon'} = \emptyset$ for
some $\epsilon >0$ and for all $t \in \Lambda_{0,\epsilon}.$  Then
$R_{DN}(t)$ is an invertible family of pseudodifferential operators
with parameter (cf \cite{BFK2} (3.13)) of order 1 and weight $\chi.$
}
\endproclaim
{\bf Remark}
For the convenience of the reader who is only interested
in the results as stated above, Y.Lee has written an easily
accessible version of [BFK2] (cf [Lee]).

\vfil\eject

\heading
4. Torsions and Witten deformation of the analytic torsion
\endheading

\subheading{4.1 Reidemeister and analytic torsion in the von
Neumann sense}

Let $M$ be a closed manifold of dimension $d$
and $\Cal W$ an $(\Cal A,\Gamma^{op})$-Hilbert module of finite type with
$\Gamma=\pi_1(M)$
the fundamental group of $M$. Let $p:\Cal E \to M$ be the
bundle of $\Cal A$-Hilbert modules over $M$ associated to $\Cal W$
as described in section 1. The fiber of this bundle is isomorphic to
the  $\Cal A$-Hilbert module $\Cal W.$  The smooth bundle $p:\Cal E\to
M$ is equipped with a flat {\it canonical} connection.
Both its
Hermitian structure $\mu$ and fiberwise $\Cal A$-action $\rho$ are  left
invariant by the
parallel transport induced by
the canonical connection.

Let $h:M\to \Bbb R$ be a smooth Morse function. For convenience we assume
that $h$ is
self-indexing, i.e. $h(x)=\text{index}  (x)$ for any critical point $x$ of $h$.
Let $g'$ be a
Riemannian metric so that $\tau = (h,g')$ is a generalized
triangulation.
This means that for any two critical points $x$
and $y$ of $h$, the unstable
manifold $W^-_x$
and the stable  manifold $W^+_y$, associated to the vector field
${\roman {grad}}_{g'}h$, intersect transversally and, in a neighborhood of any
critical point $x$ of $h$, there exist coordinates $y_1,...,y_d,$
with respect to which $h$ is of the form $h(x)= k - (y_1^2 + ...+y_k^2)/2 +
(y_{k+1}^2+...+y_d^2)/2$
with $k = \text{index}(x)$ and the metric $g$ is Euclidean (cf. Introduction).
Let $\tilde M\to M$ be the universal covering of $M$
and $\tilde h$ and  $\tilde g'$ be the lifts of $h$ and $g'$ on $\tilde M$.
Denote by ${\roman {Cr}}_q(h) \subset M$ resp. ${\roman {Cr}}_q(\tilde h)
\subset \tilde M$ the
set of critical  points of index
$q$ of $h$ resp. $\tilde h$ and let ${\roman {Cr}}(\tilde h) =
\cup_q {\roman {Cr}}_q(\tilde h).$
Clearly the group $\Gamma$ acts freely on ${\roman {Cr}}_q(\tilde h),$
for any $q,$
and the quotient
set can be identified with ${\roman {Cr}}_q(h).$

For each $\tilde x \in {\roman {Cr}}(\tilde h)$
choose an orientation $O_{\tilde x}$ for the unstable manifold $W^-_{\tilde
x}$  and denote $O_h := \{O_{\tilde x};
\tilde x \in {\roman {Cr}}(\tilde h)\}$.
To the quadruple $(M,\tau,O_h, {\Cal W})$ we
associate a cochain complex of finite type over the von Neumann algebra $\Cal
A$, $\Cal C(M,\tau,O_h)=\{\Cal C^q, \delta_q\}.$ The components $\Cal C^q$ are
the $\Cal A$-Hilbert module of finite type,
$\Cal C^q:=\Gamma(\Cal E|_{{\roman {Cr}}_q(h)})
=\bigoplus_{x\in {\roman {Cr}}_q(h)} \Cal E_x$
which can be identified with the module of $\Gamma$-equivariant
maps $f:{\roman {Cr}}_q(\tilde h)\to \Cal W.$
To define the maps $\delta_q$ a few remarks are in order.
The orientation on $\tilde M$ together with the
orientations $O_h$
induce orientations on the stable
manifolds $W^+_{\tilde x}$
which in turn  permit us to define the following functions
$m_q:{\roman {Cr}}_q(\tilde h)\times
{\roman {Cr}}_{q-1}(\tilde h)\to \Bbb Z$
$$
m_q(\tilde x,\tilde y):=\text{intersection number }(W^-_{\tilde x},
W^+_{\tilde y}).
$$
Notice that the functions $m_q$ have the following properties:
\newline (In1) $m_q(\tilde x, \tilde y)=m_q(g\tilde x, g\tilde y)$, for all
 $g\in \pi_1(M);$
\newline (In2)
$\{\tilde x\in {\roman {Cr}}_q(\tilde h); m_q(\tilde x, \tilde y) \ne 0 \}$ is
finite  for any $\tilde y\in {\roman {Cr}}_{q-1}(\tilde h);$
\newline (In3)
$\{\tilde y\in {\roman {Cr}}_{q-1}(\tilde h); m_q(\tilde x, \tilde y) \ne 0\}$
is finite for any $\tilde x\in {\roman {Cr}}_q(\tilde h);$
\newline (In4)
$\sum_{\tilde y\in {\roman {Cr}}_{q-1}(\tilde h)} m_q(\tilde x,\tilde y)\cdot
m_{q-1}(\tilde y,\tilde z)=0$ for any
 $\tilde x\in {\roman {Cr}}_q(\tilde h)$ and any
 $\tilde z\in {\roman {Cr}}_{q-2}(\tilde
h).$

Properties (In1)-(In3) imply that for any $\Gamma$-equivariant map
$f: {\roman {Cr}}_{q-1}(\tilde h) \to \Cal W,$
we can define the $\Gamma$-equivariant map
$\delta_{q-1}(f): {\roman {Cr}}_q(\tilde h) \to \Cal W$ by the formula
$$
\delta_{q-1}(f)(\tilde x)=
\sum_{\tilde y\in {\roman {Cr}}_{q-1}(\tilde h)}m_q(\tilde x,\tilde
y)f(\tilde y).\tag{4.1}
$$
By property (In4), $\delta_q\cdot \delta_{q-1}=0.$

One defines $\log T_{{\roman {comb}}}(M,\tau)\in \bold D$ by
$$\log T_{{\roman {comb}}}(M,\tau)
:=\log T(\Cal C(M,\tau,O_h,\Cal W)) \tag {4.2}$$ (cf section 1).

One can show that $\log T_{{\roman {comb}}}$ is independent of
the choice of the orientations $O_h.$

Let $(M,g)$ be a Riemannian manifold and $\Cal W$ a $(\Cal
A,\Gamma^{op})$-Hilbert module of finite type.  Let $\Lambda^q(M;\Cal
E)=C^\infty(\Cal E\otimes \Lambda^q(T^*M))$ be the space of
smooth $q$-forms with values in $\Cal W$ where $T^*M$ denotes the
cotangent bundle of $M$ and $p: \Cal E \to M$ is a smooth bundle of
$\Cal A$-Hilbert modules of finite type with fiber $\Cal W$.  The
Riemannian metric $g$ induces the Hodge operators
$R_q:\Lambda^q(T^*M)_x\to\Lambda^{d-q}(T^*M)_x~~$ $(x\in M)$ and the
Hermitian structure $\mu$ on $\Cal E$ together with the Hodge operators
induce a Hermitian structure on $\Cal E\otimes\Lambda^q(T^*M)$ given by
$(s,s' \in C^{\infty}(\Cal E); w,w' \in C^{\infty}(\Lambda^q(T^*M)) )$
$$
(s\otimes w, s'\otimes w')(x)=\mu_x(s(x),s'(x))R_q(w(x)\wedge R_qw'(x)).
$$
As a consequence $\Cal E\otimes \Lambda^q(T^*M)$ is smooth bundle
of $\Cal A-$Hilbert modules. The canonical  connection in $p:\Cal E\to
M$ can  be interpreted as a first order differential operator $_{\Cal
W} d_q:\Lambda^q(M;\Cal  E) \to \Lambda^{q+1}(M;\Cal E).$
As the canonical connection is flat,
$_{\Cal W} d_{q+1}._{\Cal W}d_q=0$ for any $q.$  Notice that $_{\Cal
W}d_q$ is an $\Cal A$-linear, differential operator and if the action
of $\Gamma$ on $\Cal W$ is trivial, $_{\Cal W}d$ is the usual exterior
differential ${\roman {Id}}\otimes d.$ In case there is no risk of ambiguity we
will write $d$ instead of $_{\Cal W}d$ and continue to call it
exterior differential.

The formal adjoint of $_{\Cal W} d_q$ with respect to the above
defined Hermitian structure is a first order differential operator
$_{\Cal W} d_q^*:\Lambda^{q+1}(M;\Cal E) \to \Lambda^q(M;\Cal E)$ and is
again $\Cal A$-linear. Introduce the  Laplacians, acting on
$q$-forms,
$$\Delta_q=d_q^*d_q+d_{q-1}d_{q-1}^*.$$
The operators  $\Delta_q$ are essentially selfadjoint, nonnegative,
elliptic and $\Cal A$-linear. The space $\Lambda^q(M;\Cal E)$
can be equipped
with the scalar product
$$
\langle u_1,u_2\rangle_{r}~ =~ \langle({\roman {Id}}+
\Delta_q)^{r/2}(u_1),({\roman {Id}}+\Delta_q)^{r/2}(u_2)\rangle
\tag {4.3}$$
where
$$
\langle({\roman {Id}}+\Delta_q)^{r/2}(u_1),({\roman
{Id}}+\Delta_q)^{r/2}(u_2)\rangle ~$$
$$ = ~ \int_M (({\roman {Id}}+\Delta_q)^{r/2}(u_1),({\roman
{Id}}+\Delta_q)^{r/2}(u_2))(x)d{\roman {vol}}_g.
$$
The completion of $\Lambda ^q(M,\Cal E)$
with respect to the scalar product $\langle.,.\rangle_{r}$ is an
$\Cal A$-Hilbert module
\newline $H_{r}(\Lambda^q(M;\Cal E))$, the space of forms of degree
$q$ in Sobolev space of order $r$.
In the case where $r=0$, we write also $L_2(\Lambda^q(M;\Cal E)).$
Obviously, these Hilbert modules are not of finite type.
Note that the operators $({\roman {Id}}+\Delta_q)^{r/2}$
define isometries  between $H_{r'}(\Lambda ^q(M;\Cal E))$ and
 $H_{(r'-r)}(\Lambda ^q(M;\Cal E)).$ Let $\Cal H_q$ be the $\Cal A-$Hilbert
module of harmonic q-forms
 $$\Cal H_q=
\{\omega\in L_2(\Lambda^q(M;\Cal E)); \Delta_q(\omega)=0\}.$$
Since $\Delta_q$ is elliptic, $\Cal H_q \subset \Lambda^q(M;\Cal E).$
The integration $Int^{(q)}$ on the q-cells of the generalized
triangulation $\tau$, which are given by the unstable manifolds of
 ${\roman {grad}}_{g'}h$, defines an $\Cal A-$
linear map
$$
Int^{(q)}:\Lambda^{q}(M;\Cal E) \to  \Cal C^q
$$
so that $\delta_q Int^{(q)} = Int^{(q)} d_q.$
Denote by $\pi _q$ the canonical projection
$\pi_q : \Cal C^q \to {\roman {Null}} (\Delta_q^{{\roman {comb}}}).$
By a theorem of Dodziuk $[Do]$ of de-Rham type, the map $\pi_q Int^{(q)},$
restricted to $\Cal H_q,$ is an
isomorphism of Hilbert modules. Denote its inverse by $\theta_q.$  Since
${\roman {Null}} (\Delta_q^{{\roman {comb}}})$ is an $\Cal A-$Hilbert
 module of finite type so is $\Cal  H_q.$
Define $T_{{\roman {met}}}$ as the positive real number, viewed as an element
 in $\bold D$ (cf. Introduction) by
$$
\log T_{{\roman {met}}}(M,g,\Cal W,\tau):=\frac {1}{2}
 \sum_q (-1)^q {\roman {logdet}}_N( \theta^*_q \theta_q)
.\tag {4.4}$$

The Reidemeister torsion $T_{{\roman {Re}}}(M,g,\Cal W,\tau)\in \bold D$ is
defined
(cf [CM],[LR]) by
$$
\log T_{{\roman {Re}}}(M,g,\Cal W,\tau)=
 \log T_{{\roman {comb}}}(M,\Cal W,\tau)+
 \log T_{{\roman {met}}}(M,g,\Cal W,\tau) \tag {4.5} $$
and the analytic torsion
$T_{{\roman {an}}}(M,g,\Cal W)\in \bold D$ (cf [Lo],[Ma]) by
$$
\log T_{{\roman {an}}}(M,g,\Cal W)=
\frac {1}{2} \sum_q (-1)^{q+1} q{\roman {logdet}}_N(\Delta_q).
\tag {4.6}$$

Following Gromov-Shubin \cite{GS}, for $\lambda \geq 0$, we introduce
the functions  $F_q(\lambda):= F_{{\underline d}_q}(\lambda) =
\sup \{ {\roman {dim}}_N {\Cal L} ; {\Cal L} \in \Cal P_q (\lambda)
\}$ where $\Cal P_q(\lambda)$
consists of all $\Cal A$-invariant closed subspaces
${\Cal L} \subset \overline{d_{q-1}(\Lambda^{q-1}(M;\Cal E))}
 \subset L_2(\Lambda^q(M;\Cal E)),$
so that for any $\omega \in {\Cal L}, \ \omega$ is in the domain of
definition of $d_q$ and
$$\| d_q \omega \| \leq \lambda^{1/2} \|\omega \|.\tag {4.7}$$
Note that a subspace ${\Cal L}$ satisfying (4.7) is in fact contained in
$\Lambda^{q,+}(M;\Cal E)$ where
$$\Lambda^{q,+}(M;\Cal E)=
\overline{d_{q-1}(\Lambda^{q-1}(M;\Cal E))}\cap \Lambda^q(M;\Cal E).
\tag {4.8}$$

These functions are elements in the space $\bold F$ (cf section 1). By
arguments of Gromov-Shubin  which we recalled in section 1.2, the
spectral functions $N_k(\lambda) = N_{\Delta_k}(\lambda)$ of the Laplace
operator $\Delta_k$ are given by
 $\beta_k + F_{k-1}(\lambda) + F_k(\lambda).$

\proclaim {\bf Definition 4.1}
\newline (1) The system $(M,\tau,\Cal W)$ is said to be of
$c-determinant$ class iff for $0 \leq k \leq d,$
$$\int_{0^+}^1 \log \lambda dN_{\Delta_k^{{\roman {comb}}}}(\lambda)
> - \infty.$$
\newline (2) The system $(M,g,\Cal W)$ is said to be of
$a-determinant$ class iff for $0 \leq k \leq d,$
$$\int _{0+}^1 \log \lambda dN_{\Delta_k} (\lambda) > - \infty $$
or, equivalently,
$$ \int_{0}^1 \log \lambda dF_k(\lambda) > - \infty.$$
\endproclaim

\subheading{4.2 Witten's deformation of the analytic torsion}

Let $\omega\in \Lambda^1(M)$ be a smooth closed 1-form on $M$.
Introduce a perturbation
\newline $(\Lambda^q(M;\Cal E),_{\Cal W}d^{\omega}_q)$ of the deRham complex
 $(\Lambda^q(M;\Cal E),_{\Cal W}d_q)$ with
$$d^{\omega}_q :=_{\Cal W}d^{\omega}_q:=_{\Cal W}d_q+\omega\wedge (.).
\tag {4.9}$$
The formal adjoint of $d^{\omega}_q$ with respect to the
Hermitian structure on $\Cal E\otimes\Lambda^q(T^*M),$ introduced in
section 4.1, is a first
order
$\Cal A$-linear, differential operator
$$(d^{\omega}_q)^*:\Lambda^{q+1}(M;\Cal
E) \to \Lambda^q(M;\Cal E).$$
Introduce the perturbed Laplacians, acting on $q$-forms,
$$
\Delta^{\omega}_q=(d^{\omega}_q)^*d^{\omega}_q+
d^{\omega}_{q-1}(d^{\omega}_{q-1})^*.
\tag {4.10}$$
The operators $\Delta^\omega_q$ are $\Cal A$-linear, elliptic
operators which are positive and essentially selfadjoint.
They are zero'th order perturbations of the Laplacians $\Delta_q$
defined above.
The case $\omega = tdh$ where
$h:M\to \Bbb R$ is a smooth function
was considered  by Witten $cf [Wi]$.
The multiplication by $e^{th}$ defines, for any r, a linear operator on
 $H_r(\Lambda^q(M;\Cal E))$, which is an isomorphism of
 $\Cal A-$Hilbert modules and we have $d_q(t) = e^{-th} d_q e^{th}.$
 We call the operators
 $\Delta_q(t) = \Delta_q^{tdh}$ the Witten Laplacians associated to $h$.
 More generally, we will refer to the complex
$(\Lambda^q(M;\Cal E), d^{\omega}_q(t))$
with $d^{\omega}_q(t)=d^{t\omega}_q$ as the
Witten complex.
Define the perturbed analytic torsion $T_{{\roman {an}}}(M,g,\Cal W,\omega)$
as an
element in the vector space $\bold D$
$$
\log T_{{\roman {an}}}(M,g,\Cal W,\omega) := \frac {1}{2} \sum_q (-1)^q
q{\roman {logdet}}_N(\Delta^{\omega}_q)
$$
and the Witten deformation of the analytic torsion
$T_{{\roman {an}}}(M,g,\Cal W,\omega)$
$$\log T_{{\roman {an}}} (M,g,\Cal W,\omega)(t):=
\log T_{{\roman {an}}} (M,g,\Cal W,t\omega).\tag {4.11}$$

{\bf Remark}
If $(M,g, \Cal W)$ is of $a-determinant$ class and $\omega = dh$ then
the Witten deformation satisfies
$\log T_{{\roman {an}}}(M,g,\Cal W, t\omega) \in \Bbb R \subset
\bold D$, for any $t.$
This can be verified as follows: define functions
 $F_{{\underline d}_k^{tdh}}(\lambda)$ as above replacing $d_k$ by $d_k^{tdh}.$
As $(L_2(\Lambda^k(M;\Cal E)),d_k)$ and $(L_2(\Lambda^k(M;\Cal E)), d_k(t))$
are isomorphic, one concludes, according to results of Gromov-Shubin,
that $F_{{\underline d}_k^{th}}(\lambda) \overset {d}\to
{\sim} F_{{\underline d}_k}(\lambda)$
and therefore, as $\Delta_k$ is of determinant class, so is $\Delta_k(t)$.

\subheading{4.3 Product formulas}

For $i=1,2,$ let $\Cal A_i$ be finite von Neumann algebras, $(M_i,g_i,\tau_i)$
closed
Riemannian manifolds of dimension $d_i$ (even or odd),
equipped with the generalized triangulations $\tau_i=(h_i,g'_i).$ Let
$\Cal W_i$ be ($\Cal A_i,\Gamma_i^{op}$)-Hilbert
modules of finite type $\Gamma_i=\pi_1(M_i)$, and  $\omega_i\in
\Lambda^1(M_i)$  closed 1-forms $(i=1,2).$
Introduce
$\Cal A:=\Cal A_1\otimes\Cal A_2$,
$\Cal W:=\Cal W_1\otimes\Cal W_2$,  $M=M_1\times M_2$,  $g=g_1\times g_2$,
$\tau:=(h=h_1+h_2, g'=g'_1\times g'_2)$ and $\omega=\omega_1\otimes 1+
1\otimes\omega_2.$ Further denote by $\Cal E \to M$ and $\Cal E_i \to M_i$
($i=1,2$) the bundles associated to $\Cal W_i$ and $\Cal W$.

\proclaim{\bf Proposition 4.1(Product formula) } (cf [CM],[Lo],[LR])
{ With the hypotheses above:
\roster
\item $$\multline
\log T_{{\roman {an}}}(M,g,\Cal W,\omega) \\
 =  \chi (M_1)\cdot\log T_{{\roman {an}}}(M_2,g_2,\Cal
W_2,\omega_2)+\chi(M_2)\cdot\log T_{{\roman {an}}}(M_1,g_1,\Cal
W_1,\omega_1)  \endmultline \tag {4.12}$$

\item $$\multline
\log T_{{\roman {Re}}}(M,g,\Cal W,\tau) \\
 =  \chi(M_1)\cdot \log
T_{{\roman {Re}}}(M_2,g_2,\Cal W_2,\tau_2)+\chi(M_2)\cdot
\log T_{{\roman {Re}}}(M_1,g_1,\Cal W_1,\tau_1)
\endmultline  \tag {4.13}$$
\endroster}
\endproclaim
{\sl Proof:} (2) follows from Corolary 1.22 and Proposition 1.9.
To prove (1) observe that
$$L_2(\Lambda^r(M,\Cal E))=\oplus_{p+q=r}L_2(\Lambda^p(M_1,
\Cal E_1))\otimes L_2(\Lambda^q(M_2,\Cal E_2))$$
and note that $\Delta_q=\oplus_{p+r=q}\Delta_{(p,r)}$ with
$$\Delta_{(p,r)}=(\Delta'_p\otimes {\roman {Id}}) +
({\roman {Id}}\otimes\Delta_r'')$$
is an $\Cal A$-linear, elliptic, differential operator where $\Delta_p'$
and $\Delta_r''$ denote the Laplacians corresponding to $\Cal E_1 \to
M_1,$ respectively, $\Cal E_2 \to M_2.$  Notice that
$e^{-t\Delta_{p,r}}=e^{-t\Delta_p'}\otimes e^{-t
\Delta_r''}$ is of trace class in the von Neumann sense.  As in
(1.30) introduce  $$ \zeta_M(\lambda ,s) = \frac {1}{2} \sum_{q\ge
1}(-1)^qq \frac{1}{\Gamma(s)}
\int_0^{\infty} t^{s-1}{\roman {tr}}_N e^{-t(\Delta_q+\lambda)}dt.\tag {4.14}$$
As a consequence of a theorem of deRham type due to Dodziuk [Do] one obtains
$\chi (M)=\sum (-1)^q {\roman {dim}}_N({\overline {\Cal H}_q(M;\Cal E)})$.
To prove (1) it suffices to verify that
$$\align
\zeta_M(\lambda,s)= & \zeta_{M_1}(\lambda,s)\cdot \chi(M_2)+
\zeta_{M_2}(\lambda,s)\cdot \chi(M_1).\tag {4.15}
\endalign$$
In order to apply the line of arguments of the proof of Proposition 1.21  one
needs only to prove that
$$ {\roman {tr}}_N e^{-t\Delta^+_{q+1}} = {\roman {tr}}_N e^{-t\Delta^-_q}
\tag {4.16}$$
where $\Delta_q^+$ respectively $\Delta_q^+$ denote the restriction of
$\Delta_q$ to $\Lambda^{q,+}(M;\Cal E)$ respectively
\newline $\Lambda^{q,-}(M;\Cal E).$
Equation (4.16) follows from the observation that the spectral
projector $P^+_{q+1}(\lambda)$ and $P^-_q(\lambda)$ associated to
$\Delta^+_{q+1}$ respectively $\Delta_q^-$ are intertwined by $d_q$ and
therefore ${\roman {tr}}_N P^+_{q+1}(\lambda) =
{\roman {tr}}_N P_q^-(\lambda).$
 $\diamondsuit$
\vfil\eject

\subheading{5. Witten's deformation of the deRham complex}

In this section we discuss Witten's deformation of the deRham complex of
differential forms with coeffiecients in a Hilbert bundle $\Cal E \to M$
of finite type and extend the analysis of Helffer-Sj\"ostrand [HS1] to this
more general setting.

Assume that $(M,g)$ is a closed Riemannian manifold and
let $h: M \to \Bbb R$ be a Morse function, so that $\tau = (h,g)$ is
a generalized triangulation.
Let $\Cal W$ be a $(\Cal A,\Gamma^{op})-$Hilbert module of finite type
with $\Gamma=\pi_1(M).$
To simplify the exposition we assume that $\Cal W$ is a free $\Cal A$- Hilbert
module.
Denote by $\Cal E\to M$ the bundle of $\Cal A$-Hilbert
modules associated to $\Cal W$.
Let $x_{q;j} \in {\roman {Cr}}_q(h)$ be a critical point of index $q$ and
$U_{qj}$
an
open neighborhood of $x_{q;j}.$
\proclaim {Definition 5.1}
$U_{qj}$ is said to be an H-neighborhood of $x_{q;j}$ if there is a ball
$B_{2\alpha} := \{ x \in \Bbb R^d ; |x| < 2\alpha \}$ and diffeomorphisms
$\phi : B_{2\alpha} \to U_{qj}$ and $\Phi : B_{2\alpha} \times \Cal W \to
\Cal E |_{U_{qj}}$ with the following properties:
\newline (i) $\phi (0)= x_{q;j};$
\newline (ii) when expressed in the coordinates of $\phi, h$ is of the form
$$h(x) = q - (x_1^2 + ...+x_q^2)/2 + (x_{q+1}^2 + ... +x_d^2)/2;$$
\newline (iii) the pull back $\phi^*(g)$ of the Riemannian metric $g$ is
the Euclidean metric;
\newline (iv) $\Phi$ is a trivialization of $\Cal E |_{U_{qj}}.$

For later use we define $U'_{qj}:= \phi (B_{\alpha}).$

 A collection $(U_x)_{x \in Cr(h)}$ of H-neighborhoods
is called a system of H-neighborhoods if, in addition, they are pairwise
disjoint.
\endproclaim

 As in section 4,
denote by $\Lambda^q(M;\Cal E):=C^{\infty}(\Cal E\otimes\Lambda^q(T^*(M)))$
 the $\Cal A-$module of smooth
$q$-forms  with values in $\Cal E$ and by $L_2(\Lambda^q(M;\Cal E))$
its $L_2-$completion, which is an $\Cal A-$Hilbert module.  We write
$\Lambda^q(M;\Cal W)$ for
$\Lambda^q(M;\Cal E)$ when $\Cal E=M\times\Cal W$ is the trivial bundle and
$\Lambda^q(M;\Bbb R)$
for the space
of smooth $q$-forms on $M.$
Consider the Witten Laplacian
$\Delta_q(t):\Lambda^q(M;\Cal E)
\to \Lambda^q(M;\Cal E)$ associated to $h$
and observe that
$$\Delta_q(t)=\Delta_q + t^2||\nabla h||^2 + tL_q \tag {5.1}$$
where $L_q$ is a zero'th order differential $\Cal A$-operator on
 $\Lambda^q(M;\Cal E)$, hence given by a bundle
endomorphism, and where $|| \nabla h||^2$ is a scalar valued function on $M$
given by $|| \nabla h ||^2(x) = \sum_{1 \leq i,j \leq d} g^{ij}(x)
\frac {\partial h}{\partial x_i} \frac {\partial h}{\partial x_j}$ with
 $(g^{ij}(x))$ denoting the inverse
of the metric tensor $g$ when expressed in local coordinates. $\Delta_q$ is
a nonnegative, selfadjoint, elliptic differential $\Cal A$-operator.
Let $\Lambda^q(M;\Cal E)_{{\roman {sm}}}$  be the image (which depends on $t$)
 of the spectral projector $Q_q(1,t)$ of $\Delta_q(t),$
corresponding to the interval $(- \infty,1].$ This space consists of smooth
q-forms and is an $\Cal A$-Hilbert module.

The purpose of this section is to study
the complex $(\Lambda^*(M;\Cal E)_{{\roman {sm}}}, d_*(t))$
for $t$ sufficiently large and to precisely formulate
and prove that this family of complexes converges to the  cochain complex
 $\Cal C^*(M,\tau,O_h)$, introduced in section 4, when
$t\to \infty. $ In the case $\Cal A = \Bbb R$ and $\Cal W = \Bbb R$ this was
done  by Helffer
 and Sj\"ostrand [HS1].
Their arguments are still valid
 in the general case. Bismut and Zhang [BZ] verified this
in the case  $\Cal A=\Bbb R.$ Here we outline the proof for an arbitrary finite
von
Neumann algebra $\Cal A$, refering to [BZ] for those
 details whose verifications are
the same as in the case  $\Cal A = \Bbb R$.

Consider $h_k :\Bbb R^d \to \Bbb R$ defined by
$h_k(x)=k+ \frac {1}{2}(-\sum_1^k |x_i|^2+\sum_{k+1}^d|x_i|^2)$ and denote by
$\tilde \Delta_q:\Lambda^q(\Bbb R^d;\Bbb R)\to \Lambda^q(\Bbb R^d;\Bbb R)$
the flat Laplacian on $q$-forms on $\Bbb R^d$ and by
$\tilde \Delta_{q;k}(t):\Lambda^q(\Bbb R^d;\Bbb R)\to
\Lambda^q(\Bbb R^d;\Bbb R)$
the Witten Laplacian associated to $h_k$. A straightforward calculation shows
that
$$\tilde \Delta_{q;k}(t)= \tilde \Delta_q+t^2|x|^2-t(d-2k)
 +2t(N_{q;k}^+ -N_{q;k}^-),  \tag {5.2} $$
where $N_{q;k}^+$ and $N_{q;k}^-$ are the
number operators introduced in [HS1] (cf also[BZ]), defined by
$$N^+_{q;k}(dx_{i_1}\wedge ... \wedge dx_{i_q}) = \# \{k+1 \leq i_j \leq d\}
dx_{i_1}\wedge ... \wedge dx_{i_q}$$
and $N^-_{q;k}:= q{\roman {Id}} - N^+_{q;k}.$
 Denote by $\tilde\omega_q(t)\in \Lambda^q (\Bbb R^d;\Bbb R)$
the $q$-form defined by
$$\tilde\omega_q(t):=(t/\pi)^{d/4}e^{-t|x|^2/2} dx_1\wedge
...\wedge dx_q. \tag {5.3} $$

For $\eta > 0,$ let $\nu_{\eta}:\Bbb R\to [0,1]$ be a smooth map equal
 to $1$ on the interval $(-\infty,\eta/2)$ and equal to
$0$ on the interval $(\eta,\infty).$ For $\epsilon > 0,$ which we will choose
 later at our convenience,
define  $\tilde\psi_q(t) \in \Lambda^q(\Bbb R^d;\Bbb R)$ by
$$\tilde\psi_q(t):=\beta(t)^{-1}\nu_\epsilon(|x|)\tilde\omega_q(t)\tag {5.3'}$$
where $\beta(t)=||\nu_\epsilon(|x|)\tilde\omega _q(t)||_2$.
 With respect to the scalar  product in
$\Lambda^q(\Bbb R^d;\Bbb R)$ induced by
the flat metric of $\Bbb R^d,$
$\langle\tilde\omega_k(t),\tilde\omega_k(t)\rangle = 1$ and
$\langle\tilde\psi_k(t),\tilde\psi_k(t)\rangle = 1.$
Consider $\Delta_q=\tilde \Delta_q\otimes {\roman {Id}}$ and
 $\Delta_{q;k}(t)=\tilde \Delta_{q;k}(t)\otimes {\roman {Id}}$,
defined on $\Lambda^q(\Bbb R^d;\Cal W).$ Both are
nonnegative, essentially selfadjoint, elliptic
$\Cal A$-operators and have the following properties:
\newline (HO1) $\text{spec}\Delta_{q;k}(t)$ is discrete and contained in
$2t\Bbb Z_{ \geq 0};$ each eigenvalue has infinite multiplicity if
$\text{dim}_{\Bbb R}\Cal A =\infty. $
\newline (HO2) $\text{Null}(\Delta_{q;k}(t))=0$ if $k\ne q$;
$\text{Null}(\Delta_{q;q}(t))$ is an $\Cal A$-Hilbert module
isometric to $\Cal W.$
\newline (HO3) Assume that $\{ v_1,...,v_l\}$ is an orthonormal basis
 of $\Cal W,$ i.e. a collection of orthonormal
vectors which generate $\Cal W,$ as an $\Cal A$-Hilbert module and such that
for any
$a,b\in \Cal A,$
$$\langle av_i, bv_j\rangle = \langle a,b\rangle \delta_{ij}.  \tag {5.4}$$
Then $\omega_{q,i}:=\tilde\omega_q(t)\otimes v_i, 1\leq i \leq l,$
 is an orthonormal basis for
$\text{Null}(\Delta_{q;q}(t)).$
Similarly  $\psi_{q,i}:=\tilde\psi_q(t)\otimes v_i,$ $1 \leq i \leq l,$
 provide an orthonormal basis for the $\Cal A$-Hilbert submodule
generated by them.
A straightforward calculation, using (5.2), (5.3) and
 $(5.3')$ and (HO1), shows that there exist constants
$C(\epsilon),C_0(\epsilon)> 0$, so that, for $1 \leq i \leq l,$ and for $t$
 sufficiently large,
$$\langle\Delta_{q;q}(t)\psi_{q,i},
\Delta_{q,q}(t)\psi_{q,i}\rangle=O(e^{-C(\epsilon)t}),\tag {5.5}$$
$$\langle\Delta_{q;k}(t)\psi_{q,i},\psi_{q,i}\rangle
\geq  C_0(\epsilon)t ~~~ ( k \ne q)   \tag {5.6}$$
and for any $\omega\in \Lambda^q(\Bbb R^d;\Cal W)$ with
$\langle\omega,\psi\rangle=0$ for $\psi$
in the Hilbert module generated by $\psi_{q;i}, 1 \leq i \leq l,$
$$\langle\Delta_{q;q}(t)\omega,\omega\rangle \geq  C_0(\epsilon)t ||\omega||^2.
\tag {5.7}$$
 For any two points $y,z \in M$ denote by $d(y,z)$ the distance induced
 by the metric g and by  $d_A(y,z)$ the distance induced by the Agmon metric
$g_A=|\nabla h|^2g.$ Let $x \in {\roman {Cr}}_k(h)$ and $U_x$ be an
H-neighborhood
as defined above.
For $y\in U_x$ one has $d_A(x,y)= |y|^2/2$ and $d(x,y)=|y|.$
Choose $\epsilon > 0$ so that
the balls $B(x;4\epsilon )=\{y\in M; d(x,y)\le 4\epsilon\}$, centered at
 critical points $x,$
are pairwise  disjoint,
and $B(x;3\epsilon)\subset U_x.$
Choose once and for all a base point $x_0 \in M,$ an orthonormal basis
$e_1,...,e_l$ of $\Cal E_{x_0},$ and choose for each
 critical point $x=x_{q;j} \in {\roman {Cr}}_q(h)$
a homotopy class $[\gamma_x]$ of paths, joining $x_0$ and $x$ (choose
$\gamma_{x_0} = \{x_0\}).$ Denote by $e_{q;j,1},...,e_{q;j,l}$ the orthonormal
basis of $\Cal E_x$ obtained from $e_1,...,e_l$ by parallel transport
along $\gamma_x,$ provided by the canonical flat connection on $\Cal E.$
Using the parallel transport,
one can identify $\Cal E|_{U_x}$ with
$U_x\times\Cal W$
and, using a system of H-neighborhoods $U_x,$ one can identify the forms
$\omega \in \Lambda^q(M;\Cal E)$ having support
in $U_x$ with  forms in $\Lambda^q(\Bbb R^d;\Cal W).$
By extending $\psi_{x,i}(t),$ defined by this identification on $U_x,$ by
zero to all of $M$ one obtains a form in $\Lambda^q(M;\Cal E),$ which we again
denote by $\psi_{x,i}(t).$
 The forms $\psi_{x,i}(t) (1 \leq i \leq l, x \in {\roman {Cr}}_q(h))$
satisfy (5.4),
and therefore provide
an orthonormal basis for the
$\Cal A$-Hilbert module which they generate.

\proclaim {Proposition 5.2} {\sl For any $q$ there exist  positive constants
$C',C'',$ and $t_0$ so that
 \newline ${\roman {spec}}(\Delta_q(t)) \cap (e^{-tC'},C''t)=\emptyset$
for $t > t_0.$
}
\endproclaim

{\sl Proof} $\ $ In a first step we prove that for $t > t_0$,
 with $t_0$ sufficiently large, there exists
a pair of orthogonal closed subspaces
$W_1=W_1(t), W_2= W_2(t)$ of $L_2(\Lambda^q(M;\Cal E))$
with $W_1 \subset \Lambda^q(M;\Cal E)$ so that the following
properties hold:
(1) $ W_1\cap W_2=\{0\};$ (2) $ W_1 + W_2=L_2(\Lambda^q(M;\Cal E));$
 (3) $ \|\Delta_q(t)\omega\| \leq  e^{-t2C'}\| \omega \|$
for $\omega \in W_1;$ and
(4) $\langle\Delta_q(t)\omega,\omega\rangle \geq 2 C'' t \langle\omega,
\omega\rangle$ for
 $\omega \in W_2 \cap \Lambda^q(M;\Cal E).$

In a second step we show that, using step 1, Proposition 5.2 follows.
Let us prove step 2 first. We argue by contradiction. Assume that there
exists a sequence $t_j \rightarrow \infty$ and real numbers
$\mu_j \in {\roman {spec}}\Delta_q(t_j) \cap (e^{-t_{j}C'},C''t_j).$
For each $j \geq 1,$ one can find an approximate eigenfunction $u_j$ in
the domain of $\Delta_q(t_j), \|u_j\|=1,$ satisfying
$$\| \Delta_q(t_j)u_j - \mu_ju_j \| \leq e^{-4C't_j}.$$
Decomposing $u_j = v_j + w_j \in W_1(t_j) \oplus W_2(t_j),$ one verifies,
using the fact that $\Delta_q(t)$ is selfadjoint,
$$|\langle \Delta_q(t_j)u_j,v_j\rangle | \leq \|\Delta_q(t_j)v_j\|\|v_j\|
+\|w_j\|\|\Delta(t_j)v_j\|$$
as well as
$$\mu _j \|v_j\|^2=\langle \mu_ju_j,v_j\rangle \leq |\langle
\Delta_q(t_j)u_j,v_j\rangle |
+ |\langle\Delta_q(t_j)u_j - \mu_ju_j,v_j\rangle |.$$
Together with property (3) these two inequalities imply
$$\mu_j \|v_j\|\|w_j\| \leq e^{-2C't_{j}}(\|v_j\|\|w_j\| + \|w_j\|^2)
+e^{-4C't_j}\|w_j\|.$$
It remains to prove step 1. Define $W_1:= W_1(t)$
to be the $\Cal A$-Hilbert module
generated by $\psi_{x,i}(t) ( 1 \leq i \leq l, x\in {\roman {Cr}}_q(h))$ and
$W_2:= W_2(t)$ its orthogonal complement in $L^2(\Lambda^q(M;\Cal E))$.
Clearly properties (1) and (2) are satisfied.
Further note that an element $\omega\in W_1$ has a representation
$\omega=\sum_{1 \leq i \leq l, x\in {\roman {Cr}}_q(h)}
a_{x,i}\psi_{x,i}(t)$ with $a_{x,i}\in \Cal A$ and that $\Delta_q (t)$, when
 restricted to $U_x$ with $x \in {\roman {Cr}}_k(h)$ and expressed in local
 coordinates introduced above,
 coincides with $\Delta_{q;k} (t).$ Therefore, in view of (5.4),(5.5) and the
support properties of $\psi_{q,i},$ we have, with $C=C(\epsilon)$ as in (5.5),
$$\langle\Delta_q(t) \omega, \omega\rangle =
 \sum_{1 \leq i \leq l, x \in {\roman {Cr}}_q(h)}
\langle a_{x,i}\Delta_q(t) \psi_{q,i}, a_{x,i} \Delta_q(t)\psi_{q,i}\rangle$$
$$\leq \sum_{i,x} ||a_{x,i}||^2 e^{-t2C} \leq ||\omega||^2 e^{-t2C}.$$
It remains to check the estimate (4).
 Denote by $\chi_x:M\to \Bbb R$ the smooth cut-off function with
support in $U_x$ defined by $\nu_{2\epsilon}$ and introduce $\chi = \sum_{x \in
Cr(h)} \chi_x.$
 For $\omega \in W_2 \cap
 \Lambda^q(M;\Cal E),$
define $\omega_1=\chi\omega$ and
$\omega_2=(1-\chi)\omega$ and observe that the support of $\omega_2$ is
 disjoint
 from the support
of any element in $W_1$; therefore $\omega_2\in W_2 \cap \Lambda^q(M;\Cal E)$
 and hence
$\omega_1\in W_2 \cap \Lambda^q(M;\Cal E).$
Since $\Delta_q(t)$ is essentially selfadjoint  one obtains
$$\langle\Delta_q(t)\omega,\omega\rangle =  \langle\Delta_q(t)\omega_1,
\omega_1\rangle +
2\langle\Delta_q(t)\omega_1,\omega_2\rangle +
\langle\Delta_q(t)\omega_2,\omega_2\rangle. \tag{5.10}$$
We show that there exist positive constants $t_0,C_1,C_2,C_3,C_4 $
 depending only on the
geometry of $(M,\Cal E\to M)$
and the chosen $\epsilon ,$ so that for any
$\omega\in W_2 \cap \Lambda^q(M;\Cal E)$ and $t > t_0$ the
following estimates
hold:
$$\langle\Delta_q(t)\omega_2,\omega_2\rangle \geq  \langle\Delta_q\omega_2,
\omega_2\rangle + C_1t^2||\omega_2||^2-
C_2t||\omega_2||^2;  \tag{5.11}$$
$$\langle\Delta_q(t)\omega_1,\omega_1\rangle \geq  C_3t||\omega_1||^2;
\tag{5.12}$$
$$\langle\Delta_q(t)\omega_1,\omega_1\rangle \geq \langle\Delta_q\omega_1,
\omega_1\rangle -C_2t
||\omega_1||^2. \tag{5.13}$$
For $\alpha > 0,$ $ \langle\Delta_q(t)\omega_1,\omega_2\rangle$ is bounded
 from below by
$$ -2C_4(1+\alpha^{-2})(||\omega_1||^2 + ||\omega_2||^2)
-C_4\alpha^2\langle\Delta_q\omega_2,\omega_2\rangle  -C_4\alpha^2
\langle\Delta_q\omega_1,
\omega_1\rangle.  \tag{5.14}$$

Note that (5.12) and (5.13) imply that for any $0 \leq \delta \leq 1$
$$\langle\Delta_q(t)\omega_1,\omega_1\rangle \geq  (1-\delta)
\langle\Delta_q\omega_1,\omega_1\rangle
 +t(\delta C_3 - (1-\delta)C_2)||\omega_1||^2.  \tag{5.15}$$

To prove  (5.11) choose $C_1:=\inf_{z\in M\backslash\cup_{x
\in {\roman {Cr}}(h)} U_x}||\nabla h(z)||^2$ and
$C_2=\sup_{x\in M}||L(x)||.$ The estimate (5.11) then follows from (5.1).

To prove (5.12) it suffices to notice that the support of  $\omega_1$
 is contained in
$\cup_{x\in {\roman {Cr}}(h)} U_x$ and $\omega_1$ is orthogonal
to $\psi_{x,i} (x \in {\roman {Cr}}_q(h), 1 \leq i \leq l).$
Thus (5.12) follows from (5.7) with $C_3 := C_0(\epsilon).$

Formula (5.13) is a direct consequence of (5.1).

To find the lower bound (5.14) note that
 $|\langle L\omega_1,\omega_2\rangle | \leq C_2|\langle \omega_1,
\omega_2 \rangle |$  and, using
that supp($\omega_2$)
 does not intersect any of the $U_x's,$
$\langle||\nabla h||^2 \omega_1, \omega_2\rangle \geq C_1 \langle\chi(1-\chi)
\omega, \omega\rangle
 \geq 0.$
Combining with (5.1) one concludes
$$\langle\Delta_q(t) \omega_1, \omega_2\rangle = \langle\Delta_q \omega_1,
\omega_2\rangle +
t^2 \langle ||\nabla h||^2 \omega_1, \omega_2\rangle + t\langle L
 \omega_1, \omega_2\rangle$$
$$\geq \langle\Delta_q\omega_1,\omega_2\rangle +(C_1t^2-C_2t)
\langle\omega_1,\omega_2\rangle.$$
For $t > C_2/{C_1}$ one thus obtains
$$\langle\Delta_q(t)\omega_1,\omega_2\rangle \geq \langle\Delta_q\omega_1,
\omega_2\rangle. \tag{5.16}$$
Therefore, the lower bound (5.14) follows from Lemma 5.3 below.

To complete the proof of property (4) combine (5.10) with the estimates
(5.11), (5.14), and (5.15) to obtain for $0 < \delta <1,$ $ \alpha > 0$ and
 $t > C_2/C_1,$
$$ \langle\Delta_q(t)\omega, \omega\rangle
 \geq (1-2C_4 \alpha^2) \langle\Delta_q \omega_2, \omega_2\rangle
+ (1- \delta - 2C_4\alpha^2)\langle\Delta_q \omega_1, \omega_1\rangle$$
$$+ (C_1t^2 - C_2t - 4C_4(1+\alpha^{-2}))||\omega_2||^2
+ (t(\delta C_3 - (1-\delta)C_2) - 4C_4(1+\alpha^{-2}))||\omega_1||^2.$$
First choose $0 < \delta < 1$ sufficiently close to 1 so that
$\delta C_3 - (1 - \delta)C_2 > 0.$ Then choose $\alpha > 0$ sufficiently
small so that $1 - \delta - 2C_4\alpha^2 > 0.$
Together with $2(||\omega_1||^2 + ||\omega_2||^2) \geq ||\omega||^2$
this establishes property (4)
for $t > t_0$ if $t_0 > C_2/C_1$ is chosen sufficiently large.$\diamondsuit$

\proclaim {Lemma 5.3} {\sl Let the q-forms $\omega,$ $\omega_1$ and $\omega_2$
be defind as above. Then there exists a constant $C_4 > 0$
 so that, for any $\alpha > 0$,
$$\langle\Delta_q \omega_1, \omega_2\rangle
\geq -C_4(1+\alpha^{-2})||\omega||^2
- C_4 \alpha^2\langle\Delta_q \omega_2, \omega_2\rangle - C_4\alpha^2
\langle\Delta_q\omega_1,
 \omega_1\rangle.$$
}
\endproclaim
{\sl Proof} $\ $ Write $\Delta_q = d_{q-1}d_{q-1}^* + d_q^*d_q$ where
$d_{q-1}^* = -(-1)^{dq}R_{d-q+1}d_{d-q}R_q,$ and $\ast = R_q$ denotes
the Hodge $\ast$ operator. Using that $\omega_1 = \chi \omega$ and $\omega_2 =
(1- \chi)\omega$ one obtains
$$\langle \Delta_q \omega_1, \omega_2\rangle = \langle d\omega_1,
d\omega_2\rangle
+ \langle d\ast \omega_1, d\ast \omega_2\rangle  \geq A + B,$$
where

$$A := \langle d\chi \wedge \omega,u(1-\chi)d\omega\rangle +
\langle d \chi\wedge
\ast \omega,u(1-\chi)d\ast \omega\rangle, $$
$$B := -\langle\chi d\omega,ud\chi\wedge\omega\rangle
- \langle\chi d\ast \omega,
ud\chi\wedge *\omega\rangle,$$
where $u$ is the characteristic function of $M\backslash {supp} \chi.$

In order to estimate the expressions $A$ and $B$ we introduce  the constant
\newline $C_5 := \sup_{1 \leq k \leq d} |||K_k|||,$ where
 $K_k:L_2(\Lambda^k(M;\Cal E)) \to L_2(\Lambda^{k+1}(M;\Cal E))$ is the
 exterior multiplication by $d\chi$. Note that $|\|K_k\| | = |\|K_k^*\| |$
 where
$K_k^*$ denotes the adjoint of $K_k.$ A straightforward calculation yields

$$|A|\leq C_5||\omega||( ||(1-\chi)d\omega||+   ||(1-\chi)d*\omega||) $$
$$\leq C_5||\omega||( ||d\omega_2||+ ||d\chi\wedge\omega|| +||d\ast \omega_2||+
||d\chi\wedge\ast \omega||)$$
$$\leq C_5||\omega||( ||d \omega_2||+ ||d \ast \omega_2||+ 2C_5||\omega||) $$
$$\leq \sqrt2 C_5||\omega||\langle\Delta_q\omega_2,\omega_2\rangle^{1/2}+
2C_5^2||\omega||^2.$$
Thus for any $\alpha > 0$
$$|A| \leq (2C_5^2 + C_5\alpha^{-2})||\omega||^2
 + \alpha^2\langle\Delta_q\omega_2, \omega_2\rangle.$$
A similar computation leads to
$$|B| \leq (2C_5^2 + C_5\alpha^{-2})||\omega||^2
+ \alpha^2 \langle\Delta_q\omega_1, \omega_1\rangle.$$
Choosing $C_4$ appropriately leads to the claimed statement. $\diamondsuit$

Proposition 5.2 yields, for $t$ sufficiently large,a decomposition of
$(\Lambda^q(M;\Cal E),d_q(t))$
$$(\Lambda^q(M;\Cal E),d_q(t)) = (\Lambda^q(M;\Cal E)_{{\roman {sm}}},d_q(t))
\oplus
(\Lambda^q(M;\Cal E)_{{\roman {la}}}, d_q(t))$$
where $\Lambda^q(M;\Cal E)_{{\roman {sm}}}$ is the image (depending on t,) of
$Q(1,t),$ the spectral projection of $\Delta_q(t)$ corresponding
to the interval $(-\infty , 1],$ and $\Lambda^q(M;\Cal E)_{{\roman {la}}}$
denotes
the orthogonal complement of $\Lambda^q(M;\Cal E)_{{\roman {sm}}}.$
Accordingly,
one can decompose $\Delta_q(t) = \Delta_{q,{\roman {sm}}}(t) +
\Delta_{q,{\roman {la}}}(t)$
where $\Delta_{q,{\roman {sm}}}(t)$ denotes the restriction of $\Delta_q(t)$ to
$\Lambda^q(M;\Cal E)_{{\roman {sm}}}$ and, similarly,
$\Delta_{q,{\roman {la}}}(t)$ denotes
the restriction to $\Lambda^q(M;\Cal E)_{{\roman {la}}}.$

Now assume that $(M,\Cal W)$ is of determinant class. Then, for
t sufficiently large,
\newline ${\roman {logdet_N}} \Delta_q(t),$
${\roman {logdet_N}} \Delta_q(t)_{{\roman {sm}}},$ and ${\roman {logdet_N}}
\Delta_q(t)$
are all real numbers and, accordingly,  we write
$$\log T_{{\roman {an}}}(t) = \log T_{{\roman {sm}}}(t) +
\log T_{{\roman {la}}}(t).$$

Let $x = x_{q;j} \in {\roman {Cr}}_q(h).$
Choose $\epsilon > 0$ and let $\{ e_{q;j,i} \}$ be the orthonormal basis
of $\Cal E_x$ as defined above.
Define the $\Cal A$- linear maps
$J_x(t): \Cal E_x \to L_2(\Lambda^q(M;\Cal E))$ by
$$J_x(t)\left( \sum_i a_i e_{q;j,i} \right) := \sum_i a_i\psi_{x,i}.
\tag{5.20}$$
We point out that the $\psi_{x,i}$'s, and thus $J_x(t),$ depend on the
choice of $\epsilon$
and notice that $J_x(t)$ is an $\Cal A$-linear isometry. Let
$J_q(t):\sum_{x\in {\roman {Cr}}_q(h)}\Cal E_x
\to L_2(\Lambda^q(M;\Cal E))$ be the sum  $J_q(t) :=
\sum_{x \in {\roman {Cr}}_q(h)} J_x(t)$.
As the images of $J_x(t)$  have disjoint
support, the map $J_q(t)$ is also an isometry.
Recall that we have denoted by $Q_q(1,t)$ the spectral projector of
 $\Delta_q(t)$ corresponding to the interval $(-\infty , 1].$
Introduce the map
$$H_q(t):= (Q_q(1,t)J_q(t))^*(Q_q(1,t)J_q(t)), \tag{5.21}$$
where $\ast$ denotes the adjoint of an operator.  $H_q(t)$ is a
selfadjoint, nonnegative, bounded, $\Cal A-$linear operator
on $\sum_{x \in {\roman {Cr}}_q(h)} \Cal E_x.$

\proclaim {Proposition 5.4}
For $\epsilon > 0$ sufficiently small, there exists a constant $c > 0$
so that
$$(Q_q(1,t) J_q(t)v-J_q(t)v)(x)=O(e^{-ct}||v||) \tag{5.22}$$
uniformly in $x\in M$ and $v \in \sum_{x \in {\roman {Cr}}_q(h)} \Cal E_x$ and
 $$ H_q(t)= {\roman {Id}} + O(e^{-ct}). \tag {5.23}$$
\endproclaim

Observe that the composition $Q_q(1,t)J_q(t)H_q(t)^{-1/2}$ is an
$\Cal A$-linear
isometry from
the $\Cal A-$Hilbert module $\sum_{x\in {\roman {Cr}}_q(h)}\Cal E_x $ to
$\Lambda^q(M;\Cal E)_{{\roman {sm}}}.$

{\sl Proof} We proceed as in \cite{BZ, p.128}. In view of Proposition 5.2,
for $t > t_0, \ Q_q(1,t)$  is given by the Riesz projector
$$ Q_q(1,t)=\frac {1}{2\pi i}\int_{S^1}(\lambda-\Delta_q(t))^{-1}d\lambda
 \tag{5.24}$$
where $S^1$ is the unit circle in $\Bbb C,$ centered at the origin.
The operator $Q_q(1,t)J_q(t) - J_q(t)$ can therefore be represented by
a Cauchy integral whose integrand is given by
$$
(\lambda-\Delta_q(t))^{-1}J_q(t) - \lambda ^{-1}J_q(t) =
\lambda^{-1}(\lambda-\Delta_q(t))^{-1}\Delta_q(t) J_q(t).
\tag{5.25}
$$
By Proposition 5.2 there exists, for any Sobolev norm $||.
||_r,$ a constant $c_r > 0$ so that
$$||\Delta_q(t) J_q(t)(v)||_r = O(e^{-c_{r}t}||v||),  \tag{5.26}$$
uniformly in $v \in \sum_{x \in {\roman {Cr}}_q(h)} \Cal E_x.$

By proceeding as in \cite{BZ, p 128-129}
one can show that there exists $C_r$ so that  for $t$ sufficiently large and
any $\lambda\in S^1$
$$|||(\lambda-\Delta_q(t))^{-1})|||_{r \to r} \leq C_r t^r. \tag{5.27}$$

Combining (5.26) and (5.27), one obtains, for $c' < c_r,$
 uniformly for $y \in M$ and $v \in \sum_{x \in {\roman {Cr}}_q(h)} \Cal E_x,$
$$||(\lambda-\Delta_q(t))^{-1}\Delta_q(t)J_q(t)v||_r
= O(e^{-c't})||v||.  \tag{5.28}$$
Choose $r > d/2$ and use the Sobolev embedding theorem to obtain (5.22)
 from (5.24), (5.25) and (5.28). (5.23) follows
 immediately from (5.22).$\diamondsuit$

Let us now consider the cochain complex $\Cal C(M,\tau,O_h,\Cal W),$
which has been
introduced in section 4.
Define  $E_{q;j,i} \in \Cal C^q$  for
$1\leq j \leq m_q$ and $1\leq i \leq l$ by
$$
E_{q;j,i}(x_{q;j'})=
\left\{\aligned e_{q;j,i}~~~ \text { if } j'=j \\0 ~~~~\text
{ if } j' \ne j.\endaligned \right. \tag{5.29}
$$
We see that $E_{q;j,i}$ is bounded as follows:
Assume that $\Cal W'$ is a free $\Cal A$-Hilbert module of
finite type with an orthonormal basis $v_1,...,v_l$
and $f:\Cal W\to \Cal W'$ is a bounded,
 $\Cal A$-linear map. Then
$$ ||f||\leq  l^{1/2} \sup\{||f(v_i)||, 1 \leq i \leq l\} \leq l^{1/2} ||f||.
 $$
With respect to this basis the differential $\delta_q$ can be written as
$$
\delta_q(E_{q;j,i})=\sum_{1\leq j'\leq m_{q+1},1\leq i'\leq
l}\gamma_{q;ji,j'i'}E_{q;j',i'}
\tag{5.30}
$$

Denote by $\Cal G^q$ the $\Cal A$-Hilbert module (which depends on $t$)
 generated by
$J_q(t)(e_{q;j,i})$ with $ 1 \leq j \leq m_q, 1 \leq i \leq l;$
 this is $W_1(t)$ introduced in Proposition 5.2. Recall that, given two closed
subspaces $W_1$ and $W_2$,
of a Hilbert space, the semi-distance  between them is defined by
 $sdist(W_1,W_2):= |||Prj_{W_1}-Prj_{W_2}\cdot
Prj_{W_1}||| = |||Prj_{W_1}-Prj_{W_1}\cdot Prj_{W_2}|||,$ where $Prj_{W_i}$
denotes the orthogonal projector on the subspace $W_i.$
Following Helffer and Sj\"ostrand [HS1] (p 262) we write
 $A(t)=\tilde O(e^{-tc})$ for a quantity $A(t)$ depending on the choice of
the $\epsilon-$dependent collection $(U_x)_{x \in {\roman {Cr}}(h)}$ of
H-neighborhoods,
if for any $\delta > 0$ there exists
 $\epsilon_{\delta} > 0$ so that, for any collection
 $(U_x)_{x \in {\roman {Cr}}(h)}$
of H-neighborhoods with $\epsilon \le \epsilon_{\delta},$
$A(t)=O(e^{-t(c-\delta)})$.

\proclaim {Proposition 5.5}
$sdist(\Cal G^q,\Lambda^q(M;\Cal E)_{{\roman {sm}}})=\tilde O(e^{-tS_q})$ where
$S_q=\inf_{x,y\in {\roman {Cr}}_q(h)} d_A(x,y).$
 \endproclaim

Here $d_A(x,y)$ denotes the Agmon distance associated to $(M,g,h)$ as
defined in  \cite{HS1}.
Proposition 5.5 is a generalization of Proposition (1.7) of \cite{HS1}
or Theorem (8.15) of \cite{BZ} and can be proved in the same way as in
\cite{HS1} once one generalizes Proposition 2.5 in \cite{HS2} as follows
\proclaim {Proposition 5.6}
 Let $\Cal K$ be an $\Cal A$-Hilbert module and $\Cal K'$ an $\Cal A-$
Hilbert submodule of $\Cal K$ with  orthonormal basis $\psi_1,...,\psi_N.$
Suppose that for any $a,b\in \Cal A,$
$\langle a\psi_i,b\psi_j \rangle = \langle a,b \rangle \delta_{ij}.$
Let $f:\Cal K \to \Cal K$ be a
selfadjoint, nonnegative, $\Cal A$-linear operator
with
$(\alpha,\beta) \cap {\roman {spec}}( f) = \emptyset$ for some real numbers
$\alpha$ and $\beta$, where $0 < \alpha < \beta.$
Suppose that $f(\psi_i)=\mu_i\psi_i +r_i$
where  $r_i \in \Cal K$ with $||r_i|| < \epsilon$ and
$\mu_i$ are real numbers satisfying $0 \leq \mu_i < \alpha$. Let ${\Cal
K}_{{\roman {sm}}}$
denote the range of the spectral projection of $f$
associated to the interval $[0,\alpha]$. Then
${\Cal K}_{{\roman {sm}}}$ is an $\Cal A$-Hilbert module and
$$sdist(\Cal K',{\Cal K}_{{\roman {sm}}}) \leq
\frac{N^{1/2}\epsilon}{\beta -\alpha}.$$
\endproclaim

{\sl Proof} $\ $ Write $(f-\lambda) \psi_i = (\mu_i-\lambda)\psi_i +r_i$
and note that for
$\lambda \in \Bbb C\backslash(\text {spec} f \cup \{ \mu_1,...,\mu_N\})$
$$(f-\lambda)^{-1}\psi_i=(\mu_i-\lambda)^{-1}\psi_i-
(\mu_i-\lambda)^{-1}(f-\lambda)^{-1}r_i.$$
Denote by $\gamma_R$ the oriented boundary of $[-(\beta-\alpha)/2,
(\beta+\alpha)/2] \times
i[-R,R]$ and by
$P_{\Cal K'}$  the orthogonal projection
on $\Cal K'.$  Applying Cauchy's formula one obtains
$$(P_{{\Cal K}_{{\roman {sm}}}} \psi_i-\psi_i)
= -\frac{1}{2\pi i}
\int_{\gamma_R}(\mu_i-\lambda)^{-1}(f-\lambda)^{-1}r_i d\lambda.$$
Letting $R \to \infty,$ the above integral becomes
$$
\frac{1}{2\pi i}\int^{-(\beta-\alpha)/2+i\infty}_{-(\beta-\alpha)/2-i\infty}
(\mu_i-\lambda)^{-1}
(f-\lambda)^{-1}r_i d\lambda -
\frac{1}{2\pi i}\int^{(\alpha+\beta)/2+i\infty}_{{(\alpha+\beta)/2-i\infty}}
(\mu_i-\lambda)^{-1}(f-\lambda)^{-1}r_i d\lambda
$$
For $\lambda= -(\beta+\alpha)/2+it, \text{ or } \lambda=(\alpha+\beta)/2+it,
$ with
$-\infty < t < \infty$, one obtains
$$||(\mu_i-\lambda)^{-1}(f-\lambda)^{-1}r_i||\leq
\frac{\epsilon}{(\beta-\alpha)^2/4+t^2}.$$ Hence
$$||P_{{\Cal K}_{{\roman {sm}}}}\psi_i-\psi_i|| \leq
\frac {\epsilon}{2\pi} \int^{\infty}_{-\infty}
\frac{dt}{(\beta-\alpha)^2/4+t^2}=
\frac{\epsilon}{\beta-\alpha}.~~~\diamondsuit$$

Let
$$\alpha_q(t):=\sup\{0,{\roman {spec}}( Q_q(1,t)\Delta_q(t))\};~~ $$
$$\beta_q(t):=\inf\{1+{\roman {spec}}\left(({\roman
{Id}}-Q_q(1,t))\Delta_q(t)\right)\} -1.$$

\proclaim {Theorem 5.7}(\cite{HS1},\cite{BZ})
\newline (1) For $t\to\infty,$~~~  $\alpha_q(t)\to 0$ and $\beta_q(t)\to
\infty.$
\newline (2) There exists a constant $t_1$ so that for $t >t_1$ the elements
$$
\varphi_{q;j,i}(t)=Q_q(1,t)J_q(t)H_q(t)^{-1/2}(e_{q;j,i}) \tag{5.31}
$$
form an orthonormal basis for $ \Lambda^q(M;\Cal E)_{{\roman {sm}}}.$ Hence
 $\Lambda^q(M;\Cal E)_{{\roman {sm}}}$
is a  free $\Cal A-$Hilbert module
of rank $l\times \#{\roman {Cr}}_q(h).$
\newline (3) There exist $\eta > 0$ and $C > 0$ such that for $t$
sufficiently large, $1 \leq r \leq l,$
$$\sup_{ x \in M \setminus U_{qj}} \parallel \varphi _{q;j,r} \parallel
\leq C e ^{-\eta t}.$$
\newline (4) Let $W^-_{q;j}$ denote the unstable manifold with respect to
the flow corresponding to $grad _g h$ at the critical point $x_{q;j}$ of $h$.
Choose a system of H-neighborhoods $(U_{x_{q;j}})$ so that
$U_{x_{q;j}}\cap W^-_{q;j'} = \emptyset$ for $j'\ne j.$
When expressed in local coordinates on $U_{x_{q;j}} \cap W^-_{q;j}$ the
q-forms $\varphi_{q;j,i}(t)$ satisfy the following estimate:
$$
\varphi_{q;j,i}(t)=(t/4)^{(d/4)}e^{-t|x|^2/2}(dx_1\wedge ...\wedge dx_q\otimes
e_i + O(t^{-1})).
$$
\newline (5) Representing $d_q(t)$ with respect to this basis
$$
d_q(t)\varphi_{q;j,i}(t)=\sum_{1\leq j'\leq m_{q+1},1\leq
i'\leq l}\eta_{q;ji,j'i'}(t)
\varphi_{q;j',i'}(t)
$$
the coefficients $\eta_{q;ji,j'i'}$ satisfy
$$\eta_{q;ji,j'i'}(t)=e^{-t}(t/\pi)^{1/2}(\gamma_{q;ji,j'i'}
+O(t^{-1/2})).$$
\endproclaim

{\sl Proof} Statement (1) follows from Proposition 5.2.
Concerning statement (2), note that the $\Cal A-$Hilbert module $W_1$ of
 finite type
as defined in the proof of Proposition 5.2 is free, of
 rank $l\times \#{\roman {Cr}}_q(h)$ and  contained in
$\Lambda^q(M; \Cal E)_{{\roman {sm}}}.$ Therefore, should  $W_1$ be not equal
to
 $\Lambda^q(M; \Cal E)_{{\roman {sm}}}$ one could conclude that
$\Lambda^q(M;\Cal E)_{{\roman {sm}}}\cap  W_2 \ne \{ 0 \},$ where $W_2$ is
the orthogonal complement of $W_1$ as defined in the proof of Proposition 5.2.
 In view of Proposition 5.2 this is, however, not possible.
To verify the estimates (3), (4) and (5) one follows the arguments
 in \cite{HS1} (Proposition 1.7,Theorem 2.5 and Proposition 3.3) or \cite{BZ}
 (Theorem 8.15, Theorem 8.27,Theorem 8.30).
$\diamondsuit$

We need an application of the above results (cf \cite{BZ}):
\proclaim {Corollary 5.8}
$$ Int^{(q)} ( e^{ht} \phi_{q;j,r}(t))=
 \left( \frac {t}{\pi} \right) ^{(d-2q)/4} e^{qt}
(E_{q;j,r} + O(t^{-1})).$$
\endproclaim

{\sl Proof} We must show that for any cell $W^-_{q;j'}$
$$ \int_{W^-_{q;j'}} \phi_{q;j,r} (t) e^{ht} =
\left( \frac {t}{\pi} \right)^{(d-2q)/4} e^{tq} (\delta_{jj'} e_{q;j,r}
+ O(t^{-1})).$$
First, note that, due to Theorem 5.7 and to the choice of $U_{qj'},$
it suffices to consider the case where $j=j'.$ Moreover, it
suffices to estimate
$$ \int_{W^-_{q;j} \cap U_{qj}} \phi_{q;j,r}(t) e^{ht}.$$
Note that on $W^-_{q;j} \cap U_{qj},$ the function $e^{ht}$ is of the form
$$ e^{ht} = e^{qt} e^{-t(\sum_1^q x^2_k)/2}.$$
By Theorem 5.7, we conclude that
$$\int_{W^-_{q;j} \cap U_{qj}} \phi_{q;j,r} (t) e^{ht} =
\left( \frac {t}{\pi} \right) ^{d/4} e^{qt} \int_{W^-_{q;j} \cap U_{qj}}
e^{-t \sum_1^q x_k^2} (dx_1 \wedge ... \wedge x_q e_{q;j,r} + O(t^{-1}))$$
$$= e^{qt} \left( \frac {t}{\pi}
 \right) ^{d/4} \left( \frac {t}{\pi} \right)^{-q/2}
(e_{q;j,r} + O(t^{-1})).\ \ \ ~~~~\diamondsuit$$

Finally we state and prove Proposition 2, which is a generalized version
of Proposition $1$, mentioned in the introduction, and which is due
to Gromov-Shubin (cf also $[Ef]$).
\proclaim {Proposition 2} (\cite{Ef},\cite{GS})  Let ${\Cal W}$ be an ${\Cal
A}$- Hilbert module of finite type, not necessarily free.  Then the
following statements are true:
\roster
\item  Suppose $g$ is a Riemannian metric and $\tau = (h,g')$ is a generalized
triangulation of $M.$ Then the system $(M,g,\Cal W)$ is of
$a-determinant$ class iff $(M, \tau, \Cal W)$ is of $c-determinant$
class.

\item  If $M_1$ and $M_2$ are two homotopy equivalent connected manifolds
and $\tau _1$ and $\tau_2$ are generalized triangulations of $M_1,$
respectively
$M_2$, then $(M_1, \tau_1, \Cal W)$ is of $c-determinant$ class iff
$(M_2,\tau_2,\Cal W)$ is of $c-determinant$ class.
\endroster
\endproclaim

{\sl Proof} (1) Results of Novikov-Shubin (\cite{NS1,2}) imply that
$(M,g,\Cal W)$ is of $a-determinant$ class iff $(M,g',\Cal W)$ is.
Since, for fixed t, multiplication by $e^{th}$ provides a bounded
isomorphism of $L_2(\Lambda^*(M; \Cal E)) \to L_2(\Lambda^*(M;\Cal E))$
which intertwines $d_k(t)$ with $d_k$ it follows from results of
Gromov-Shubin (\cite{GS}) that $(M,g',\Cal W)$ is of $a-determiant$
class
$$\int ^1_{0+} \log \lambda dN_{\Delta_{k}(t)}(\lambda) > -\infty$$
for all values of $t$.
By Theorem 5.7 and Corollary 5.8 this is equivalent to saying that
 $(\Lambda^k(M;\Cal E)_{{\roman {sm}}} ,d_k(t))$ and
$(\Lambda ^k(M;\Cal E)_{{\roman {sm}}} ,\tilde d_k(t))$
are  of determinant class where
\newline $\tilde d_k(t):= e^{t} (\frac {t}{\pi})^{\frac {-1}{2}}d_k(t).$
 Notice that
$$f_k(t): \Lambda^k(M;\Cal E)_{{\roman {sm}}} \to \Cal C^k,$$
defined by
$$f_k(t) = \left( \left( \frac {\pi}{t} \right) ^{\frac {d-2k}{4}} e^{-tk}
\right)
Int^{(k)} e^{th}$$
establishes an isomorphism between $(\Lambda^k(M;\Cal E)_{{\roman {sm}}},
\tilde d_k(t))$
and $\Cal C^k$ and therefore, by Proposition 1.18, one concludes
that $(M,g',\Cal W)$ is of $a-determinant$ class iff $(M,\tau,\Cal W)$
is of $c-determinant$ class. Statement (2) is a direct consequence
of Proposition 1.18.
\vfil\eject

\subheading{ 6.1 Asymptotic expansion of Witten's deformation of the analytic
torsion}

Let $(M,g)$ be a Riemannian manifold with fundamental group $\Gamma =
\pi_1(M)$ and $h : M \longrightarrow {\Bbb R}$ a Morse function so that
$\tau=(h,g)$
is a generalised triangulation.
Let ${\Cal A}$ be a finite von Neumann algebra and ${\Cal W}$
an $(\Cal A,\Gamma^{op})$-Hilbert module of finite type.  The canonical bundle
$p:\Cal E\longrightarrow M$ associated to
$\Cal W$ is equipped with a canonical flat
connection which in turn induces (via parallel transport) a Hermitian
structure $\mu$ on ${\Cal E} \longrightarrow M.$ Throughout
this subsection we assume that $(M,{\Cal W})$ is of determinant class.

\proclaim{Definition}{A function $a: {\Bbb R} \longrightarrow {\Bbb R}$
is said to
have an asymptotic expansion for $t \longrightarrow \infty$ if there
exists a sequence $i_1 > i_2 > \dots > i_N = 0$ and constants
$(a_k)_{1\leq k \leq N}, \ (b_k)_{1\leq k \leq N}$ such that $$
\align
a(t) = &
\sum_1^N a_k t^{i_k} + \sum_1^N b_k t^{i_k}\log{t} + o(1).\tag{6.1}
\endalign$$}
\endproclaim

For convenience we denote by ${\roman {FT}}(a(t))$ the coefficient
$a_N$ in the asymptotic expansion of $a(t)$ corresponding to $t^0.$

Denote by $\beta_q$ the Betti numbers and by $\chi(M,\tau) = \sum_q
(-1)^q \beta_q$ the Euler-Poincare characteristic of the cochain
complex ${\Cal C}(M,\tau,O_h).$

Recall that in section 5 we introduced $T_{{\roman {an}}}(h,t), \
T_{{\roman {sm}}}(h,t) $ and $T_{{\roman {la}}}(h,t)$ and in
section 4 we introduced $T_{{\roman {Re}}}(\tau), \ T_{{\roman
{comb}}}(\tau) $ and $T_{{\roman {met}}}(\tau).$  In this section we
prove the following

\proclaim{Theorem A}{Let $(M,g)$ be a closed Riemannian manifold of
odd dimension,
${\Cal W}$ an $({\Cal A},\Gamma^{{\roman {op}}})$-Hilbert module of
finite type with $l = {\roman {dim}}_{{\roman {N}}} \Cal W$ and $h:M
\longrightarrow {\Bbb R}$ a Morse function.  Assume that $(M,{\Cal
W})$ is of determinant class and that $\tau =
(h,g)$ is a generalized triangulation.  Then the following statements
are true:
\roster

\item The functions $\log T_{{\roman {an}}}(h,t), \ \log T_{{\roman
{sm}}}(h,t)$ and $\log T_{{\roman {la}}}(h,t)$ admit asymptotic
expansions for $t \longrightarrow \infty.$

\item The asymptotic expansion of $\log T_{{\roman {an}}}(h,t)$ is of
the form
$$\align
\log T_{{\roman {an}}}(h,t) = & \log T_{{\roman {an}}}(h,0) - \log
T_{{\roman {met}}}(\tau)  + \tag{6.2} \\
                               & \frac{1}{2}\left(\sum_{q=0}^d
(-1)^{q+1} q \beta_q \right) (2t- \log t + \log \pi) + O(t^{-1}).
\endalign$$

\item The asymptotic expansion of $\log T_{{\roman {sm}}}(h,t)$ is of
the form
$$\align
 \log T_{{\roman {comb}}}(\tau) +  \frac{1}{2}\left(\sum_{q=0}^d
(-1)^{q+1} (q \beta_q -q m_q l)
\right) (2t- \log t + \log \pi) + o(1).  \tag{6.3}
\endalign$$
\endroster
}
\endproclaim
As argued in Introduction it suffices to prove the statements for
$\Cal W$ a free $\Cal A$-module.
We begin by deriving an alternative formula for the analytic
torsion (cf \cite{RS} ,\cite{Ch} and \cite{BFK1}).
  The space of $q$-forms can be
decomposed into orthogonal subspaces:
$$\align
\qme = & \qmep \oplus \qmem \oplus {\Cal H}_t^q \tag{6.4}
\endalign$$
where
$$\align
\qmep = & closure ( d_{q-1}(t)  \qmme) ; \tag{6.5}\\
\qmem = & closure ( d_{q}(t)^*  \qmpe) ; \tag{6.6}\\
{\Cal H}_t^q = & \{\omega \in \qme ; \Delta_q(t) \omega = 0\}. \tag{6.7}
\endalign$$
Note that the spaces $\Lambda^{\pm,q}_t(M;{\Cal E})$ are invariant
with respect to the Laplacian $ \Delta_q(t).$  Denote by $
\Delta_q^\pm(t)$ the restriction of $ \Delta_q(t)$ to
$\Lambda^{\pm,q}_t(M;{\Cal E})$ which are given by $\Delta_q^+(t) =
d_{q-1}(t) d_{q-1}(t)^* $ and $\Delta_q^-(t) = d_{q}(t)^*d_{q}(t). $
The operator $d_q(t)$ maps the space $\qmem$ injectively onto a dense
subspace of
 $\Lambda^{+,q+1}_t(M;{\Cal E})$ and it intertwines $\Delta_q^-(t) $
and $\Delta_{q+1}^+(t) .$  As a consequence, $d_q(t)$ intertwines the
spectral projectors $Q_q^-(\lambda,t)$ and $Q_{q+1}^+(\lambda,t),$
$$\align
d_q(t) Q_q^-(\lambda,t) = & Q_{q+1}^+(\lambda,t)d_q(t). \tag{6.8}
\endalign$$
This implies
$$\align
N_q^-(\lambda,t) = & \trn( Q_q^-(\lambda,t)) \\
                 = & \trn( Q_{q+1}^+(\lambda,t)) =
N_{q+1}^+(\lambda,t). \tag{6.9}
\endalign$$
Note that
both $\Delta_q^+(t) $ and $\Delta_q^-(t) $ are of determinant class i.e
$\int_{0^{+}}^{1} \log \lambda dN^{\epsilon}_A(\lambda) > - \infty, \;
\epsilon=+,-.$
Using the heat kernel representation of the zeta function we obtain
(cf \cite{Lo}):
$$\align
\int_1^\infty \frac{dx}{x} \trn \left(e^{-x \Delta_q^-(t)}\right) = &
\int_1^\infty \frac{dx}{x} \trn \left(e^{-x \Delta_{q+1}^+(t)}\right)
\tag{6.10.A}
\endalign$$
and for $\Re s$ sufficiently large,
$$\align
\frac{1}{\Gamma(s)}\int_0^1 dx x^{s-1} \trn\left( e^{-x
\Delta_q^-(t)}\right) = & \frac{1}{\Gamma(s)} \int_0^1
dx x^{s-1} \trn\left( e^{-x \Delta_{q+1}^+(t)}\right). \tag{6.10.B}
\endalign$$
Formulas (6.10) are now used to write
$$\align
\log T_{{\roman {an}}}(h,t) = & \frac{1}{2} \sum_{q=0}^d (-1)^q
{\roman {logdet}}_{N}
\Delta_q^-(t) \tag{6.11} \\
                            = &   \frac{1}{2} \sum_{q=0}^d (-1)^q
{\roman {logdet}}_{N}
\Delta_{q+1}^+(t).
\endalign$$
Our first goal is to compute the variation $\frac{d}{dt}\log
T_{{\roman {an}}}(h,t)$ of $\log T_{{\roman {an}}}(h,t).$  For this
purpose, we again use the heat kernel representation of the zeta
function and write (cf \cite{Lo})
$$\align
{\roman {logdet}}_{N} \Delta_q^+(t) = & - \left. \frac{\partial }{\partial
s}\right|_{s=0} \left(\frac{1}{\Gamma(s)}\int_0^1  x^{s-1} \trn \left(e^{-x
\Delta_q^+(t)}\right) dx \right) \tag{6.12} \\
                       & -\int_1^\infty  \frac{1}{x} \trn
\left(e^{-x \Delta_q^+(t)}\right)dx .
\endalign$$
To analyze the $t$-dependence of ${\roman {logdet}}_{N} \Delta_q^+(t)$ we
treat the
two terms on the right hand side of (6.12) seperately.  To illustrate
the new difficulties which arise (as compared with the classical
situation) we point out that the differentiability of $\int_1^\infty
x^{-1} \trn e^{-x \Delta_q^+(t)}dx$ with respect to $t$ is far from
being obvious.

We begin by computing $\frac{d}{dt}  ( \trn e^{-x \Delta_q^+(t)})$ and
note that $ \Delta_q^+(t):\qmep \to \qmep$ where the
space $\qmep = \overline{d_{q-1}(t) \qmme} = e^{-th}\Lambda^{+,q}(M; {\Cal E})$
depends on $t.$  It is therefore convenient to introduce $\tdq =
e^{th} \Delta^+_q(t) e^{-th} : \Lambda^{+,q}(M; {\Cal E})
\longrightarrow \Lambda^{+,q}(M; {\Cal E})$ which is isospectral with
$ \Delta^+_q(t).$  Hence, $ \trn e^{-x \Delta_q^+(t)} =  \trn e^{-x
\tdq}.$  Now one computes $\frac{d}{dt}  \trn e^{-x \tdq}$ using
Duhamel's principle and the identity $\tdq = e^{2th}(d_{q-1}d_{q-1}^*
+ 2tdh \wedge d_{q-1}^*)e^{-2th}:$
$$\align
\frac{d}{dt} \left( \trn e^{-x \tdq}\right) = & -x \trn\left(
\frac{d}{dt}\left(\tdq e^{-x \tdq}\right)\right) \\
                               = & \trn \left(2[h,-x\tdq] e^{-x\tdq}
\right) \\
   & - 2x\trn \left(e^{2th}dh \wedge d_{q-1}^*e^{-2th}
e^{-x\tdq}\right)
\endalign$$
where $[A,B]$ denotes the commutator of the two operators $A$ and $B.$
Therefore
$$\align
 \trn \left(2[h,-x\tdq] e^{-x\tdq}\right) = & 0.
\endalign$$
Using that $e^{th}d_{q-1}^* e^{-th} = d_{q-1}(t)^*$ and that $e^{-th}
e^{-x \tdq} e^{th} = e^{-x \Delta_q^+(t)}$ we obtain
$$\align
\frac{d}{dt} \left( \trn \left( e^{-x \Delta_q^+(t)}\right) \right) =
& -2x \trn \left(dh\wedge d_{q-1}(t)^*  e^{-x \Delta_q^+(t)}\right).
\endalign$$
Further observe that, despite the fact that $d_q(t):\qmep
\longrightarrow \Lambda_t^{-,q+1}(M;{\Cal E})$ is
\newline not invertible (it might not be onto) we can form
$$\align
d_{q-1}(t)^* = & d_{q-1}(t)^{-1}  d_{q-1}(t) d_{q-1}(t)^* =
d_{q-1}(t)^{-1}  \Delta_q^+(t)
\endalign$$
where the domain of definition of $ d_{q-1}(t)^{-1} $ is the range of
$ d_{q-1}(t).$  We note that
$$\align
dh \wedge d_{q-1}(t)^{-1}  \Delta_q^+(t)  = & (d_{q-1}(t) h
d_{q-1}(t)^{-1}  -h )\Delta_q^+(t).
\endalign$$
This leads to the following formula
$$\align
\frac{d}{dt} \left( \trn \left(  e^{-x \Delta_q^+(t)}\right) \right) =
& -2x\trn \left(d_{q-1}(t) h d_{q-1}(t)^{-1} \Delta_q^+(t)  e^{-x
\Delta_q^+(t)}\right) \\  & + 2x \trn \left(h \Delta_q^+(t)  e^{-x
\Delta_q^+(t)}\right).
\endalign$$
Next we observe that
$$\align
 2x \trn\left(h \Delta_q^+(t)  e^{-x \Delta_q^+(t)}\right)= & -2x
\frac{d}{dx} \left( \trn\left(h e^{-x \Delta_q^+(t)}\right) \right)
\endalign$$
and that
$$\align
\trn \left(d_{q-1}(t) h d_{q-1}(t)^{-1} \Delta_q^+(t)  e^{-x
\Delta_q^+(t)}\right)  = & \trn \left( h d_{q-1}(t)^{-1} \Delta_q^+(t)  e^{-x
\Delta_q^+(t)} d_{q-1}(t)\right) \\
                = & \trn \left( h  \Delta_{q-1}^-(t)  e^{-x
\Delta_{q-1}^-(t)}\right) \\
                = & -\frac{d}{dx} \left( \trn \left( h  e^{-x
\Delta_{q-1}^-(t)}\right) \right).
\endalign$$
We have therefore proved that
$$\align
\frac{d}{dt} \left( \trn \left(  e^{-x \Delta_q^+(t)}\right) \right) =
2x\frac{d}{dx} \left( \trn \left( h  e^{-x
\Delta_{q-1}^-(t)}\right) \right) -2x \frac{d}{dx} \left( \trn \left(
h  e^{-x \Delta_{q}^+(t)}\right) \right).
\endalign$$
This leads to
$$\align
\frac{1}{2} \sum_{q=0}^d (-1)^{q+1} \frac{d}{dt} \left( \trn \left(
e^{-x \Delta_q^+(t)}\right) \right) = &
x \sum_{q=0}^d (-1)^{q+1} \frac{d}{dx} \left( \trn \left( h  e^{-x
\Delta_{q-1}^-(t)}\right) \right) \\  & -x \sum_{q=0}^d (-1)^{q+1}
\frac{d}{dx} \left( \trn \left( h  e^{-x \Delta_{q}^+(t)}\right)
\right) \\
  = &  x \frac{d}{dx} \left( \sum_{q=0}^d (-1)^{q}  \trn \left( h  e^{-x
\Delta_{q}(t)}({\roman {Id}} - Q_q(0,t)) \right) \right). \tag{6.13}
\endalign $$
The above formula is used to prove that $-\int_1^\infty \frac{1}{x}
\frac{1}{2} \sum_q (-1)^{q+1} \trn \left(e^{-x \Delta_q^+(t)}\right) dx$ has a
continuous derivative with respect to $t.$

By the Leibniz rule for improper integrals, it suffices to verify that
\newline $f(x,t) = - \frac{1}{x} \frac{1}{2} \sum_q (-1)^{q+1} \trn \left(e^{-x
\Delta_q^+(t)}\right) $ and $\frac{\partial f}{\partial t}(x,t)$ are both
continuous and the integrals $\int_1^\infty f(x,t) dx$ and
$\int_1^\infty \frac{\partial f}{\partial t}(x,t) dx$  both converge
uniformly with respect to $t$ ($t$ varying in a compact interval).
Clearly $f(x,t)$ is continuous and, by the above formula,
$$\frac{\partial f}{\partial t}(x,t) = -\frac{d}{dx}  \sum_{q=0}^d
(-1)^{q} \trn \left(he^{-x \Delta_q(t)}({\roman {Id}} - Q_q(0,t))\right). $$
The uniform convergence of the integrals $\int_1^\infty f(x,t) dx$ and
$\int_1^\infty \frac{\partial f}{\partial t}(x,t) dx$ follows from

\proclaim{Lemma 6.1}{Let $I$ be an arbitrary compact interval contained in
$[0,\infty).$  Then
\roster
\item  $ \lim_{x \longrightarrow \infty} \int_1^x \frac{1}{x}
 \trn (e^{-x \Delta_q^+(t)}) dx$ converges uniformly for $t \in I.$

\item  $\lim_{x \longrightarrow \infty} \int_1^x
\frac{\partial f}{\partial t}(x,t)dx $ converges uniformly for $t \in
I.$
\endroster
}
\endproclaim

{\sl Proof }  (1)  Note that the integrand $\frac{1}{x} \trn \left(e^{-x
\Delta_q^+(t)}\right)$ is positive.  Therefore
$$\multline
0 \leq \int_u^\infty \frac{1}{x} \trn \left(e^{-x \Delta_q^+(t)}\right) dx =
\int_{0^+}^\infty dN_{\Delta_q^+(t)}(\mu) \int_u^\infty
\frac{1}{x}e^{-\mu x} dx \\
  \leq  \int_{u^{-\frac{1}{2}}}^\infty dN_{\Delta_q^+(t)}(\mu) \frac{e^{-\mu
u}}{\mu u} + \int_{0^+}^{u^{-\frac{1}{2}}} dN_{\Delta_q^+(t)}(\mu) \left(
\log{(\mu u)} e^{-\mu u}  + \int_{\mu u}^\infty e^{-s}\log s ds \right) \\
  \leq  \frac{e^{-u^{\frac{1}{2}}}}{u^{\frac{1}{2}}}\int_{0^+}^\infty
dN_{\Delta_q^+(t)}(\mu) + C \int_{0^+}^{u^{-\frac{1}{2}}}
dN_{\Delta_q^+(t)}(\mu)
\endmultline  \tag{6.14} $$
where $C>0$ is a bound for the function $ \log{(\mu u)}
e^{-\mu u} + \int_{\mu u}^\infty e^{-s}(\log s) ds$ for $\mu$ in the
interval $[0,1].$  Statement (1) follows from (6.14) for $u
\longrightarrow \infty$ and from the fact that
$N_{\Delta_q^+(t)}(\mu)$ is right continuous with respect to $\mu,$
uniformly in $t$ for $t$ in $I.$

(2)  From (6.13), $\int_1^u \frac{\partial f}{\partial t}(x,t) =
-  \sum_{q=0}^d (-1)^{q} \trn (he^{-x \Delta_q^+(t)}({\roman {Id}} -
Q_q(0,t)))|_1^u. $  Therefore it suffices to prove that, for $0 \leq q
\leq d $ and uniformly for $t$ in $I,$
$$\align
\lim_{x \longrightarrow \infty}  \trn \left(he^{-x
\Delta_q(t)}({\roman {Id}} - Q_q(0,t))\right) = 0. \tag{6.15}
\endalign$$
This can be seen as follows:
$$\align
| \trn \left(he^{-x \Delta_q(t)}({\roman {Id}} - Q_q(0,t))\right)| \leq & \| h
\|_{L^\infty} \int_{0^+}^{\infty} e^{-x\lambda}
dN_q(\lambda,t) \\
   \leq &  \| h \|_{L^\infty} \left( \int_{0^+}^{x^{-\frac{1}{2}}}
 e^{-x\lambda}
dN_q(\lambda,t) + \int_{x^{-\frac{1}{2}}}^\infty
e^{-x\lambda}  dN_q(\lambda,t) \right) \\
  \leq &  \| h \|_{L^\infty} \left( (N_q(x^{-\frac{1}{2}},t) -
N_q(0,t))  + e^{-x^{\frac{1}{2}}}\int_{0}^\infty
dN_q(\lambda,t) \right)
\endalign$$
and (6.15) follows from the fact that $N_q(\lambda,t)$ is right
continuous with respect to $\mu,$ uniformly for $t$ in $I.$

We have shown that $-\int_1^\infty   \frac{1}{x} \frac{1}{2} \sum_q
(-1)^{q+1} \trn \left(e^{-x \Delta_q^+(t)}\right)dx$ has a continuous
derivative with respect to $t:$
$$\multline
\frac{d}{dt}\left(-\int_1^\infty   \frac{1}{x} \frac{1}{2} \sum_q
(-1)^{q+1} \trn \left(e^{-x \Delta_q^+(t)}\right)dx \right) \\ =  \sum_{q=0}^d
(-1)^q  \trn \left(he^{-x \Delta_q(t)}({\roman {Id}} - Q_q(0,t))\right).
\diamondsuit
\endmultline \tag{6.16} $$
We now analyze the $t$-derivative of $\frac{\partial }{\partial
s}|_{s = 0} \frac{1}{\Gamma(s)} \int_0^1 x^{s-1}  \frac{1}{2} \sum_q
(-1)^{q+1} \trn (e^{-x \Delta_q^+(t)})dx .$

For $\Re s > \frac{d}{2},$ we integrate by parts in (6.12) to obtain
$$
\multline
-\frac{d }{dt} \left( \frac{1}{\Gamma(s)} \int_0^1 x^{s-1}  \frac{1}{2} \sum_q
(-1)^{q+1} \trn \left(e^{-x \Delta_q^+(t)}\right)dx\right) \\
 =  -\frac{1}{\Gamma(s)} \int_0^1 x^{s}\frac{d}{dx} \left( \sum_{q=0}^d
(-1)^q  \trn \left(
he^{-x \Delta_q(t)}({\roman {Id}} - Q_q(0,t))\right) dx\right) \\
  =  - \frac{1}{\Gamma(s)} \sum_{q=0}^d (-1)^q  \trn \left(he^{-
\Delta_q(t)}({\roman {Id}} - Q_q(0,t))\right) \\
 + \frac{s}{\Gamma(s)}
\int_0^1 x^{s-1}\frac{d}{dx}  \sum_{q=0}^d (-1)^q  \trn (he^{-x
\Delta_q(t)}({\roman {Id}} - Q_q(0,t))) dx.
\endmultline
$$
The last two functions both have a meromorphic extension to the $s$-plane
which is regular at $s=0.$  With $\frac{1}{\Gamma(s)} =
\frac{s}{\Gamma(s+1)}$ and $\Gamma(1) = 1,$ we obtain
$$\multline
\frac{d }{dt}\left(-\frac{\partial }{\partial s}|_{s=0}
\frac{1}{\Gamma(s)} \int_0^1 x^{s-1}  \frac{1}{2} \sum_q (-1)^{q+1}
\trn \left(e^{-x \Delta_q^+(t)}\right)dx \right) \\ =  - \sum_{q=0}^d
(-1)^q  \trn \left(he^{-
\Delta_q(t)}({\roman {Id}} - Q_q(0,t))\right) \\
  + {\roman
{F.p.}}_{s=0}\frac{1}{\Gamma(s)}\int_0^1 x^{s-1} \sum_{q=0}^d (-1)^q
\trn \left(he^{-x \Delta_q(t)}({\roman {Id}} - Q_q(0,t))\right) dx.
\endmultline \tag{6.17}$$
Combing (6.16) and (6.17) we conclude that $\log T_{\roman {an}}(h,t)$
is continuous and has a continuous derivative with respect to $t.$
Moreover,
$$\align
\frac{d}{dt}\log T_{\roman {an}}(h,t) = &
{\roman {F.p.}}_{s=0}\frac{1}{\Gamma(s)}\int_0^1 x^{s-1}
\sum_{q=0}^d (-1)^q  \trn \left(he^{-x \Delta_q(t)}({\roman {Id}} -
Q_q(0,t))\right) dx. \tag{6.18}
\endalign$$
Next, $ \trn \left(he^{-x \Delta_q(t)}({\roman {Id}} -
Q_q(0,t))\right) =  \trn \left(he^{-x \Delta_q^+(t)}\right) -  \trn
\left(h Q_q(0,t)\right).$  Further, $${\roman
{F.p.}}_{s=0}\frac{1}{\Gamma(s)}\int_0^1 x^{s-1} dx = {\roman
{F.p.}}_{s=0}\frac{s}{\Gamma(s+1)}\frac{1}{s} = 1$$and the heat kernel
expansion for the Schwartz kernel $K_q(y,y',x,t)$ of
$e^{-x\Delta_q(t)}$ on the diagonal $y = y'$ is of the form
$$\align
K_q(y,y',x,t) = & \sum_{j=0}^{N-1} x^{\frac{j-d}{2}} l_j(y,t) +
O(x^{\frac{N-d}{2}},t) \tag{6.19}
\endalign$$
where $l_j(y,t)$ are densities defined on $M$ with values in ${\Cal
B}.$  By a standard parity argument one concludes that $l_d(\cdot,t) =
0$ and argues as in the classical case to conclude that
$$\multline
{\roman {F.p.}}_{s=0}\frac{1}{\Gamma(s)}\int_0^1 x^{s-1}
\sum_{q=0}^d (-1)^q  \trn (he^{-x \Delta_q(t)}({\roman {Id}} -
Q_q(0,t))) dx \\ =   \sum_{q=0}^d (-1)^{q+1}  \trn (hQ_q(0,t))
=   \sum_{q=0}^d (-1)^{q+1}  \trn \left(Q_q(0,t) h
Q_q(0,t)\right).
\endmultline \tag{6.20}  $$
We have proved the following
\proclaim{Proposition 6.2}{ $\frac{d}{dt}\log T_{\roman {an}}(h,t) =
\sum_{q=0}^d (-1)^{q+1}  \trn \left(Q_q(0,t) h Q_q(0,t)\right) .$}
\endproclaim

Next, we express the terms $  \trn \left(Q_q(0,t) hQ_q(0,t)\right) $ in a more
explicit way.  It is convenient to introduce $P_q(t) = Q_q(0,t).$
Consider $K_q(t) :{\Cal H}^q_t(M;{\Cal E}) \longrightarrow \qh$
defined by
$$\align
K_q(t)(\omega) := & P_q(0) e^{th} \omega. \tag{6.21}
\endalign$$
Using the decomposition ($\omega \in \qh$) $e^{-th}\omega =
e^{-th}\omega_+(t) + \omega_0(t) \in \qmep \oplus \qht$ where
$\omega_+(t) \in \qmep$ and $\omega_0(t) \in \qht,$ one verifies that
$P_q(t) e^{-th}$ is the right inverse of $K_q(t).$  Therefore,
$K_q(t)$ is an isomorphism.  This implies that $K_q'(t) = (K_q(t)
K_q(t)^*)^{\frac{1}{2}} $ is a selfadjoint, positive, ${\Cal
A}$-linear operator on $\qh$ and thus admits a determinant with
${\roman {det}}_{N} K_q'(t) > 0.$  Note that $K_q(t)^*$ is given by
$P_q(t)e^{th}$
and thus $K_q(t)K_q(t)^*$ can be written as
$$\align
K_q(t)K_q(t)^* = & P_q(0) e^{th} P_q(t) e^{th} P_q(0). \tag{6.22}
\endalign$$

\proclaim{Lemma 6.3}{$\trn \left(P_q(t) h P_q(t) \right) = \frac{d}{dt}
{\roman {logdet}}_{N} (K_q(t)
K_q(t)^*)^{\frac{1}{2}}.$}
\endproclaim

{\sl Proof}  $\ $ Using Proposition 1.9 we note that
$$\align
\frac{d}{dt} {\roman {logdet}}_{N} (K_q(t)K_q(t)^*)^{\frac{1}{2}} = &
\frac{1}{2}\frac{d}{dt} {\roman {logdet}}_{N} (K_q(t) K_q(t)^*)  \tag{6.23}  \\
  = & \frac{1}{2} \trn \left( \frac{d}{dt} (K_q(t) K_q(t)^*) (K_q(t)
K_q(t)^*)^{-1} \right) .
\endalign$$
Using (6.22) and writing $\dot{P}_q(t) = \frac{d}{dt} P_q(t),$ we obtain
$$\align
\frac{d}{dt} (K_q(t)K_q(t)^*) = &  P_q(0) h e^{th} P_q(t) e^{th} P_q(0)
+  P_q(0) e^{th} \dot{P}_q(t) e^{th} P_q(0) + \\
   & + P_q(0) e^{th} P_q(t)h e^{th} P_q(0). \tag{6.24}
\endalign$$
To compute $\dot{P}_q(t) = \frac{d}{dt}(P_q(t)^2) = \dot{P}_q(t) P_q(t) +
P_q(t) \dot{P}_q(t) $ we consider the orthogonal decomposition $\qme =
\qht \oplus \qmep \oplus \qmem.$  An element $\omega \in \qme$ can be
uniquely written as
$$\align
\omega = & \omega_0(t) +e^{-th} \omega_+(t) + e^{th} \omega_-(t) \tag{6.25}
\endalign$$
where $\omega_\pm(t) \in \Lambda^{\pm,q}(M;{\Cal E})$ and $\omega_0(t)
= P_q(t) \omega(t).$  We conclude that
$$\align
0 = \frac{d}{dt} \omega = &  \dot{\omega}_0(t) +e^{-th}
\dot{\omega}_+(t) + e^{th} \dot{\omega}_-(t) \\
     &  -he^{-th} \omega_+(t) + he^{th} \omega_-(t).  \tag{6.26}
\endalign$$
Note that $\dot{\omega}_\pm(t) \in \Lambda^{\pm,q}(M;{\Cal E})$ and
therefore $e^{-th}\dot{\omega}_+(t) \in \qmep$ and $e^{th}\dot{\omega}_-(t)
\in \qmem.$  Applying $P_q(t)$ to (6.26) leads to
$$\align
0 = & P_q(t)  \dot{\omega}_0(t)  -P_q(t) he^{-th} \omega_+(t) + P_q(t)
he^{th} \omega_-(t)  .
\endalign$$
Denoting by $P_q^\pm (t)$ the orthogonal projectors $\qme
\longrightarrow \Lambda^{\pm,q}_t(M;{\Cal E})$ we therefore obtain
$$\align
P_q(t) \dot{P}_q(t) = & P_q(t)h P_q^+(t) - P_q(t) h P_q^-(t). \tag{6.27}
\endalign$$
Observe that the projectors $P_q(t)$ and therefore $\dot{P}_q(t)$ are
selfadjoint to conclude that
$$\align
\dot{P}_q(t) P_q(t) = & P_q^+(t) h P_q(t) - P_q^-(t) h P_q(t). \tag{6.28}
\endalign$$
Combining (6.27) and (6.28) we obtain
$$\align
 \dot{P}_q(t) = & \dot{P}_q(t) P_q(t) + P_q(t) \dot{P}_q(t) \\
  = & P_q(t)h P_q^+(t) + P_q^+(t)h P_q(t) - P_q(t) h P_q^-(t) - P_q^-(t) h
P_q(t). \tag{6.29}
\endalign$$
We apply formula (6.29) to rewrite (6.24),
$$\align
P_q(0) e^{th} \dot{P}_q(t) e^{th} P_q(0) = &  P_q(0)  e^{th} P_q(t)h P_q^+(t)
e^{th} P_q(0)  +  P_q(0)  e^{th} P_q^+(t)h P_q(t)e^{th}P_q(0) \\
  & - P_q(0)  e^{th} P_q(t)h P_q^-(t)e^{th}P_q(0) - P_q(0)  e^{th}
P_q^-(t)h P_q(t)e^{th}P_q(0) . \tag{6.30}
\endalign$$
To simplify (6.30), notice that $d_q(t)^* = e^{th}d_q^* e^{-th}$ and
therefore $e^{th} \qh \subset \qht \oplus \qmem$ which implies that $ P_q^+(t)
e^{th} P_q(0) = 0.$  Taking the adjoint, we conclude that $ P_q(0)
e^{th} P_q^+(t) = 0.$  Thus, the first two terms on the right hand
side of (6.30) are zero.  Moreover $ P_q^-(t)e^{th}P_q(0) = ({\roman
{Id}} - P_q(t))e^{th}P_q(0)$ as well as $ P_q(0)e^{th}P_q^-(t) =
P_q(0)e^{th}({\roman {Id}} - P_q(t)).$  Applying these considerations
to (6.24) yields
$$\align
\frac{d}{dt} (K_q(t)K_q(t)^*) = &  P_q(0) h e^{th} P_q(t) e^{th} P_q(0)
-  P_q(0) e^{th} P_q(t)h ({\roman {Id}} - P_q(t)) e^{th} P_q(0) \\
  & +  P_q(0) e^{th} P_q(t)h e^{th} P_q(0) -
P_q(0) e^{th}  ({\roman {Id}} - P_q(t))h P_q(t) e^{th} P_q(0)\\
 = &2 P_q(0)  e^{th} P_q(t)h P_q(t) e^{th} P_q(0) \\
 = & 2 K_q(t) P_q(t)h P_q(t)  K_q(t)^*. \tag{6.31}
\endalign$$
Substituting (6.31) into (6.23) leads to
$$\align
\frac{d}{dt} {\roman {logdet}}_{N}  (K_q(t)K_q(t)^*)^{\frac{1}{2}}
= & \trn ( K_q(t)
P_q(t)h P_q(t)  K_q(t)^*) (K_q(t)K_q(t)^*)^{-1}) \\
  = &  \trn ( K_q(t)P_q(t)h P_q(t)  K_q(t)^{-1}) \\
  = &  \trn ( P_q(t)h P_q(t) )
\endalign$$
which concludes the proof of the lemma. $\diamondsuit$

Using that $K_q(0) = {\roman {Id}}$ and therefore that ${\roman {det}}_{N}
(K_q(0)K_q(0)^*)^{\frac{1}{2}} = 1, $ Proposition 6.2 together with
Lemma 6.3 lead to
$$\align
\log T_{{\roman {an}}}(h,t) = & \log T_{{\roman {an}}}(h,0) + \sum_{q=0}^d
(-1)^{q+1} \int_0^t \frac{d}{dt} {\roman {logdet}}_{N}
(K_q(t)K_q(t)^*)^{\frac{1}{2}} dt \\
 = & \log T_{{\roman {an}}}(h,0) + \sum_{q=0}^d
(-1)^{q+1} {\roman {logdet}}_{N}
(K_q(t)K_q(t)^*)^{\frac{1}{2}} . \tag{6.32}
\endalign$$
In section 4.1 we introduced the ${\Cal A}$-linear isomorphisms
$$\align
 & \theta_q :{\roman {Null}}\Delta_q^{{\roman {comb}}} \longrightarrow
\qh
\endalign$$
and the metric part of the Reidemeister torsion $T_{{\roman
{met}}}(\tau)$ defined by $$ \log T_{{\roman {met}}}(M,g,{\Cal W},\tau) =
\frac{1}{2} \sum_{q=0}^d
(-1)^{q}  {\roman {logdet}}_{N}(\theta_q^* \theta_q).$$By
applying the analysis of the Witten deformation of the deRham complex
by Helffer-Sj\"ostrand (cf section 5) we obtain

\proclaim{Lemma 6.4}{For $t$ sufficiently large, the following
statements hold:

$$\align {\roman {logdet}}_{N}(K_q(t)K_q(t)^*)^{\frac{1}{2}} = &
{\roman {logdet}}_{N} (\theta_q^* \theta_q)^{\frac{1}{2}} \\
   &  + q
\beta_q t +\beta_q\left(\frac{d-2q}{4}\right) \log \left(\frac{t}{\pi}
\right) + O(t^{-1}) \tag{6.33} \endalign$$
and
$$\align \sum_{q=0}^d (-1)^{q+1}{\roman
{logdet}}_{N}(K_q(t)K_q(t)^*)^{\frac{1}{2}} = & - \log T_{{\roman
{met}}}(M,g,{\Cal W},\tau) +  \sum_{q=0}^d (-1)^{q+1} q\beta_q t \\
 & + \sum_{q=0}^d (-1)^{q+1} \beta_q \frac{d-2q}{4} \log
\left(\frac{t}{\pi}  \right) + O(t^{-1}). \tag{6.34} \endalign$$

}
\endproclaim

{\sl Proof } $\ $ Summing with respect to $q,$ statement (6.34) follows
directly from statement (6.33) and the definition $ \log T_{{\roman
{met}}}(M,g,{\Cal W},\tau) = \frac{1}{2} \sum_{q=0}^d
(-1)^{q}  {\roman {logdet}}_{N}(\theta_q^* \theta_q).$  To prove
(6.33) we use Stokes'
 theorem to write
$$\align
K_q(t) = & \theta_q K_q'
\endalign$$
where $K_q'(t) = \pi_q K_q''(t) I_q(t),$ the map $\pi_q : {\Cal
C}^q \longrightarrow {\roman {Null}}\Delta_q^{{\roman {comb}}}$
denotes orthogonal projection, $I_q(t) : \qht \longrightarrow
\Lambda^q(M; {\Cal E})_{{\roman {sm}}} $ denotes inclusion and
$$\align
K_q''(t) :\Lambda^q(M; {\Cal E})_{{\roman {sm}}} \longrightarrow {\Cal
C}^q ,   \ \  & \omega(t) \longrightarrow  {Int}^{(q)}(e^{th}\omega(t)).
\endalign$$
Note that
$$\align
{\roman {logdet}}_{N}(K_q(t)K_q(t)^*)^{\frac{1}{2}} = &
{\roman {logdet}}_{N} (\theta_q^* \theta_q)^{\frac{1}{2}} +
 \frac{1}{2}{\roman {logdet}}_{N}(K_q'(t)K_q'(t)^*).\tag{6.35}
\endalign$$
To compute ${\roman {logdet}}_{N}(K_q'(t)K_q'(t)^*)$ introduce the
scaled version of $K_q''(t),$
$$\align
K_q'''(t) := &  \left(\frac{\pi}{t}\right)^{\frac{(d-2q)}{4}} e^{-tq}
K_q''(t) \tag{6.36}
\endalign$$
Then, with $K_q''''(t) := \pi_q(t) K_q'''(t) I_q(t),$
$$\align
{\roman {logdet}}_{N}(K_q'(t)K_q'(t)^*)  = & \beta_q
\frac{(d-2q)}{4} \log \left( \frac{t}{\pi} \right) + \beta_q
2qt + {\roman {logdet}}_{N}(K_q''''(t)K_q''''(t)^*). \tag{6.37}
\endalign$$
In view of Corollary 5.8, we can apply Proposition 1.17 to the map
$K_q''''(t)$ and therefore obtain
$$\align
 {\roman {logdet}}_{N}(K_q''''(t)K_q''''(t)^*) = & O(t^{-1}). \tag{6.38}
\endalign$$
Formula (6.33) now follows from (6.35), (6.37) and (6.38).

\proclaim{Lemma 6.5}{For $t \longrightarrow \infty,$
$$\align
\log T_{{\roman {sm}}}(h,t) = & \log T_{{\roman {comb}}}(\tau) \\
  &   +\frac{1}{2} \left(\sum_{q=0}^d(-1)^{q+1} q (\beta_q - m_q l)\right)
(2t -\log t + \log \pi) +o(1).
\endalign$$
}
\endproclaim
{\sl Proof}  $\ $ Recall that $\log T_{{\roman {sm}}}(h,t)$ is a real
number defined by
$$\align
\log T_{{\roman {sm}}}(h,t) = & \frac{1}{2}
\left(\sum_{q=0}^d(-1)^{q+1} q {\roman {logdet}}_{{\roman N}}
\Delta_q(t)_{{\roman {sm}}}\right) \tag{6.46}
\endalign$$
and that, for any $0 < C < \infty,$
$$\align
{\roman {logdet}}_{{\roman N}} \Delta_q(t)_{{\roman {sm}}}  = &
-\left. \frac{\partial }{\partial s} \right| _{s=0} \frac{1}{\Gamma(s)}
 \int_0^C
\frac{dx}{x} x^s \left( \trn \left( e^{-x  \Delta_q(t)_{{\roman
{sm}}}}\right) - \beta_q \right) \\
  & - \int_C^\infty
\frac{dx}{x}  \left( \trn \left( e^{-x  \Delta_q(t)_{{\roman
{sm}}}}\right) - \beta_q \right) \tag{6.47}
\endalign$$
where $\beta_q = {\roman {dim}}_{\roman N} {\roman {Null}}
\Delta_q^{{\roman {comb}}} =  {\roman {dim}}_{\roman N} {\roman {Null}}
\Delta_q(t).$   In view of Theorem 5.4, (5), introduce
$\widetilde{\Delta}_q(t)_{{\roman {sm}}} : = \frac{\pi}{t} e^{2t}
\Delta_q(t)_{{\roman {sm}}} .$  By a change of variable of
integration, $y = \frac{t}{\pi} e^{-2t}x,$ we obtain
$$\align
 \int_{C\frac{\pi}{t} e^{2t}}^\infty
\frac{dx}{x}  \left( \trn \left( e^{-x  \Delta_q(t)_{{\roman
{sm}}}}\right) - \beta_q \right) = &   \int_{C}^\infty
\frac{dy}{y}  \left( \trn \left( e^{-y  \widetilde{\Delta}_q(t)_{{\roman
{sm}}}}\right) - \beta_q \right)\\
 = & \int_{0^+}^\infty \left( \int_C^\infty \frac{1}{y} e^{-\mu y} dy
\right) dF_{\widetilde{\Delta}_q(t)_{{\roman {sm}}}}(\mu) \\
 = & \int_{0^+}^\infty \left( \int_C^\infty  e^{-\mu y} dy
\right) F_{\widetilde{\Delta}_q(t)_{{\roman {sm}}}}(\mu) d\mu
\endalign$$
where $F_{\widetilde{\Delta}_q(t)_{{\roman {sm}}}}$ has been introduced
in (1.9) and the last equality follows from integration by parts.
 From Theorem 5.7,(5) and Proposition 1.18 we conclude that there exists
 $t_0 > 0$ such that, for $t \geq t_0$ and $ 0 \leq q \leq d, \
F_{\widetilde{\Delta}_q(t)_{{\roman {sm}}}} (\mu) \leq
F_{\Delta_q^{{\roman {comb}}}}(10 \mu)  .$
For $t \geq t_0$ and $ 0 \leq q \leq d,$
the above computations lead to
$$\align
 \int_{C\frac{\pi}{t}e^{2t}}^\infty
\frac{dx}{x}  \left( \trn \left( e^{-x  \Delta_q(t)_{{\roman
{sm}}}}\right) - \beta_q \right) \leq &  \int_{\frac{C}{10}}^\infty
\frac{dx}{x}  \left( \trn \left( e^{-x  \Delta_q^{{\roman
{comb}}}}\right) - \beta_q \right) . \tag{6.48}
\endalign$$
Taking into account that $(M,{\Cal W})$ is of determinant class, one
can choose $C > 0$ such that
$$\align
 \int_{\frac {C}{10}}^\infty
\frac{dx}{x}  \left( \trn \left( e^{-x  \Delta_q^{{\roman
{comb}}}}\right) - \beta_q \right) \leq & \epsilon . \tag{6.49}
\endalign$$
Therefore, for all $t \geq t_0, \ 0 \leq q \leq l,$ it suffices to
consider
$$\align
- \left. \frac{\partial }{\partial s} \right| _{s=0}
 \frac{1}{\Gamma(s)}
\int_0^{C\frac{\pi}{t}e^{2t}} \frac{dx}{x} x^s \left( \trn \left( e^{-x
\Delta_q(t)_{{\roman {sm}}}}\right) - \beta_q \right).  & \tag{6.50}
\endalign$$
By Theorem 5.7, $\Delta_q(t)_{{\roman {sm}}},$ when expressed in
convenient coordinates, is of the form
\newline $\frac{t}{\pi}
e^{-2t}(\Delta_q^{{\roman {comb}}} +O(t^{-\frac{1}{2}})).  $  By a
change of variable of integration, $y = \frac{t}{\pi}e^{-2t}x,$ the
expression (6.50) can be written as a sum of two terms, $I_q +II_q,$
where
$$\align
I_q = & - (\log \pi - \log t + 2t ) III_q(s)_{|_{s=0}} \\
II_q = & - \left. \frac{d}{ds} \right| _{s=0} III_q(s)
\endalign$$
and
$$\align
III_q(s) = &\frac{1}{\Gamma(s)} \int_0^{C}
\frac{dy}{y}y^s  \left( \trn \left( e^{-y  (\Delta_q^{{\roman
{comb}}} +O(t^{-\frac{1}{2}})) }\right) - \beta_q \right).
\endalign$$
We first evaluate $III_q(s)$ at $s=0.$  Use $\frac{1}{\Gamma(s)} =
\frac{s}{\Gamma(s+1)} $ and integrate by parts to obtain, for
arbitrary $\delta >0,$
$$\align
\left. III_q(s) \right| _{s=0} = &\frac{1}{\Gamma(s+1)} \int_0^{\delta}
dy \frac{d}{dy}(y^s)  \left. \left( \trn \left( e^{-y  (\Delta_q^{{\roman
{comb}}} +O(t^{-\frac{1}{2}})) }\right) - \beta_q \right) \right| _{s=0} \\
  = & \trn \left( e^{-\delta  (\Delta_q^{{\roman
{comb}}} +O(t^{-\frac{1}{2}})) }\right) -\beta_q  \\
  + & \int_0^{\delta} dy    \trn \left((\Delta_q^{{\roman
{comb}}} +O(t^{-\frac{1}{2}})) e^{-y  (\Delta_q^{{\roman
{comb}}} +O(t^{-\frac{1}{2}})) }\right) . \tag{6.51}
\endalign$$
Taking the limit as $\delta \longrightarrow 0$ in (6.51) leads to
$$\align
\left. III_q(s) \right| _{s=0} = & m_q l - \beta_q. \tag{6.52}
\endalign$$
To compute $II_q = - \left. \frac{d}{ds} \right| _{s=0} III_q(s)$ recall that
$$\align
\ldetn \Delta_q^{{\roman {comb}}} = & -\left. \frac{\partial }{\partial
s} \right| _{s=0} \frac{1}{\Gamma(s)} \int_0^{C}
\frac{dy}{y}y^s  \left( \trn \left( e^{-y  (\Delta_q^{{\roman
{comb}}} ) }\right) - \beta_q \right) \\
  & -\int_C^{\infty}
\frac{dy}{y}  \left( \trn \left( e^{-y  (\Delta_q^{{\roman
{comb}}} ) }\right) - \beta_q \right)
\endalign$$
and use the estimate $(0 \leq y \leq C)$
$$\align
\left|\trn \left( e^{-y  (\Delta_q^{{\roman {comb}}})}  \right) -
 \trn \left( e^{-y  (\Delta_q^{{\roman {comb}}} +O(t^{-\frac{1}{2}})) }\right)
\right| \leq & y O \left(t^{-\frac{1}{2}} \right) \tag{6.53}
\endalign$$
to conclude, together with (6.49) that
$$\align
\left| II_q - \ldetn \Delta_q^{{\roman {comb}}}\right| \leq & \epsilon
+ O \left(t^{-\frac{1}{2}} \right). \tag{6.54}
\endalign$$
Combining (6.46)-(6.48) and (6.52)-(6.54) we conclude that for given
$\epsilon > 0,$ there exists $t_\epsilon >0$ so that for all $t >
t_\epsilon, $
$$\align
\left| \log T_{{\roman {sm}}}(h,t) - \log T_{{\roman {comb}}}(\tau) -
\frac{1}{2} \sum_{q=0}^d(-1)^{q+1} q (\beta_q - m_q l)
(2t -\log \frac{t}{\pi}) \right| \leq & 3\epsilon.\ \  \diamondsuit
\endalign$$
{\sl Proof of Theorem A}  First note that $\log T_{{\roman {la}}}(h,t)
= \log T_{{\roman {an}}}(h,t) - \log T_{{\roman {sm}}}(h,t).$
Therefore the asymptotic expansion of $\log T_{{\roman {la}}}(h,t)$ is
obtained from the expansions of $\log T_{{\roman {an}}}(h,t)$ and
$\log T_{{\roman {sm}}}(h,t).$  The asymptotic expansion for $\log
T_{{\roman {an}}}(h,t) $ follows from (6.32) and Lemma 6.4 together
with the fact that, since $d$ is odd, $\chi(M,\tau) = \sum_{q=0}^d
(-1)^q \beta_q = 0.$  The asymptotic expansion for $\log T_{{\roman
{sm}}}(h,t) $ is contained in Lemma 6.5.

\subheading{6.2  Comparison theorem for Witten's deformation of the
analytic torsion}

The family of operators $\Delta_q(t)$ is a family with parameter of
order 2 and weight 1 (cf \cite{Sh1} and section 3).
The family $\Delta_q(t)$ fails to be elliptic
with parameter precisely at the critical points of the Morse function
$h.$  We can therefore use the Mayer-Vietoris type formula for
determinants (cf  section 3) to localize the failure of
the family $\Delta_q(t)$ to be elliptic with parameter and to obtain a
relative result which compares the asymptotic expansions of Witten's
deformation of the analytic torsion corresponding to two different
systems $(M^d,h,g,{\Cal W})$ and $({M'}^d ,h',g',{\Cal W})$ where the manifolds
$M$ and $M'$ have the same fundamental group $\Gamma$ and $(h,g)$
and $(h',g')$ are generalized triangulations.

\proclaim{Theorem B}{Let $d$ be odd.  Suppose that $\tau = (h,g)$ and
$\tau' = (h',g')$ are generalized triangulations with $\# {\roman
{Cr}}_q(h) = \# {\roman {Cr}}_q(h') , \ (0 \leq q \leq d),$ and that
$(M,{\Cal W})$ and $(M',{\Cal W})$ are of determinant class.  Then the
following statements hold:
\roster

\item  The free term ${\roman {FT}}(\log T_{{\roman {la}}}(h,t) - \log
T_{{\roman {la}}}(h',t) )$ of the asymptotic expansion of
\newline $\log
T_{{\roman {la}}}(h,t) - \log T_{{\roman {la}}}(h',t)$ is given by
$$\multline
{\roman {FT}}(\log T_{{\roman {la}}}(h,t) - \log
T_{{\roman {la}}}(h',t) ) \\ =  \int_{M \setminus {\roman {Cr}}(h)}
a_0(h,\epsilon = 0) -  \int_{M' \setminus {\roman {Cr}}(h')}
a_0(h',\epsilon = 0)
\endmultline \tag{6.55}$$
where the densities $a_0(h,\epsilon = 0) $ and $a_0(h',\epsilon = 0) $
are smooth forms of degree $d$ and are given by explicit local
formulas and the difference appearing in the right hand side
of (6.55) is taken in the sense explained in the remark below, (6.56).

\item  Due to the assumption that $d$ is odd,
$$\align
a_0(h,\epsilon = 0,x) + a_0(d-h,\epsilon = 0, x) = & 0.
\endalign$$
\endroster
}
\endproclaim

{\bf {Remark}}  The integral $ \int_{M \setminus {\roman {Cr}}(h)}
a_0(h,\epsilon = 0)$ need not be convergent and the difference on the
right hand side of (6.55) should be understood in the following sense:
In view of the definition of a generalized triangulation, there exist
neighborhoods $V$ of ${\roman {Cr}}(h)$ and $V'$ of ${\roman
{Cr}}(h')$ and a smooth bundle isomorphism $F:{\Cal E}_{|_V}
\longrightarrow {\Cal E}'_{|_{V'}}$ so that $f$ and $F$
intertwine the functions $h$ and $h',$ the metrics $g$ and $g'$ and
the Laplace operators $\Delta_q$ and $\Delta_q'.$  Define
$$\multline
 \int_{M \setminus {\roman {Cr}}(h)}
a_0(h,\varepsilon = 0) -  \int_{M' \setminus {\roman {Cr}}(h')}
a_0(h',\varepsilon = 0)\\ = \int_{M \setminus V} a_0(h,\varepsilon = 0) -
\int_{M' \setminus V'} a_0(h',\varepsilon = 0).
\endmultline  \tag{6.56}$$
Clearly, the definition is independent of the choice of $V$ and $V'.$

As an application of Theorem A and Theorem B we obtain the
following result:

\proclaim{Corollary C}{ Let $M$ and $M'$ be two closed manifolds with the same
fundamental group and the same dimension $d$ and let $\Cal W$ be an
$(\Cal A,\Gamma^{op})$-
Hilbert module of finite type. Suppose that $\tau = (h,g)$ and
$\tau' = (h',g')$ are generalized triangulations with $\# {\roman
{Cr}}_q(h) = \# {\roman {Cr}}_q(h') , \ (0 \leq q \leq d),$ and that
$(M,{\Cal W})$ and $(M',{\Cal W})$ are of determinant class.  Let
$T'_{{\roman {an}}} = T_{{\roman {an}}}(M',g',\Cal W) $ and
$T'_{{\roman {Re}}} = T_{{\roman {Re}}}(M',\tau ',\Cal W) .$   Then
$$\align
\log T_{{\roman {an}}} - \log T'_{{\roman {an}}} = & \log T_{{\roman
{Re}}}(\tau) - \log T_{{\roman {Re}}}(\tau').
\endalign$$
}
\endproclaim

Let $(M,\tau=(h,g))$ be a manifold equipped with a generalized triangulation.

Let $x_{q;j}\in {\roman {Cr}}_q(h)$ be a critical point of $h$ of index $q$ and
$U_{qj}$ a system of $H$-neighbourhoods of $x_{q;j}$ (cf Definition 5.1).
Introduce the manifolds
$$
M_I:=M\backslash\cup_{q,j}U'_{qj};\quad M_{II}:=\cup_{q,j}\overline{U_{qj}'},
$$
where $U'_{qj}$ is defined as in Definition 5.1. Both manifolds
$M_I$ and $M_{II}$ have the same boundary, given by a disjoint union of
spheres of dimension $d-1$.

Fix $\varepsilon>0$ and consider the operator $\Delta_q(t)+\varepsilon$. Its
symbol with respect to arbitrary coordinates $(\varphi,\psi)$ of
$(M,\Cal E\to M)$ is of the form
$$
a_2(x,\xi)+t^2||\nabla h||^2+a_1(x,\xi)+tL(x)+\epsilon\tag{6.57}
$$
where $a_i:B_{2\alpha}\times\bold R^d\to End(\Lambda^q(\bold R^d)
\otimes{\Cal W})$
$(i=1,2)$ are homogeneous of degree $i$ in $\xi$, where
$||\nabla h||^2:B_{2\alpha}\to\bold R$ is given by
$$
||\nabla h||^2=\sum_{1\leq i,j\leq d}g^{ij}\frac{\partial h}{\partial x_i}
\frac{\partial h}{\partial x_j}
$$
and where
$L:B_{2\alpha}\to End(\Lambda^q(\bold R^d))$ is the operator
$L=\Cal L_{\nabla h}+\Cal L_{\nabla h}^*$ of order 0 with $\Cal L_{\nabla h}$
denoting the Lie-derivative of $q$-forms along the vector field
$$
\nabla h=\sum_{i,j}g^{ij}\frac{\partial h}{\partial x_i}
\frac{\partial }{\partial x_j}.
$$
The operator $\Cal L_{\nabla h}^*$ is the adjoint of
$\Cal L_{\nabla h}$ with respect to the metric $g$ and is given by
$$
\Cal L_{\nabla h}^*=-(-1)^{q(d+q)}R_{d-q}\Cal L_{\nabla h}R_q\tag{6.58}
$$
where $R_q:\Lambda^q(B_{2\alpha})\to\Lambda^{d-q}(B_{2\alpha})$ is the Hodge
operator associated to the metric $\varphi^*g$. Recall that we have
denoted by ${\roman {Cr}}(h)$ the set of all critical points of $h$. Set
$M^*:=M\backslash {\roman {Cr}}(h)$. For an arbitrary chart $(\varphi,\psi)$ of
$(M^*,\Cal E|_{M^*}\to M^*)$, define,
as discussed in section 3, the symbol expansion
$\sum_{j\geq0}r_{-2-j}(h,\varepsilon,x,\xi,t,\mu)$ of the resolvent
$(\mu-\Delta_q(t)-\varepsilon)^{-1}$ inductively:
$$
r_{-2}(h,\varepsilon,x,\xi,t,\mu)=(\mu-a_2(x,\xi)-t^2||\nabla h||^2)^{-1}
$$
and, for $j\geq1$,
$$\align
r_{-2-j}& =
     (\mu-a_2-t^2||\nabla h||^2)^{-1}
     \sum\Sb1\leq|\alpha|\leq2\\l+|\alpha|=j\endSb
     \frac{1}{\alpha!}
     \partial_\xi^\alpha
     a_2\left( D_x \right)^\alpha r_{-2-l}\tag{6.59}\\
&+(\mu-a_2-t^2||\nabla h||^2)^{-1}
     \sum\Sb0\leq|\alpha|\leq1\\l+|\alpha|=j\endSb
    \partial_{\xi}^\alpha(a_1+tL)
    \left( D_x \right)^\alpha r_{-2-l}\\
&+(\mu-a_2-t^2||\nabla h||^2)^{-1}\varepsilon r_{-j}.
\endalign$$
Note that $r_{-2-j}$ has the following homogeneity property: for
$\lambda\in\bold R_+$
$$
r_{-2-j}(h,\varepsilon,x,\lambda\xi,\lambda t,\lambda^{1/2}\mu)=
\lambda^{-2-j}r_{-2-j}(h,\varepsilon,x,\xi,t,\mu).\tag{6.60}
$$
For later use, we introduce the densities $a_0(h,\varepsilon,x)$ on $M^*$
with values in $\Bbb R$, defined with respect to the chart
$(\varphi,\psi)$ and arbitrary $\varepsilon$ as
$$\align
a_0(h,\varepsilon,x)&=
   \left. \frac{\partial}{\partial s} \right| _{s=0}
   \left(\frac{1}{2\pi}\right)^d\frac{1}{2\pi i}
   \int_{\bold R^d}d\xi\tag{6.61}\\
   &\int_\Gamma d\mu
   \mu^{-s} {\roman {tr}}_Nr_{-2-d}(h,\varepsilon,x,\xi,t=1,\mu)\\
   &=\frac{-1}{(2\pi)^d}\int_{\bold R^d}d\xi\int_0^\infty d\mu
   {\roman {tr}}_N \left( r_{-2-d}(h,\varepsilon,x\xi,,t=1,-\mu) \right) .
\endalign$$

\proclaim{Proposition {6.6}} Assume that $(M^d,\tau=(h,g),\Cal W)$ and
$({M'}^d,\tau'=(h',g'),\Cal W)$ satisfy the hypothesis of Theorem B.
Then for any $\varepsilon>0$

(i) $\log det_N(\Delta_q(h,t)+\varepsilon)-\log det_N(\Delta_q(h',t)+
\varepsilon)$ has a
complete asymptotic expansion for $t\to\infty$ whose free term is
denoted by $\overline a_0:=\overline a_0(h, h',\varepsilon).$

(ii) The coefficient $\overline a_0$ can be represented in the form
$$
\overline a_0=\int_{M_{I}}a_0(h,\varepsilon,x)-
\int_{M'_{I}}a_0( h',\varepsilon',x')\tag{6.62}
$$
where $a_0(h,\varepsilon,x)$ and $a_0(h',\varepsilon,x')$ are the
densities introduced in \thetag {6.75} for arbitrary $\varepsilon$.

(iii) If $\dim M=d$ is odd then
$$
\overline a_0(h,h',\varepsilon)+\overline a_0(d-h,d-h',\varepsilon)=0\quad
(\text{ all }\varepsilon>0).\tag{6.63}
$$
\endproclaim

\demo{Proof} The proof is based on a Mayer-Vietoris type
formula (Theorem 3.6). Note that $\Delta_q(h,t)+\varepsilon$ is a family of
invertible, selfadjoint, elliptic operators with parameter $t$ of order 2
and weight 1 for any $\varepsilon>0$. The same is true for the operators
$(\Delta_q^I(h,t)+\varepsilon)_D$ and $(\Delta_q^{II}(h,t)+\varepsilon)_D$
obtained by
restricting $\Delta_q(h,t)+\varepsilon$ to $M_I$ and $M_{II},$ respectively,
and imposing
Dirichlet boundary conditions. We can therefore apply Theorem 3.6.
Denote by $R_{DN}(h,t,\varepsilon)$ the Dirichlet to Neumann operator defined
in section 3.3. We conclude from Theorem
3.6 (4) that $R_{DN}(h,t,\varepsilon)$ is an invertible pseudodifferential
operator with parameter of order 1 and weight 2 and from Theorem 3.6
(2) we conclude that $R_{DN}(h,t,\varepsilon)$ is elliptic with parameter $t$.
According to Theorem 3.4, $\log\det_N R_{DN}(h,t,\varepsilon)$  has an
asymptotic expansion for $t\to\infty$. Inspecting the principal symbol
of $(\Delta_q^I(h,t)+\varepsilon)_D$ one observes that $(\Delta_q^I(h,t)+
\varepsilon)_D$
is a family of invertible, selfadjoint differential operators with
parameter of order 2 and weight 1 which is elliptic with parameter.
 From Theorem 3.5 we therefore conclude that
$\log\det(\Delta_q^I(h,t)+\varepsilon)_D$ admits a complete asymptotic
expansion as $t\to\infty$. Finally $(\Delta_q^{II}(h,t)+\varepsilon)_D$ is a
family of invertible selfadjoint operators with parameter of order 2
and weight 1, which is, however, not elliptic with parameter.

Of course the same considerations can be made for the system
$( M',h',g')$ to conclude that
$\log\det_N R_{DN}( h',t,\varepsilon)$ and
$\log\det_N (\Delta_q^I(h',t)+\varepsilon)_D$ have both asymptotic
expansions for $t\to\infty$. Applying the Mayer-Vietoris type formula
(Theorem 3.6 (3)) for $\log det_N (\Delta_q(h,t)+\varepsilon)$ and
$\log det_N(\Delta_q(h',t)+\varepsilon)$ we obtain for the difference
$$\align
\log det_N & (\Delta_q(h,t)+\varepsilon)-\log det_N (\Delta_q(h',t)+
\varepsilon)\tag{6.64}\\
& =\log det_N(\Delta_q^I(h,t)+\varepsilon)_D-\log det_N(\Delta_q^I(h',t)+
\varepsilon)_D\\
& +\log det_N (\Delta_q^{II}(h,t)+\varepsilon)_D-\log det_N
(\Delta_q^{II}(h',t)+\varepsilon)_D\\
& +\log det_N R_{DN}(h,t,\varepsilon)-\log det_N R_{DN}(h',t,\varepsilon)\\
& +\log\bar c(h,t,\varepsilon)-\log\bar c(h',t,\varepsilon).
\endalign$$
Note that $M_{II}$ and $M'_{II}$ are isometric and
$\Cal E_{|M_{II}}$ as well as $\Cal E'_{|M'_{II}}$ are trivial.
Consequently
$$
\log det_N (\Delta_q^{II}(h,t)+\varepsilon)_D=\log det_N (\Delta_q^{II}(h',t)+
\varepsilon)_D.
$$
Due to our definition of $H$-coordinates the isometry between $M_{II}$ and
$\tilde M_{II}$ extends to neighbourhoods of $M_{II}$ and
$\tilde M_{II}$. As a consequence we conclude from
Theorem 3.6(4)and
Theorem 3.4 that  $\bar c(h,t,\varepsilon)=\bar c(h',t,\varepsilon)$ and
that $\log det_N R_{DN}(h,t,\varepsilon)$ and
\newline $\log det_N R_{DN}(h',t,\varepsilon)$
have identical asymptotic expansions.
We have therefore proved that
$$
\log det_N (\Delta_q(h,t)+\varepsilon)-\log det_N (\Delta_q( h',t)+
\varepsilon)
$$
has an asymptotic expansion as $t\to\infty$ which is
identical to the complete asymptotic expansion for
$\log\det_N(\Delta_q^I(h,t)+\varepsilon)_D
-\log\det_N(\Delta_q^I(h',t)+\varepsilon)_D$. According to Theorem 3.5 the
free term in the asymptotic expansions of both
$\log\det_N(\Delta_q^I(h,t)$ $+\varepsilon)_D$ and
$\log\det_N(\Delta_q^I(h',t)+\varepsilon)_D$ consists of a boundary
contribution and a contribution from the interior. Recall that
$\partial M_I$ and $\partial M'_I$ are isometric and that in
collar neighbourhoods of $\partial M_I$ and of $\partial M'_I$
the symbols of $(\Delta_q^I(h,t)+\varepsilon)_D$ and
$(\Delta_q^I(h',t)+\varepsilon)_D$ are identical when expressed in
(H)-coordinates. Therefore the boundary contributions are the same
and the free
term in the asymptotic expansion of
$\log det_N(\Delta_q^I(h,t)+\varepsilon)-\log det_N(\Delta_q^I(h',t)
+\varepsilon)$ is
given by
$$
\overline a_0=\int_{M_I}a_0(h,\varepsilon,x)-
\int_{M'_I}a_0(h',\varepsilon,x') \tag{6.65}
$$
where the densities $a_0(h,\varepsilon,x)$ and
$a_0(h',\varepsilon, x')$ are given by \thetag {6.61}.

Noting that $a_0(h,\varepsilon,x)$ and $a_0(h',\varepsilon,x')$ are
identical on
$M_{II}\backslash {\roman {Cr}}(h)\cong  M'_{II}\backslash
{\roman {Cr}} ( h')$
statement (ii) follows. Towards (iii), observe that if $M$ is of odd
dimension, the quantity
\newline $r_{-d-2}(h,\varepsilon,x,\xi,t,\mu)$ defining
$a_0(h,\varepsilon,x)$ satisfies, according to \thetag {6.59} and
\thetag {6.60} ,
$$
r_{-d-2}(d-h,\varepsilon,x,\xi,t,\mu)=r_{-d-2}(h,\varepsilon,x,\xi,-t,\mu)
\tag{6.66}
$$
and
$$
r_{-d-2}(h,\varepsilon,x,-\xi,-t,\mu)=-r_{-d-2}(h,\varepsilon,x,\xi,t,\mu).
\tag{6.67}
$$
Therefore
$r_{-d-2}(h,\varepsilon,x,\xi,t,\mu)+r_{-d-2}(d-h,\varepsilon,x,\xi,t,\mu)$
is an
odd function of $\xi$. Integrating over $|\xi|=1$ we conclude that
$a_0(h,\varepsilon,x)+a_0(d-h,\varepsilon,x)=0$.\qed
\enddemo
For any $\varepsilon > 0$
introduce the following perturbed version of $\log T(h,t)$
$$
A(h,t,\varepsilon):=\frac{1}{2}\sum_{q=0}^d(-1)^{q+1}q\log det_N (\Delta_q(h,t)
+\varepsilon).\tag{6.68}
$$
Note that $A(h,t,\varepsilon)$ can be written as a sum
$$
A(h,t,\varepsilon)=A_{{\roman {sm}}}(h,t,\varepsilon)+
A_{{\roman {la}}}(h,t,\varepsilon)\tag{6.69}
$$
where $A_{\roman{sm}}$ is defined in a fashion analogous to
 $\log T_{{\roman {sm}}}(h,t),$
$$
A_{{\roman {sm}}}(h,t,\varepsilon):=\frac{1}{2}
\sum_{q=0}^d(-1)^{q+1}q\log det_N(\Delta_q^{{\roman {sm}}}(h,t)+
\varepsilon)
$$
with
$$
\Delta_q^{sm}(h,t):=\Delta_q(h,t)|_{\Lambda^q(M;\Cal E)_{sm}}
$$
and $A_{la}(t,h,\varepsilon)$ is given by
$A(h,t,\varepsilon)-A_{sm}(h,t,\varepsilon)$. Observe
that the spectrum of the operator $\Delta_q^{sm}(h,t)$ tends to 0 as
$t\to\infty$ and therefore, by Theorem 5.7 (5)
$$
\log det_N (\Delta_q^{sm}(h,t)+\varepsilon)=m_q\log\varepsilon+
O\left(\frac{1}{\varepsilon}te^{-2t}\right)
$$
for $t\to\infty$. This shows that
$A_{sm}(h,t,\varepsilon)-A_{sm}(h',t,\varepsilon)$ is exponentially small as
$t\to\infty$ and hence, for any fixed $\varepsilon>0$, it has a trivial
 asymptotic expansion for $t\to\infty$. In view of
\thetag {6.69} and Proposition 6.6 we conclude that for any $\varepsilon>0,$
$A(h,t,\varepsilon)-A(h',t,\varepsilon)$ and
$A_{la}(h,t,\varepsilon)-A_{la}(h',t,\varepsilon)$ have complete
asymptotic expansions for $t\to\infty$ and, moreover, these
expansions are identical. In particular we conclude that the free
terms of the two expansions are identical
$$
FT(A_{la}(h,t,\varepsilon)-A_{la}(h',t,\varepsilon))
=FT(A(h,t,\varepsilon)-A(h',t,\varepsilon)).
$$
Using Proposition 6.6 (ii) and the fact that the densities
$a_0(h,\varepsilon,x)$ and $a_0(h',\varepsilon,x)$ (6.61) are
continuous in $\varepsilon$ we obtain

\proclaim{Lemma {6.7}}

(i) For any $\varepsilon>0$, $A_{la}(h,t,\varepsilon)-A_{la}(h',t,\varepsilon)$
has a
complete asymptotic expansion for $t\to\infty$ which is identical to
the asymptotic expansion for $A(h,t,\varepsilon)-A(h',t,\varepsilon)$.

(ii) The limit
$$
\lim_{\varepsilon\to0} \ FT(A_{la}(h,t,\varepsilon)-A_{la}(h',t,\varepsilon))
$$
exists and is given by
$$\align
\lim_{\varepsilon\to0}&
FT(A_{la}(h,t,\varepsilon)-A_{la}(h',t,\varepsilon))\tag{6.70}\\
&=\int_{M_{I}}a_0(h,\varepsilon=0,x)-
  \int_{M'_{I}}a_0(h',\varepsilon=0,x').
\endalign$$
\endproclaim

We proceed to investigate the left hand side of \thetag {6.70}.
For this we need the following estimate for the spectral
function $N_q(t,\lambda)$ of $\Delta_q(t)=\Delta_q(h,t)$.

\proclaim{Lemma {6.8}} There exists a constant $C>0$ independent
of $t$ and $\lambda \geq 1$ such that, for $t$ sufficiently large,
$$
N^{{\roman {la}}}_q(t,\lambda)\leq C\lambda^d.
$$
\endproclaim

\demo{Proof} First note that for $t$ sufficiently large,
$\Delta_q(t)\geq\Delta_q(0)=\Delta_q$ on $M_{I}$.
By Weyl's law, we conclude
$$
N_q^{I}(t,\lambda)\leq N_q^{I}(0,\lambda)\leq C_1\lambda^{d/2}
$$
where $N_q^{I}(t,\lambda)$ is the spectral function for
the operator $\Delta_q(t)$ restricted to $M_{I}$ with
Neumann boundary conditions (Neumann spectrum). Recall that
$M_{II}=\cup_{k,j}U_{kj}$. On each of the discs $U_{kj}$, $\Delta_q(t)$,
when expressed in (H)-coordinates, is the direct sum of shifted
harmonic oscillators of the form
$$
H_t:=-\frac{d^2}{dx^2}+t^2x^2+tc
$$
with $-\alpha<x<\alpha$ ($\alpha$ as in \thetag {6.61}). Following $[CFKS,
p.218]$ introduce the scaling operator $S_t$ defined by
$$
S_tf(x):=t^{1/2}f(tx).
$$
Then $S_{t^{1/2}} \cdot tK \cdot S_{t^{1/2}}^{-1}=H_t$ where
$$
K:=-\frac{d^2}{dx^2}+x^2+c.
$$
Therefore the Neumann spectrum  of $H_t$ on the interval $-\alpha<x<\alpha$
is the same as the Neumann spectrum of $tK$ when considered on the
interval $-\sqrt t\alpha<x<\sqrt t\alpha$. Denote by $N_{tK;\sqrt t}(\lambda)$
the spectral function of the operator $tK$ on the interval
$-\sqrt t\alpha<x<\sqrt t\alpha$ with Neumann boundary conditions
 and by $N_{tK;\sqrt t}^D(\lambda)$ the
spectral function of the operator $tK$ on the interval $-\sqrt t
\alpha<x<\sqrt t\alpha$ with Dirichlet spectrum.
 Note that for all $t\geq0$ and $\lambda$ sufficiently
large
$$
N_{tK;\sqrt t}(\lambda)\leq N_{tK;\sqrt t}^D(\lambda)+1\leq2
N_{K;\sqrt t}^D(\lambda/t).
$$
Comparing the Dirichlet problem for $K$ on $-\sqrt t\alpha\leq x\leq
\sqrt t\alpha$
with the one on the whole real line we conclude that $N_{K;\sqrt t}^D
(\lambda/t)\leq C_2\lambda/t\leq C_2\lambda$ for $t\geq1$ with a constant
$C_2>0$
independent of $\lambda$ and $t$. Hence we have shown that the spectral
function $N_q^{II}(t,\lambda)$ of the operator
$\Delta_q(t)$ on $M_{II}$ with Neumann boundary conditions can be estimated by
$$
N_q^{II}(t,\lambda)\leq C_3\lambda^d
$$
for a constant $C_3>0$ independent of $t$ and $\lambda$. The subadditive
property of the Neumann spectral function implies that
$$
N_q(t,\lambda+0)\leq N_q^I(t,\lambda+0)+N_q^{II}(t,\lambda+0)\leq C\lambda^d
$$
for some constant $C>0$ independent of $t$ and $\lambda$ and for $t$
sufficiently large.\qed
\enddemo
Let us introduce the trace of the heat kernel of
$\Delta^{{\roman {la}}}_q(t)$,
$$
\theta_q(t,\mu)=\int_1^\infty e^{-\mu\lambda }dN_q^{\roman {la}}(t,\lambda).
$$

\proclaim{Corollary {6.9}}

(i) There exists a  constant $C>0$ independent of $t$ and $\mu$ such
that, for $t$ sufficiently large and $\mu >0$
$$
\theta_q(t,\mu)\leq C\mu^{-d}.\tag{6.71}
$$

(ii) There exist constants $C>0$ and $\beta >0$ independent of $t$ and
$\mu$ such
that, for $t$ sufficiently large, and $\mu\geq1/\sqrt t$
$$
\theta_q(t,\mu)\leq Ce^{-\beta t\mu}.\tag{6.72}
$$
\endproclaim

\demo{Proof}

(i) By Proposition 5.2 there exists a constant $C_1>0$
such that $\roman {spec}(\Delta^{la}_q(t))\subset [C_1t,\infty)$ and
therefore
$$
\theta_q(t,\mu)=\int_{C_1t}^\infty e^{-\mu\lambda}dN_q(t,\lambda).
$$
Integrating by parts we obtain
$$
\theta_q(t,\mu)\leq\mu\int_{C_1t}^\infty e^{-\mu\lambda}N_q(t,\lambda)d
\lambda.\tag{6.73}
$$
By Lemma 6.8, one then concludes
$$
\theta_q(t,\mu)\leq\frac{C}{\mu^d}
\int_{C_1t\mu}^\infty e^{-\lambda}\lambda^dd\lambda
\leq \tilde C/\mu^d.
$$

(ii) From \thetag {6.72} and Lemma 6.8 we obtain
$$
\theta_q(t,\mu)\leq C\mu e^{-C_1t\mu/2}\int_{C_1t}^\infty e^{-\mu\lambda/2}
\lambda^dd\lambda\leq\frac{\tilde C}{\mu^d}e^{-C_1t\mu/2}.
$$
By choosing $\beta<C_1/2$ and $C>0$ sufficiently large we obtain
(ii).\qed
\enddemo

Recall from Theorem A that $\log T_{{\roman {la}}}(h,t)$  has an asymptotic
expansion for $t\to\infty$.

\proclaim{Proposition {6.10}}
$$
\lim_{\varepsilon\to 0}FT(A_{la}(h,t,\varepsilon)-A_{la}(h',t,\varepsilon))=
FT(\log T_{la}(h,t))-FT(\log T_{la}(h',t)).
$$
\endproclaim

\demo{Proof} We verify that the function, defined for $\varepsilon>0$ and $t$
sufficiently large by
$$
H(t,\varepsilon):=A_{la}(h,t,\varepsilon)-A_{la}(h',t,
\varepsilon)+\log T_{la}(h,t)
-\log T_{la}(h',t)
$$
is of the form
$$
H(t,\varepsilon)=\sum_{k=1}^d\varepsilon^kf_k(t)+g(t,\varepsilon)\tag{6.74}
$$
where $g(t,\varepsilon)=o(1)$ uniformly in $\varepsilon$. The statement of
the proposition can be deduced from \thetag {6.74} as follows: Recall
that for $\varepsilon>0$, $H(t,\varepsilon)$ has an asymptotic expansion for
$t\to\infty$. As $g(t,\varepsilon)=o(1)$ uniformly in $\varepsilon$ we
conclude that for any $\varepsilon>0$, $\sum_{k=1}^d\varepsilon^kf_k(t)$
has an asymptotic
expansion for $t\to\infty$. By taking $d$ different values
$0<\varepsilon_1<\dots<\varepsilon_d$ for $\varepsilon$ and using that
the Vandermonde
determinant is nonzero
$$
\det\pmatrix
\varepsilon_1 & \dots & \varepsilon_1^d\\
\vdots && \vdots\\
\varepsilon_d & \dots & \varepsilon_d^d
\endpmatrix
\neq 0
$$
we conclude that for any $1\leq k\leq d$, $f_k(t)$ has an asymptotic
expansion for $t\to\infty$ and that for any $\varepsilon>0$
$$
FT(H(t,\varepsilon))=\sum_{k=1}^d\varepsilon^kFT(f_k(t)).
$$
Hence $\lim_{\varepsilon\to 0}FT(H(t,\varepsilon))$ exists and
$\lim_{\varepsilon\to 0}FT(H(t,\varepsilon))=0.$
To prove \thetag {6.74}we introduce the
zeta function $\zeta_{q,la}$ of
$\Delta^{{\roman {la}}}_q(t)+\varepsilon$,
$$
\zeta_{q,la}(t,\varepsilon,s)=\frac{1}{\Gamma(s)}
\int_0^\infty\mu^{s-1}\theta_q(t,\mu)
e^{-\varepsilon\mu}d\mu\tag{6.75}
$$
with $\theta_q (t,\mu)$ given as above. The integral in
\thetag {6.75} can be split into two parts
$$
\zeta_{q,la}^I(t,\varepsilon,s)=
\frac{1}{\Gamma(s)}\int_{1/\sqrt t}^\infty\mu^{s-1}\theta_q(t,\mu)
e^{-\varepsilon\mu}d\mu\tag{6.76}
$$
and
$$
\zeta_{q,la}^{II}(t,\varepsilon,s)=
\frac{1}{\Gamma(s)}\int_0^{1/\sqrt t}\mu^{s-1}\theta_q(t,\mu)
e^{-\varepsilon\mu}d\mu.\tag{6.77}
$$
First let us consider
$$
\zeta_{q,la}^I(t,\varepsilon,s)-\zeta_{q,la}^I(t,\varepsilon=0,s)=
\frac{1}{\Gamma(s)}\int_{1/\sqrt t}^\infty\mu^s\theta_q(t,\mu)
\frac{e^{-\varepsilon\mu}-1}{\mu}d\mu.\tag{6.78}
$$
Note that
$$
\zeta_{q,la}^I(t,\varepsilon,s)-\zeta_{q,la}^I(t,\varepsilon=0,s)
$$
is, by Corollary 6.9 (ii), an entire function of $s$.
Noting that
$$
\left. \frac{d}{ds}\left(\frac{1}{\Gamma(s)}\right) \right| _{s=0} = 1
$$
and that $1-e^{-\varepsilon\mu}\leq\varepsilon\mu$, we obtain
$$\align
|\frac{d}{ds} _{s=0}&
(\zeta_{q,la}^I(t,\varepsilon,s)-\zeta_{q,la}^I(t,\varepsilon=0,s))|\\
& =|\int_{1/\sqrt t}^\infty\theta_q(t,\mu)
\frac{e^{-\varepsilon\mu}-1}{\mu}d\mu|\\
&\leq\varepsilon C\int_{1/\sqrt t}^\infty e^{-\beta t\mu}d\mu
=\frac{\varepsilon C}{\beta t}e^{-\beta\sqrt t}
\endalign$$

where we have used Corollary 6.9. To analyze the term
$$
\frac{d}{ds}_{s=0}(\zeta_{q,la}^{II}(t,\varepsilon,s)-\zeta_{q,la}
^{II}(t,\varepsilon=0,s)),
$$
first expand $(e^{-\varepsilon\mu}-1)/\mu$
$$
(e^{-\varepsilon\mu}-1)/\mu=\sum_{k=1}^d\frac{(-1)^k}{k!}\varepsilon^k\mu^{k-1}
+\varepsilon^{d+1}\mu^de(\varepsilon,\mu)
$$
where the error term is given by
$$
e(\varepsilon,\mu)=
\left(\sum_{k=d+1}^\infty\frac{(-1)^k}{k!}\varepsilon^k\mu^{k-1}\right)
/\varepsilon^{d+1}\mu^d.
$$
Note that, by Corollary 6.9,
 $\mu^d\theta_q(t,\mu)\leq C$ and
therefore
$$
 \int_0^{1/\sqrt t}\mu^s\theta_q(t,\mu)
 \varepsilon^{d+1}\mu^de(\varepsilon,\mu)d\mu
$$
is a meromorphic function of $s$,with $s=0$ a regular point and, for
$t$ sufficiently large
$$
\left. \left| \frac{d}{ds} \right| _{s=0} \left(\frac{1}{\Gamma(s)}
 \int_0^{1/\sqrt t}\mu^s\theta_q(t,\mu)
 \varepsilon^{d+1}\mu^de(\varepsilon,\mu)d\mu\right) \right|
 \leq\varepsilon^{d+1}C/\sqrt t
$$
where $C$ is independent of $t$ and $\varepsilon$, $0\leq\varepsilon\leq1$.
Finally,
recall that $\theta_q(t,\mu)$ admits an expansion for $\mu\to 0+$ of the
form
$$
\theta_q(t,\mu)=\sum_{j=0}^dC_j(t)\mu^{(j-d)/2}+\theta_q'(t,\mu)
$$
where $\theta_q'(t,\mu)$ is continuous in $\mu\geq 0$. Therefore, for
$1\leq k\leq d$,
$$
\frac{1}{\Gamma(s)}\int_0^{1/\sqrt t}\mu^s\theta_q(t,\mu)\frac{(-1)^k}{k!}
\varepsilon^k\mu^{k-1}d\mu
$$
is analytic in $s$ at $s=0$ and
$$
\sum_{k=1}^d \left. \frac{d}{ds} \right| _{s=0} \left( \frac{1}{\Gamma(s)}
\int_0^{1/\sqrt t}\mu^s\theta_q(t,\mu)\frac{(-1)^k}{k!}
\varepsilon^k\mu^{k-1}d\mu \right)
$$
is of the form $\sum_{k=1}^d\varepsilon^kf_k(t)$. This establishes
\thetag {6.74}.\qed
\enddemo

\demo{Proof of Theorem B} From Theorem A we know that
$\log T_{la}(h,t)-\log T_{la}(h',t)$ has an asymptotic expansion
for $t\to\infty$. By Proposition 6.10, the free term of the asymptotic
expansion $\overline a_0$ is given by
$$
\overline a_0=\lim_{\varepsilon\to0}FT(A_{la}(h,t,
\varepsilon)-A_{la}(h',t,\varepsilon)).
$$
By Lemma 6.7 (ii) we conclude that
$$
\overline a_0=\int_{M_{I}}a_0(h,\varepsilon=0,x)-
\int_{M'_{I}}a_0(h',\varepsilon=0,x').
$$
Equation \thetag {6.69} is proved in Proposition 6.5 (iii). In view of
the equality
$$
FT(\log T_{{\roman {an}}}(h,t)-\log T_{{\roman {sm}}}(h,t))
=FT(\log T_{{\roman {la}}}(h,t))
$$
one can see that $\overline a_0$ is independent of $h$ and
$ h'$ within the class of functions $h$  and $h'$
which give rise to the same cochain
complexes $C^*(M;\tau,\Cal W)$ respectively
$C^*(M;\tau',\Cal W')$.
This combined with the locality of $a_0$ implies that
$\int_{M_{I}}a_0(h,\varepsilon=0,x)=\int_{M_{I}}a_0(h',\varepsilon=0,x')$
which in turn implies (ii).
\qed
\enddemo

\demo{Proof of Corollary C} Choose a bijection
$\Theta: {\roman {Cr}}(h)\to {\roman {Cr}}(h')$
 so that $\Theta(x_{q;j})$ is a critical point
$x'_{q;j}$ of $h'$ of index $q$. By assumption $\Theta$
extends to an isometry
$\Theta:\cup_{q,j}U_{qj}\to\cup_{q,j} U'_{qj}$ where $(U_{qj})$ and
$(U'_{qj})$ are systems of H-neighbourhoods for $h,$ respectively
$h'$. Denote by $\tau,$ respectively, $\tau'$ the
triangulation induced by $(h,g),$ respectively $(h', g')$
and by $\tau_{\Cal D},$ respectively
$\tau'_{\Cal D}$ the triangulations $\tau_{\Cal D}=(d-h,g),$ respectively
$\tau'_{\Cal D}=(d-h',g')$. It follows from Poincar\'e duality and
 $d$ odd  that  $\log T_{\roman {met}}(\tau)=
\log T_{\roman {met}}(\tau_{\Cal D})$ and
$\log T_{\roman {comb}}(\tau)= \log T_{\roman {comb}}(\tau_{\Cal D})$ .

Using Theorem A for both $h$ and $d-h$, we obtain
$$\align
2\log T_{\roman {an}} - 2\log T'_{\roman {an}}
 & =FT(\log T_{\roman {an}}(h,t)-\log T_{\roman {an}}(h',t))\\
& +FT(\log T_{\roman {an}}(d-h,t)-\log T_{\roman {an}}(d-h',t))\\
& +2\log T_{\roman {met}}(\tau) - 2\log T_{\roman {met}}(\tau')
\endalign$$
Decomposing $\log T_{\roman {an}}(h,t)=\log T_{\roman {la}}(h,t) +
\log T_{\roman {sm}}(h,t)$ and
taking into account the asymptotics of $\log T_{\roman {sm}}(h,t)$
 (cf Theorem A) we conclude that
$$\align
2\log T_{\roman {an}} - 2\log T'_{\roman {an}}&
=2\log T_{\roman {comb}}(\tau) - 2\log T_{\roman {comb}}(\tau')\\
& +FT(\log T_{\roman {la}}(h,t)-\log T_{\roman {la}}(h',t))\\
& +FT(\log T_{\roman {la}}(d-h,t)-\log T_{\roman {la}}(d-h',t))
\endalign$$
 from which the Corollary follows by \thetag {6.55} and
\thetag {6.56}.\qed
\enddemo

\subheading{ 6.3 Proof of Theorem 2}

In this subsection we provide the proof of Theorem 2 using Corollary C
of subsection 6.2 together with the product formulas for the
Reidemeister torsion and the analytic torsion established in
Proposition 4.1.

First we need the following result concerning the metric anomaly of
the analytic torsion, which is a straightforward generalization of a
classical result due to Ray-Singer \cite{RS} , and can be proved by arguments
similar to the ones presented in subsection 6.1.

\proclaim{Lemma 6.11}{Let $M^d$ be a closed manifold of odd dimension
$d$ such that $(M,{\Cal W})$ is of determinant class.  Assume that
$g(u)$ is a smooth one-parameter family of Riemannian metrics on $M.$
Then $\log T_{{\roman {an}}}(M, g(u), {\Cal W})$ is a smooth function
of $u$ whose derivative is given by
$$
\frac{d}{du}\log T_{{\roman {an}}}(M, g(u), {\Cal W}) = \frac{1}{2}
\frac{d}{du} \sum_{q=0}^d (-1)^q \ldetn (\sigma_q(u)^* \sigma_q(u)) \tag{6.79}
$$
where $\sigma_q(u)$ is the ${\Cal A}$-linear, bounded isomorphism
$$\sigma_q(u): {\roman {Null}}\Delta_q(u_0) \longrightarrow  {\roman
{Null}}\Delta_q(u) $$provided by Hodge theory and $u_0$ is arbitrary
but fixed.}
\endproclaim

Given generalized triangulations $\tau = (h,g')$ and $\tau' = (h',g'')$
of $M, \  \tau'$ is called a subdivision of $\tau$ if
\roster

\item ${\roman {Cr}}_q(h) \subset {\roman {Cr}}_q(h') \ \ (0 \leq q
\leq d)  $

\item $W_x^\pm(h',g'') \subset W_x^\pm(h,g')$ for any $x \in {\roman
{Cr}}_q(h). $
\endroster

The following result can be found in \cite{Mi2}:

\proclaim{Lemma 6.12}{Let $\tau =(h,g')$ be a generalized
triangulation, $0 \leq q \leq d-1$ an integer and $x, y$ two distinct points
in $M \setminus {\roman {Cr}}(h).$  Then there exists a generalized
triangulation $\tau' = (h',g'') $ with the following properties:

\roster

\item ${\roman {Cr}}_k(h')  = {\roman {Cr}}_k(h)$ for $k \neq q, \
q+1;$

\item ${\roman {Cr}}_q(h') = {\roman {Cr}}_q(h) \cup \{x\}; \ {\roman
{Cr}}_{q+1}(h') = {\roman {Cr}}_{q+1}(h) \cup \{y\};$

\item $\tau'$ is a subdivision of $\tau;$

\item $W_y^- \cap W_x^+$ is connected.

\endroster
}
\endproclaim

Since the Reidemeister torsion does not change under subdivision (cf
\cite{Mi1}), one obtains $$T_{{\roman {Re}}}(M,g,{\Cal W},\tau)
=T_{{\roman {Re}}}(M,g,{\Cal W},\tau').$$

{\sl Proof of Theorem 2} \
By Lemma 6.12, which concerns the metric anomaly, and in view of the
definition of $T_{{\roman met}},$  it suffices to prove Theorem 1 in
the case where $g=g', \ \tau = (g',h).$

Consider the sphere $S^6 = \{x = (x_1, \dots, x_7) \in {\Bbb R}^7;
\sum x_j^2 = 1\}$ with an arbitrary generalized triangulation $\tau_1 =
(h_1,g_1).$  Let $\tau= (h,g)$ be a generalized triangulation for $M$
and consider $M\times S^6,$ endowed with the Riemannian metric $g
\times g_1$ and the triangulation $\tau \times \tau_1 = (h+h_1, g
\times g_1).$  Note that $\Gamma = \pi_1(M) = \pi_1(M\times S^6).$
By assumption, $(M,{\Cal W})$ is of determinant class.  As $S^6$ is of
determinant class, $(M\times S^6,{\Cal W})$ is of determinant class as
well.  Moreover, by the product formulas of Proposisiton 4.1,
$$\align
\log T_{{\roman {an}}} (M\times S^6, g \times g_1, {\Cal W}) = & 2\log
T_{{\roman {an}}} (M, g , {\Cal W})  \tag{6.80}
\endalign$$
and
$$\align
\log T_{{\roman {Re}}} (M\times S^6, g \times g_1, {\Cal W},\tau \times
\tau_1) = & 2\log T_{{\roman {Re}}} (M, g , {\Cal W},\tau)  \tag{6.81}
\endalign$$
where we used that $\chi(S^6) = 2$ and that $\chi(M,{\Cal W}) = 0$ (as
$M$ is of odd dimension).

Next, consider the product $S^3 \times S^3$ of the 3-spheres, $S^3 =  \{x =
(x_1, \dots, x_4) \in {\Bbb R}^4; \sum x_j^2 = 1\},$ with an arbitrary
generalized triangulation $\tau_2 = (h_2,g_2).$  Arguing as above, we
conclude that $(M\times S^3 \times S^3,{\Cal W})$ is of determinant class and
that, by the product formulas of Proposisiton 4.1,
$$\align
\log T_{{\roman {an}}} (M\times S^3 \times S^3, g \times g_2, {\Cal W})
= & 0
\endalign$$
and
$$\align
\log T_{{\roman {Re}}} (M\times S^3 \times S^3, g \times g_2, {\Cal
W},\tau \times \tau_2) = & 0
\endalign$$
where we used that $\chi(S^3 \times S^3) = 0$ and that $\chi(M,{\Cal W})
= 0.$

Choose a subdivision $\tau' = (h',g'')$ of the generalized
triangulation $\tau \times \tau_1$ in $M \times S^6$ and a subdivision
$\tau'' = (h'',g'')$ of the generalized triangulation $\tau \times
\tau_2$ in $M \times S^3 \times S^3$ so that, for $0 \leq  q \leq 6, \
\# {\roman {Cr}}_q(h') =   \# {\roman {Cr}}_q(h'').$  This is possible
because $M\times S^6$ and $M \times S^3 \times S^3$ are both of odd
dimension and therfore $\chi(M \times S^3 \times S^3,{\Cal W}) = \chi(M
\times S^6,{\Cal W}) = 0.$  We conclude from the above, Corollary C,
Lemma 6.11 and Lemma 6.12 that
$$\multline
2 \log T_{{\roman {an}}} (M, g , {\Cal W}) - 2\log T_{{\roman {Re}}}
(M,  g , {\Cal W}, \tau)\\   =  \log T_{{\roman {an}}} (M\times S^6,
g \times g_1, {\Cal W})  - \log T_{{\roman {Re}}} (M\times S^6,
g \times g_1, {\Cal W}, \tau \times \tau_1)  \\
 =  \log T_{{\roman {an}}} (M\times S^6, g' , {\Cal W}) - \log
T_{{\roman {Re}}} (M \times S^6,  g' , {\Cal W}, \tau') \\
 =  \log T_{{\roman {an}}} (M\times S^3 \times S^3, g'' , {\Cal W}) -
\log T_{{\roman {Re}}} (M \times S^3 \times S^3,  g'' , {\Cal W},
\tau'') \\
 =  \log T_{{\roman {an}}} (M\times S^3 \times S^3, g\times g_2 ,
{\Cal W}) - \log T_{{\roman {Re}}} (M \times S^3 \times S^3,  g\times
g_2 , {\Cal W}, \tau \times \tau_2)
 =  0.
\endmultline$$
This proves Theorem 2.
\vfill \eject

\end

\enddocument